\documentclass[12pt,twoside]{amsart}
\usepackage{amssymb}
\usepackage{amscd}
\usepackage[abbrev,alphabetic]{amsrefs}
\usepackage{hyperref}
\usepackage{comment}
\usepackage{array,multirow,tabularx,longtable}
\usepackage{tikz}
\usetikzlibrary{cd}
\usepackage{here}
\usepackage{multirow}
\usepackage[margin=1.25in]{geometry}

\usepackage{fancybox}
\usepackage{ascmac}

\title[Elliptic singularities and threefold flops]
{Elliptic singularities and threefold flops in positive characteristic} 
\author{Hiromu Tanaka} 
\address{Graduate School of Mathematical Sciences, 
The University of Tokyo, 
3-8-1 Komaba, Meguro-ku, Tokyo 153-8914, JAPAN} 
\email{tanaka@ms.u-tokyo.ac.jp}
\newcommand{\aff}[0]{{\operatorname{aff}}}

\newcommand{\reg}[0]{{\operatorname{reg}}}

\newcommand{\red}[0]{{\operatorname{red}}}

\renewcommand{\Im}[0]{{\operatorname{Im}}}
\newcommand{\Tr}[0]{{\operatorname{Tr}}}
\newcommand{\Proj}[0]{{\operatorname{Proj}}}
\newcommand{\Spec}[0]{{\operatorname{Spec}}}

\newcommand{\Bs}[0]{{\operatorname{Bs}}}
\newcommand{\Supp}[0]{{\operatorname{Supp}}}
\newcommand{\Pic}[0]{{\operatorname{Pic}}}
\newcommand{\mult}[0]{{\operatorname{mult}}}
\newcommand{\Ex}[0]{{\operatorname{Ex}}}

\newcommand{\univ}[0]{{\operatorname{univ}}}
\newcommand{\gen}[0]{{\operatorname{gen}}}

\renewcommand{\min}[0]{{\operatorname{min}}}
\newtheorem{thm}{Theorem}[section]
\newtheorem{lem}[thm]{Lemma}
\newtheorem{cor}[thm]{Corollary}
\newtheorem{prop}[thm]{Proposition}
\newtheorem{claim}[thm]{Claim}    
    
\newtheorem*{claim*}{Claim}         
\newtheorem{step}{Step}

\theoremstyle{definition}
\newtheorem{ex}[thm]{Example}

\newtheorem{dfn}[thm]{Definition}

\newtheorem{rem}[thm]{Remark}
      
\newtheorem{nota}[thm]{Notation}         
\newtheorem{nasi}[thm]{}

\makeatletter
  
  \@addtoreset{equation}{thm}
  \makeatother

\newcommand{\cred}{\color{black}}
\newcommand{\cyan}{\color{cyan}}

\newcommand{\MO}{\mathcal{O}}

\newcommand{\Q}{\mathbb{Q}}
\newcommand{\Z}{\mathbb{Z}}
\newcommand{\F}{\mathbb{F}}
\newcommand{\A}{\mathbb{A}}

\renewcommand{\P}{\mathbb{P}}

\newcommand{\m}{\mathfrak{m}}

\begin{document}

\maketitle

\begin{abstract}
Let $X$ be a smooth  threefold  over an algebraically closed field of positive characteristic. 
We prove that an arbitrary flop of $X$ is smooth. To this end, we study Gorenstein curves of genus one and two-dimensional elliptic singularities defined over imperfect fields. 
\end{abstract}

\tableofcontents

\section{Introduction}

Flop is an important operation in birational geometry. 
The purpose of this paper is to establish the following theorem, 
which is a positive-characteristic analogue of \cite[Theorem 2.4]{Kol89}.

\begin{thm}[Theorem \ref{t-flop}]\label{main-flop}
Let $k$ be an algebraically closed field of characteristic $p>0$ and 
{\cred let} $Y$ be a smooth  threefold over $k$. 
Let $\psi : Y \to Z$ be a flopping contraction, 
i.e., 
$\psi$ is a projective birational $k$-morphism to a normal threefold $Z$ over $k$ such that 
$\dim \Ex(\psi)=1$ and $K_Y$ is $\psi$-numerically trivial. 
Then the flop $Y^+$ of $\psi$ exists and $Y^+$ is smooth. 
\end{thm}

By applying the same argument as in characteristic zero, 
the problem is reduced to prove the following result on Gorenstein terminal threefolds.

\begin{thm}[Theorem \ref{t-terminal3}]\label{main-terminal}
Let $k$ be an algebraically closed field of characteristic $p>0$ and 
let $Z$ be a Gorenstein terminal threefold over $k$. 
Let $P \in Z$ be a singular point. 
Then $\dim_k \m_P/\m_P^2 = 4$ and $\mult\,\MO_{Z, P} =2$, 
where $\m_P$ denotes the maximal ideal of the local ring $\MO_{Z, P}$. 
\end{thm}

In order to establish Theorem \ref{main-terminal}, 
we again follow the strategy of the proof  of characteristic zero. 
In characteristic zero, the proof of Theorem \ref{main-terminal} is carried out by taking a general hyperplane section passing through $P$. 
Although the naive analogue does not work by the failure of the Bertini theorem, 
we shall use the generic hyperplane section as a replacement, which is a hyperplane section $X$ of $Z \times_k \kappa$ for the function field $\kappa$ of $\P^n_k$ for some $n$ (cf. Subsection \ref{ss-generic}). 
Then the essentially  new part of this paper is the following theorem on two-dimensional elliptic singularities over an arbitrary field of characteristic $p>0$.

\begin{thm}[Theorem \ref{t-ell-3}, Theorem \ref{t-ell-12}]\label{main-ell}
Let $\kappa$ be a field of characteristic $p>0$ and 
let $X$ be a normal surface over $\kappa$ with a unique non-regular point $P$. 
Let $f: Y \to X$ be the minimal resolution of $X$. 
Assume that $P \in X$ is an elliptic singularity, i.e., $X$ is Gorenstein and $\dim_{\kappa} R^1f_*\MO_Y =1$. 
Set $Z$ to be the nonzero effective Cartier divisor on $Y$ satisfying $K_Y + Z \sim f^*K_X$. 
Then the following hold. 
\begin{enumerate}
    \item If $-Z^2 \geq 3$, then the blowup $\overline Y \to X$ at $P$ coincides with the canonical model of $Y$ over $X$. 
    \item If $-Z^2 =1$ or $-Z^2 =2$, then $\dim_{\kappa} (\m_P/\m_P^2) =3$ and $\mult\,\MO_{X, P} =2$, where $\m_P$ denotes the maximal ideal of the local ring $\MO_{X, P}$.
\end{enumerate}
\end{thm}

The main result (Theorem \ref{main-flop}) will be applied in the author's papers on classification of Fano threefolds in positive characteristic \cite{TanII}, \cite{TanIV}.  

\subsection{Overviews of proofs}

As explained above, the essentially new part is the result on elliptic singularities over 
an arbitrary field $\kappa$ (Theorem \ref{main-ell}). 
Most of Section \ref{s-genus-one-prime} and Section \ref{s-ell} are devoted to proving this theorem. 
We here overview of how to show Theorem \ref{main-ell}(1). 
The main part is to prove the very ampleness and the projective normality 
for a divisor on an elliptic-like curve $Z$. 


Although $Z$ is not necessarily integral, 
let us start with the case when $Z$ is integral. 

\begin{thm}[Theorem \ref{t-g1-va-prime}]\label{intro-t-g1-va-prime}
Let $\kappa$ be a field of characteristic $p>0$. 
Let $Z$ be a one-dimensional Gorenstein projective integral scheme over $\kappa$ 
with $\dim_{\kappa} H^1(Z, \MO_Z) =1$. 
Take a Cartier divisor $D$ on $Z$ such that $\deg_{\kappa} D \geq 3$. 
Then $\bigoplus_{m=0}^{\infty} H^0(Z, \MO_Z(mD))$ is generated by 
$H^0(Z, \MO_Z(D))$ as a $\kappa$-algebra.  
\end{thm}

The main steps of the proof of {\cred Theorem \ref{intro-t-g1-va-prime}} are as follows. 
\begin{enumerate}
\item[(I)] $|D|$ is base point free. 
For its induced morphism $\varphi : Z \to \P^N_{\kappa}$, 
let $\psi : Z \to W$ be the induced morphism onto its image $W := \varphi(Z) \subset \P^N_{\kappa}$. 
\item[(II)] $\psi : Z \to W$ is birational. 
\item[(III)] $\psi : Z \to W$ is an isomorphism, i.e., $|D|$ is very ample. 
\end{enumerate}

(I) We first treat the case when $Z$ is regular. 
By the Riemann-Roch theorem, we may assume that $D$ is effective. 
If $D$ is not equal to a point, then the problem is easy, because we obtain an exact sequence 
\[
H^0(Z, \MO_Z(D)) \to H^0(P, \MO_Z(D)|_P) \to H^1(Z, \MO_Z(D-P))=0
\]
for $D = P+P'+ \cdots$. 
Hence the problem is reduced to the case when $D$ is a closed point of $Z$. 
In this case, we obtain an exact sequence: 
\[
H^0(Z, \MO_Z(D)) \xrightarrow{\rho} H^0(D, \MO_Z(D)|_D) \to H^1(Z, \MO_Z) \to H^1(Z, \MO_Z(D))=0. 
\]
Note that $\dim_{\kappa} H^0(D, \MO_Z(D)|_D) = \dim_{\kappa} H^0(D, \MO_D) = \deg_{\kappa} D \geq 3$. 
It suffices to show that $\rho$ is nonzero, 
which follows from 
\[
\dim_{\kappa} H^0(D, \MO_Z(D)|_D) \geq 3 > 1 = \dim_{\kappa} H^1(Z, \MO_Z). 
\]
Thus $|D|$ is base point free if $Z$ is regular. 

We now treat the case when $Z$ is not regular. 
By the same argument as above, 
we see that the base locus ${\rm Bs}\,|D|$, which is a closed subset of $Z$, 
is disjoint from the regular locus $Z_{\reg}$ of $Z$. 
Let $\nu : Z' \to Z$ be the normalisation of $Z$. 
As a typical case, let us assume that $Z' \simeq \P^1_{\kappa'}$ for a field $\kappa'$ and 
that 
$\nu : Z' \to Z$ coincides with the pinching along 
the field extension $\kappa \hookrightarrow \kappa'$ of finite degree. 
In this case, $P := Z \setminus Z_{\reg}$ is a closed point, 
and ${\rm Bs}\,|D| \subset \{ P\}$. 
Take an open cover 
\[
Z' = \P^1_{\kappa'} = \Spec\,\kappa'[x] \cup \Spec\,\kappa'[y]
\]
with the gluing relation $xy=1$. 
We may assume that $P':= \nu^{-1}(P)_{\red}$ is the origin of $\Spec\,\kappa'[x]$. 
In this case, we obtain 
\[
Z = \Spec\,A \cup \Spec\,\kappa'[y]\qquad \text{for}\qquad    A := \kappa +x\kappa'[x]. 
\]
We may assume that $D|_{\Spec\,A} = \Spec\,(A/fA)$ for some $f \in A$ and that $f\kappa'[x] = x\kappa'[x]$, as otherwise the problem is easier. 
Then we can check that ${\rm div}_{Z'}(f)= P' -Q'$ for the origin $Q'$ of the other open subset, 
where ${\rm div}_{Z'}(f)$ denotes te principal divisor on $Z'$ defined by $f \in \kappa'[x] \subset K(Z')$. 
In particular, $f^{-1} \in \kappa'[y]$. 
Set $g := 1+f \in A$. 
Then 
\[
\frac{g}{f} = \frac{1+f}{f} = 1 + f^{-1} \in \kappa'[y]. 
\]
Therefore, we obtain $\frac{g}{f} \in H^0(X, \MO_X(D))$ as follows: 
\begin{itemize}
    \item $\frac{g}{f} \in f^{-1}A =\MO_X(D)|_{\Spec\,A}$, and 
    \item $\frac{g}{f} \in \overline{A} =\MO_X(D)|_{\Spec\,\kappa'[y]}$. 
\end{itemize}
Since $\frac{g}{f}$ does not vanish at the origin $P$, 
$|D|$ is base point free.

(II) 
By (I), we have the finite surjective morphism $\psi : Z \to W$ to a projective curve $W$ 
induced by an ample globally generated divisor $D$. 
In particular, we can write $D \sim \psi^*D_W$ for some ample Cartier divisor $D_W$ {\cred on $W$}. 
The essential part is to apply the following inequality for the $\Delta$-genus introduced by Fujita 
\cite{Fuj90}: 
\[
 \Delta(W, D_W) := \dim W +\deg D_W - h^0(W, D_W) \geq 0.
\] 
 Fujita established this inequality for the case when the base field $\kappa$ is algebraically closed. 
Although the author does not know whether this inequality extends to the general case, 
we can prove this for the case when $\dim W =1$ (Lemma \ref{l-Delta-nonneg}). 
This, together with standard argument, implies the following: 
\begin{eqnarray*}
0 & \leq & \Delta(W, D_W)\\
&=& 1 + \deg D_W - h^0(W, D_W)\\
&=& 1 + \frac{\deg D}{\deg \psi } - h^0(Z, D)\\
&=& 1 + \frac{\deg D}{\deg \psi } - \deg D\\
&=& 1 + \deg D \left( \frac{1}{\deg \psi } - 1\right)\\
&\leq & 1 + 3 \left( \frac{1}{\deg \psi } - 1\right).
\end{eqnarray*}
Hence $\deg \psi \leq \frac{3}{2}$, i.e., $\deg \psi =1$. 

(III) 
By (I) and (II), we have the induced morphisms: 
\[
\varphi : Z \xrightarrow{\psi} W \hookrightarrow \P^N_{\kappa}, 
\]
where $\varphi$ is the morphism induced by the complete linear system, $W := \varphi(Z)$, and 
$\psi : Z \to W$ is a finite birational morphism. 
One of the technical parts is to prove $H^1(W, \MO_W(1))=0$ for $\MO_W(1) := \MO_{\P^N_k}(1)|_W$. 
We here only give a sketch of the proof of  (III) by assuming $H^1(W, \MO_W(1))=0$. 
We  get the following for $\MO_Z(1) := \psi^*\MO_W(1)$: 
\begin{equation}\label{intro-e01-g1-va1}
\chi(W, \MO_W(1)) = h^0(W, \MO_W(1)) = h^0(Z, \MO_Z(1)) = \chi(Z, \MO_Z(1)). 
\end{equation}
We have the following conductor exact sequence 
for the conductor closed subschemes $C_Z \subset Z$ and $C_W \subset W$
(\ref{e3-cond}): 
\begin{equation}\label{intro-e02-g1-va1}
0 \to \MO_W \to \psi_* \MO_Z \oplus \MO_{C_W} \to \psi_*\MO_{C_Z} \to 0. 
\end{equation}
We then obtain 
\[
\chi(W, \MO_W) - \chi(Z, \MO_Z) = \chi(C_W, \MO_{C_W}) - \chi(C_Z, \MO_{C_Z}) 
\overset{{\rm (i)}}{=} \chi(W, \MO_W(1)) -\chi(Z, \MO_Z(1))\overset{{\rm (ii)}}{=} 0, 
\]
where (ii) follows from (\ref{intro-e01-g1-va1}) and (i) is obtained by applying the  tensor 
product $(-) \otimes \MO_W(1)$ to the exact sequence (\ref{intro-e02-g1-va1}). 
By $H^0(Z, \MO_Z)=H^0(W, \MO_W) =\kappa$, we get $h^1(Z, \MO_Z)=h^1(W, \MO_W)$. 
Then (\ref{intro-e02-g1-va1}) induces the following exact sequence: 
\[
0 \to  H^0(W, \MO_W) \to H^0(Z, \MO_Z) \oplus H^0(C_W, \MO_{C_W}) \to H^0(C_Z, \MO_{C_Z}) \to 0.
\]
Again by $H^0(Z, \MO_Z)=H^0(W, \MO_W) =\kappa$, 
we get the induced $\kappa$-linear isomorphism: 
\[
 H^0(C_W, \MO_{C_W}) \xrightarrow{\simeq} H^0(C_Z, \MO_{C_Z}), 
\]
which implies $\MO_{C_W} \xrightarrow{\simeq} \psi_*\MO_{C_Z}$. 
Again by the conductor exact sequence (\ref{intro-e02-g1-va1}), 
we obtain $\MO_W \xrightarrow{\simeq} \psi_*\MO_Z$, i.e., $\psi: Z \to W$ is an isomorphism. 
For more details on the proof of Theorem \ref{intro-t-g1-va-prime}, 
see Section \ref{s-genus-one-prime}. 

\medskip

Based on Theorem \ref{intro-t-g1-va-prime} explained above, 
we now overview of  how to prove Theorem \ref{main-ell}(1). 
The key result is the following theorem, which is a variant of 
Theorem \ref{intro-t-g1-va-prime}.

\begin{thm}[Theorem \ref{t-g1-va}]\label{intro-t-g1-va}
Let $\kappa$ be a field and let $P \in X$ be an elliptic singularity over $\kappa$. 
Assume that $X$ is affine, $K_X\sim 0$, 
and $X \setminus P$ is regular. 
Let $f : Y \to X$ be the minimal resolution of $X$. 
We define the effective Cartier divisor $Z$ by $K_Y +Z \sim f^*K_X$. 
Set $L_Z :=\MO_Y(-Z)$. 
If $-Z^2 \geq 3$, then  $\bigoplus_{m=0}^{\infty} H^0(Z, L_Z^{\otimes m})$  is generated by 
$H^0(Z, L_Z)$  as a $k$-algebra. 
\end{thm}

By contracting all the curves $C$ satisfying $f(C) =P$ and $K_Y \cdot C =0$, 
we obtain the canonical model $\overline Y$ of $Y$ over $X$: 
\[
f : Y \xrightarrow{g} \overline Y \xrightarrow{f} X, 
\]
where $\overline Y$ has at worst canonical singularities. 
In particular, $K_{\overline Y}$ is $\overline{f}$-ample and $K_Y \sim g^*K_{\overline Y}$. 
For $\overline Z := g_*Z$, we see that $K_{\overline Y} + \overline Z \sim \overline{f}^*K_X$, and hence $\overline Z$ is an effective Cartier divisor. 
Set $L_{\overline Z} := \MO_{\overline Y}(-\overline Z)$. We get $L_Z \simeq g^*L_{\overline Z}$. 
Note that $L_Z$ is nef, whilst $L_{\overline Z}$ is ample. 
This is the main advantage of the canonical model $\overline Y$. 
In what follows, we assume  $Y = \overline Y$ for simplicity. 
Similar to the above, 
the main steps 
are as follows. 
\begin{enumerate}
\item[(I)'] $|L_Z|$ is base point free. 
For its induced morphism $\varphi : Z \to \P^N_{\kappa}$, 
let $\psi : Z \to W$ be the induced morphism onto its image $W := \varphi(Z) \subset \P^N_{\kappa}$. 
\item[(II)'] $\psi : Z \to W$ is birational. 
\item[(III)'] $\psi : Z \to W$ is an isomorphism, i.e., $|L_Z|$ is very ample. 
\end{enumerate}
The proof of (III)' is identical to that of (III) above. 
In what follows, we overview some of the ideas for (I)' and (II)'. 

(I)' If $Z$ is an integral scheme, then we are done by (I). 
We may assume that $Z$ is not a prime divisor on $Y$. 
A key property is as follows: 
\begin{enumerate}
\item[$(*)$] 
Take a prime divisor $C \subset \Supp\,Z$ and set $Z' := Z -C \neq 0$. 
Then  the restriction map 
\[
H^0(Z, L_Z) \to H^0(Z', L_Z|_{Z'})
\]
is surjective (Proposition \ref{p-bpf-smaller}) and $|L_Z|_{Z'}|$ is very ample 
(Proposition \ref{p-g0-bpf}(2), Lemma \ref{l-ell-fund}(2)). 
\end{enumerate}
By $(*)$, $\Bs\,|L_Z|$ is disjoint from $Z'$ for every prime divisor $C \subset \Supp\,Z$, 
which implies $\Bs\,|L_Z| = \emptyset$. 

(II)' 
We only  treat the following two typical cases. 
\begin{enumerate}
\item[(a)] $Z = C_1 + C_2 + C_3$, where $C_1, C_2, C_3$ are mutually distinct prime divisors. 
\item[(b)] $Z = n C$, where $n \geq 3$ and $C$ is a prime divisor. 
\end{enumerate}

(a) Assume $Z = C_1 + C_2 + C_3$, where $C_1, C_2, C_3$ are mutually distinct prime divisors. 
Applying $(*)$ by setting $C := C_1$ and $Z' := Z - C_1 = 
{\cred C_2 + C_3}$, 
we see that $C_2$ and $C_3$ are mapping to distinct irreducible components $W_2$ and $W_3$ of $W$. 
By symmetry, $W$ has exactly three irreducible components: $W = W_1 \cup W_2 \cup W_3$ such that $\psi(Z_i) = W_i$ for each $i \in \{1, 2, 3\}$. 
Moreover, $\psi|_{Z_1 \cup Z_2}:Z_1 \cup Z_2 \to W_1 \cup W_2$ is an isomorphism, 
and hence $\psi : Z \to W$ is birational again by symmetry. 

(b) Assume $Z = n C$, where $n \geq 3$ and $C$ is a prime divisor. 
Set $Z' :=(n-1)C$ and let $W'$ be the scheme-theoretic image of $Z'$ by $\psi : Z \to W$: 
\[
\begin{tikzcd}[column sep=huge, row sep=large]
Y & Z \arrow[l, hook'] 
\arrow[d, "\psi"'] & Z' \arrow[l, hook', "j_Z"'] 
\arrow[d, "\psi'"',"\simeq"] \arrow[ll, bend right, hook']\\
& W& W'. \arrow[l, hook', "j_W"']
\end{tikzcd}
\]
By $(*)$, 
the induced morphism $\psi' : Z' \to W'$ is an isomorphism. 
Let $\eta_Z$ and $\eta_W$ be the generic points of 
$Z$ and $W$, respectively. 
Consider the following local rings at the generic points 
\[
A := \MO_{W, \eta_W},\quad 
B := \MO_{Z, \eta_Z},\quad 
R := \MO_{Y, \eta_Z}, 
\]
\[
A' := \MO_{W', \eta_{W}},\qquad 
B' := \MO_{Z', \eta_{Z}}. 
\]
Note that $R$ is a discrete valuation ring. 
Pick a generator $t \in R$ of the maximal ideal $\m_R$ of $R$, i.e., 
$\m_R = tR$. 
In particular, $B = R/t^nR$ and $B' = R/ t^{n-1}R$. 
We have the following commutative diagram consisting of the induced ring homomorphisms: 
\[
\begin{tikzcd}[column sep=huge]
R \arrow[r, twoheadrightarrow, "\pi"] \arrow[rr, twoheadrightarrow, bend left, "\pi'"] & B= R/t^nR 
\arrow[r, twoheadrightarrow, "\pi''"] & B' = R/ t^{n-1}R\\
& A \arrow[r, twoheadrightarrow] \arrow[u, hook] & A'.\arrow[u, "\simeq"'] \\
\end{tikzcd}
\]
In order to show that $\psi :  Z \to W$ is birational, it suffices to show $A=B$. 
In what follows, we identify $A'$ and $B'$ via the above isomorphism. 
For $r \in R$, we set $\overline{r} := \pi(r) \in B$.  
By the above commutative diagram,  (\#) holds. 
\begin{enumerate}
\item[(\#)] For an element $r \in R$, 
there exists $s \in R$ such that $\overline{r} - \overline{s}\overline{t}^{n-1} \in A$. 
\end{enumerate}
By {\cred (\#) and} $B =R +R\overline{t} + R\overline{t}^2 + \cdots + R \overline{t}^{n-1}$, 
it is enough to prove $\overline{r}\overline{t}, \overline{r}\overline{t}^2, \overline{r}\overline{t}^3, ... \in A$. 
As a weaker statement, we only show $\overline{t}, \overline{t}^2, \overline{t}^3,... \in A$. 
By (\#), 
we have $\overline{t} = \overline{s} \cdot \overline{t}^{n-1} +a$
for some $a \in A$ and $s \in R$. 
It follows from $\overline{t}^n=0$ and $n-1 \geq 1$ 
that  $\overline{t}^{\ell} = (\overline{t} - \overline{s} \cdot \overline{t}^{n-1})^{\ell} =a^{\ell} \in A$ for every integer $\ell \geq 2$. 
Let us show $\overline{t} \in A$. 
By 
$\overline{t} = \overline{s} \cdot \overline{t}^{n-1} +a$, 
it suffices to prove $\overline{s} \cdot \overline{t}^{n-1} \in A$. 
We have shown that  $\overline{t}^{n-1} \in A$ by $n-1 \geq 2$. 
Again by (\#), we can write $\overline{s} = \overline{s'} \cdot  \overline{t}^{n-1} +a'$ 
for some $s' \in R$ and $a' \in A$. 
Then $\overline{s} \cdot \overline{t}^{n-1} = (\overline{s'} \cdot  \overline{t}^{n-1} +a') \cdot \overline{t}^{n-1} = a' \cdot \overline{t}^{n-1}  \in A$ 
by $a' \in A$ and $\overline{t}^{n-1} \in A$. 


For more details on the proofs of 
Theorem \ref{main-ell} and Theorem \ref{intro-t-g1-va}, 
see Section \ref{s-ell}.

\vspace{3mm}

{\cred 
\textbf{Acknowledgements:} 
The author would like to thank the referee for reading the manuscript carefully and for suggesting several  improvements. 
The author was funded by JSPS KAKENHI Grant numbers JP22H01112 and JP23K03028. 
}
\section{Preliminaries}\label{s-prelim}

\subsection{Notation}

\begin{enumerate}
\item We will freely use the notation and terminology in \cite{Har77}, \cite{KM98}, and \cite{Kol13}. 
\item 
Throughout this paper, 
we work over a field $k$ of characteristic $p>0$ unless otherwise specified. 
\item For an integral scheme $X$, 
we define the {\em function field} $K(X)$ of $X$ 
as the local ring $\MO_{X, \xi}$ at the generic point $\xi$ of $X$. 
For an integral domain $A$, $K(A)$ denotes the function field of $\Spec\,A$. 

\item 
For a scheme $X$, its {\em reduced structure} $X_{\red}$ 
is the reduced closed subscheme of $X$ such that the induced closed immersion 
$X_{\red} \to X$ is surjective. 
\item For a field $\kappa$, 
we say that $X$ is a {\em variety over} $\kappa$ if 
$X$ is an integral scheme which is separated and of finite type over $\kappa$. 
We say that $X$ is a {\em curve} (resp. a {\em surface}, resp. a {\em threefold}) over $\kappa$ 
if $X$ is a variety over $\kappa$ of dimension one (resp. {two}, 
resp. {three}). 
\item A variety $X$ over a field is {\em regular} (resp. {\em normal}) if 
the local ring $\MO_{X, x}$ at any point $x \in X$ is regular (resp. an integrally closed domain). 
\item 
For a normal variety $X$ over a field 
{\cred $\kappa$}, 
we define the {\em canonical divisor} $K_X$ as a Weil divisor on $X$ such that 
$\MO_X(K_X) \simeq \omega_{X/{\cred \kappa}}$, where $\omega_{X/{\cred \kappa}}$ denotes the dualising sheaf (cf. \cite[Section 2.3]{Tan18b}). 
Canonical divisors are unique up to linear equivalence. 
Note that $\omega_{X/{\cred \kappa}} \simeq \omega_{X/{\cred \kappa'}}$ for any 
field extension ${\cred \kappa} \subset {\cred \kappa'}$ 
which induces 
a factorisation $X \to \Spec\,{\cred \kappa'} \to \Spec\,{\cred \kappa}$ (\cite[Lemma 2.7]{Tan18b}). 
\item 
We say that a morphism $\psi : V \to W$ of noetherian schemes is {\em birational} 
if there exists an open dense subset $W^{\circ}$ of $W$ such that 
$V^{\circ} := \psi^{-1}(W^{\circ})$ is an open dense subset of $V$ and 
the induced morphism $\psi|_{V^{\circ}} : V^{\circ} \to W^{\circ}$ is an isomorphism.
\item 
Given a noetherian scheme $X$ and a closed point $P$ of $X$, 
$\m_P$ denotes the coherent ideal sheaf on $X$ corresponding to $P$. 
\item 
For ring homomorphisms $k \to \kappa$ of rings and a $k$-scheme $X$, 
we set $X \times_k \kappa := X \times_{\Spec\,k} \Spec\,\kappa$. 
\item 
{\cred Given a finite set $A$, $\# A$ denotes the cardinality of $A$, 
i.e., 
the  number of elements of $A$.}
\end{enumerate}

\subsection{Fundamental cycles}

\begin{dfn}\label{d-fund-cycle}
Let $k$ be a field. 
Let $X$ be a normal surface over $k$ with a unique non-regular point $P$. 
\begin{enumerate}
\item 
For a resolution $f: Y \to X$ (i.e., a projective birational morphism from a regular surface $Y$), 
we say that $Z$ is the {\em {\cred fundamental} cycle}  of $f$ if 
\begin{enumerate}
\item $Z$ is a nonzero $f$-exceptional effective $\Z$-divisor on $Y$, 
    \item $-Z$ is $f$-nef, and 
    \item the inequality $W \geq Z$ holds for every nonzero $f$-exceptional effective $\Z$-divisor $W$ on $Y$ such that $-W$ is $f$-nef. 
\end{enumerate}
The fundamental cycle of $f$ is often denoted by $Z_f$. 
\item For the minimal resolution $f: Y \to X$ of $X$, 
$Z_f$ is called the {\em fundamental cycle of} $P \in X$. 
\end{enumerate}
\end{dfn}

\begin{prop}\label{p-fund-exist}
Let $k$ be a field. 
Let $X$ be a normal surface over $k$ with a unique non-regular point $P$. 
Let $f: Y \to X$ be a resolution of $X$. 
Then there exists the fundamental cycle $Z_f$ of $f$. 
\end{prop}

\begin{proof}
Since the intersection matrix of $\Ex(f)$ is negative definite \cite[Theorem 10.1]{Kol13}, 
we can apply the same argument as in \cite[Lemma 7.2.1 and Lemma 7.2.3]{Ish18}. 
\end{proof}

\begin{prop}\label{p-fund-blowup}
Let $k$ be a field. 
Let $X$ be a normal surface over $k$  with a unique non-regular point $P$. 
Let 
\[
f' : Y' \xrightarrow{g} Y \xrightarrow{f} X
\]
be projective birational morphisms, where both $f$ and $f'$ are resolutions of $X$. 
Then $g^*Z_f = Z_{f'}$. In particular, the positive integer $-Z_f^2$ does not depend on the choice of $f: Y \to X$. 
\end{prop}

\begin{proof}
The same argument as in \cite[Proposition 7.2.6]{Ish18} works. 
\end{proof}

In order to compute multiplicities of elliptic singularities, 
we shall later need the following result, which is mainly extracted from \cite{Wag70}.

\begin{prop}\label{p-fund-mult}
Let $k$ be a field. 
Let $X$ be a normal surface over $k$  with a unique non-regular point $P$. 
Let $f: Y \to X$ be a resolution of $X$. 
Then the following hold. 
\begin{enumerate}
    \item $\m_P \MO_Y \subset \MO_Y(-Z_f)$. 
    \item $\mult\,\MO_{X, P} \geq -Z_f^2$. 
    \item If     $\m_P\MO_Y$ is invertible, then 
    \[
    \mult\,\MO_{X, P} = -W^2, 
    \]
    where $W$ is a Cartier divisor on $Y$ satisying $\m_P\MO_Y \simeq \MO_Y(-W)$. 
\item If $V$ is an effective Cartier divisor on $Y$ such that 
\begin{enumerate}
\item $P \in f(\Supp\,V)$ and 
\item $\MO_Y(-V)$ is $f$-free,  
\end{enumerate}
then $\MO_Y(-V) \subset \m_P \MO_Y$. 
\item If $\MO_Y(-Z_f)$ is $f$-free, then  $\m_P \MO_Y = \MO_Y(-Z_f)$. 
\end{enumerate}
\end{prop}

\begin{proof}
As for (1), \cite[the proof of Proposition 2.3]{Wag70} works. 
In order to show (2), we may replace $f: Y \to X$ by a higher resolution 
(Proposition \ref{p-fund-blowup}), 
and hence we may assume that $\m_P\MO_Y$ is invertible. 
Then (2) and (3) follow from the same argument as in \cite[Theorem 2.7]{Wag70}. 

Let us show (4). 
We may assume that $X$ is affine and  $\MO_Y(-V)$ is globally generated, i.e.,  there is a surjective $\MO_Y$-module homomorphism 
\[
\theta: \MO_Y^{\oplus N} \to \MO_Y(-V), \qquad (a_1, ..., a_N) \mapsto a_1\zeta_1+ \cdots a_N\zeta_N,
\]
\[
\text{where} \qquad
\zeta_1, ..., \zeta_N \in \Gamma(Y, \MO_Y(-V)). 
\]
Via $\Gamma(Y, \MO_Y(-V)) \subset \Gamma(Y, \MO_Y) =\MO_X$, 
(a) implies $\zeta_1, ..., \zeta_N \in \m_P$. 
Fix $Q \in Y$ and $\xi \in \MO_Y(-V)_Q$, where $\MO_Y(-V)_Q$ denotes the stalk at $Q$. 
Since $\theta_Q :\MO_{Y, Q}^{\oplus N} \to \MO_Y(-V)_Q$ is surjective, 
we can write $\xi  = \sum_{i=1}^N b_i \zeta_i$ for some $b_i \in \MO_{Y, Q}$. 
This implies $\MO_Y(-V) \subset \m_P \MO_Y$. 
Thus (4) holds. 
The assertion (5) immediately follows from (1) and (4). 
\end{proof}

We shall need Artin's characterisation for rational singularities 
by using  fundamental cycles (Proposition \ref{p-Artin-rat}). 
We include a proof, because our setting is more generalised than Artin's original one. 
To this end, we start by recalling the following.

\begin{prop}\label{p-Lipman}
Let $A$ be a noetherian ring and let $f: X \to \Spec\,A$ be a proper morphism 
such that $X$ is a notherian scheme, $H^1(X, \MO_X)=0$, and $\dim f^{-1}(y) \leq 1$ for any point $y \in \Spec\,A$. 
Let $L$ be an invertible sheaf on $X$. 
Then the following hold. 
\begin{enumerate}
\item 
If $L$ is $f$-nef, i.e., $L \cdot C \geq 0$ for any one-dimensional integral closed subscheme on $X$ such that $f(C)$ is a point, then $L$ is globally generated. 
\item If $L$ is $f$-nef, then $H^1(V, L) =0$. 
\item If $L$ is $f$-ample, then $L$ is $f$-very ample. 
\end{enumerate}
\end{prop}

\begin{proof}
The assertion follows from \cite[Theorem 12.1]{Lip69}. 
\qedhere


\end{proof}

\begin{prop}\label{p-Artin-rat}
Let $k$ be a field. 
Let $X$ be a normal surface over $k$ with a unique non-regular point $P$. 
Let $f: Y \to X$ be a resolution of $X$. 
Then the following are equivalent. 
\begin{enumerate}
    \item $P \in X$ is a rational singularity, i.e., $R^1f_*\MO_Y =0$. 
    \item $H^1(W, \MO_W)=0$ for every nonzero $f$-exceptional effective $\Z$-divisor $W$ on $Y$. 
    \item $H^1(Z_f, \MO_{Z_f})=0$, where $Z_f$ denotes the fundamental cycle of $f$. 
\end{enumerate}
\end{prop}

\begin{proof}
The implication $(2) \Rightarrow (3)$ is clear. 

Let us show $(1) \Rightarrow (2)$. 
Assume (1). 
By the exact sequence 
\[
0 \to \MO_Y(-W) \to \MO_Y \to \MO_W \to 0, 
\]
we obtain 
\[
0 = R^1f_*\MO_Y \to R^1f_*\MO_W (\simeq H^1(W, \MO_W)) \to {\cred R^2f_*\MO_Y(-Y)=0},  
\]
{\cred where the equality $R^2f_*\MO_Y(-Y)=0$ follows from the fact that 
the dimension of every fibre of $f$ is $\leq 1$.} 
Thus (2) holds.


Let us show $(3) \Rightarrow (1)$. 
Assume (3). 
Set $Z := Z_f$. 
For any $n \in \Z_{>0}$, we have 
\[
0 \to \MO_Y(-nZ)|_Z \to \MO_{(n+1)Z} \to \MO_{nZ} \to 0. 
\]
It follows from $H^1(Z, \MO_Y(-nZ)|_Z)=0$ (Proposition \ref{p-Lipman}) that we obtain $H^1((n+1)Z, \MO_{(n+1)Z}) \simeq  H^1(nZ, \MO_{nZ})$. 
By induction on $n$, (3) implies  $H^1(nZ, \MO_{nZ})=0$ for every $n \in \Z_{>0}$. 
By the formal function theorem, we obtain 
\[
R^1f_*\MO_Y \otimes_{\MO_{X,P}} \widehat{\MO}_{X, P} \simeq \varprojlim_n H^1(nZ, \MO_{nZ}) =0. 
\]
Since $\MO_{X, P} \to \widehat{\MO}_{X, P}$ is faithfully flat, we get 
$R^1f_*\MO_Y =0$. Thus (1) holds. 
\end{proof}

\subsection{Conductor exact sequences}

\begin{dfn}
Let $A \subset B$ be a ring extension. 
We set 
\[
I_{B/A} := \{ b \in B \,|\, bB \subset A\}. 
\]
Then we see that $I_{B/A} \subset A$ and $I_{B/A}$ is an ideal of each of $A$ and $B$. 
{\cred Indeed, $I_{B/A}$ is clearly an ideal of $B$ 
and the inclusion $I_{B/A} \subset A$ holds, 
because $b \in I_{B/A}$ implies $b=b \cdot 1 \in bB \subset A$. 
These properties automatically imply that $I_{B/A}$ is an ideal of $A$. 
In particular, we get} 
\[
{\cred I_{B/A} = \{ a \in A \,|\, aB \subset A\}.} 
\]
We call $I_{B/A}$ the {\em conductor ideal} of $A \subset B$. 
\end{dfn}

For a ring extension $A \subset B$, we can directly check the following properties. 

\begin{enumerate}
\item[(I)] $I_{B/A} = A$ if and only if $A=B$. 
\item[(II)] Given a multiplicatively closed subset $S$ of $A$, 
we have $I_{S^{-1}B/S^{-1}A} = S^{-1}I_{B/A}$. 
\end{enumerate}

\begin{dfn}
Let $\nu: X \to Y$ be a finite morphism of noetherian schemes 
such that $\MO_Y \to \nu_*\MO_X$ is injective. 
By (II), there exists a unique coherent ideal sheaf $I_Y$ on $Y$ such that $\Gamma(Y_1, I_Y) = I_{\MO_X(\nu^{-1}(Y_1))/\MO_Y(Y_1)}$  for any affine open subset $Y_1$. 
The coherent sheaves $I_Y$ and $I_X := I_Y \MO_Y$ are called the {\em conductor ideals} (of $\nu$). 
By construction, we get $\nu_*I_X = I_Y$. 
\end{dfn}

Let $C_X$ and $C_Y$ be the closed subschemes on $X$ and $Y$ corresponding to $I_X$ and $I_Y$, respectively. 
Then we have the following cartesian diagram: 
    \begin{equation}\label{e1-cond}
    \begin{tikzcd}
    C_X \arrow[r, hook] \arrow[d, "\nu|_{C_X}"] & X \arrow[d, "\nu"]\\
    C_Y \arrow[r, hook] & Y. 
\end{tikzcd}
    \end{equation}
By $\nu_*I_X = I_Y$ and the snake lemma, 
we obtain the following commutativa diagram in which each horizontal sequence is exact and the vertical arrows are injective: 
\begin{equation}\label{e2-cond}
\begin{tikzcd}
0 \arrow[r] & \nu_*I_X \arrow[r] & \nu_*\MO_X \arrow[r] & \nu_*\MO_{C_X} \arrow[r] & 0\\
0 \arrow[r] &  I_Y \arrow[r] \arrow[u, equal] & \MO_Y \arrow[r] \arrow[u, hook] & \MO_{C_Y} \arrow[r] \arrow[u, hook]  & 0.
\end{tikzcd}
\end{equation}
By diagram chase, we obtain the following exact sequence: 
\begin{equation}\label{e3-cond}
0 \to \MO_Y \to \nu_*\MO_X \oplus \MO_{C_Y} \xrightarrow{\delta} \nu_*\MO_{C_X} \to 0, 
\end{equation}
where $\delta$ is defined by the difference. 
We shall call this {\em the conductor exact sequence}.


\subsection{Curves over an imperfect field}\label{ss-known}

In this subsection, we work over a field $k$. 
We summarise some known results on curves over $k$. 
All the results in this subsection are well known to experts.

Let $X$ be a one-dimensional projective scheme over $k$. 
Let $L$ be a Cartier divisor or an invertible sheaf on $X$. 
The {\em degree}  of $L$, denoted by $\deg_k L$ or $\deg L$, 
is defined by 
\[
\deg L := \deg_k L  := \chi(X, L) -\chi(X, \MO_X). 
\]
The equation 
\[
\chi(X, L) = \deg_k L +\chi(X, \MO_X)
\]
is called the {\em Riemann-Roch theorem}, although this is obvious under our definition of $\deg_k L$. 
Note that $\deg_k L$ actually depends on the base field $k$, i.e., 
if $k'$ is a subfield of $k$ with $[k:k']<\infty$, 
then $X$ is a one-dimensional projective scheme {\cred over} $k'$ and 
\[
\deg_{k'} L = [k:k']\deg_k L.
\]


For  a regular projective curve $X$ over $k$ and a Cartier divisor $D$ on $X$, 
$\deg D$ coincides with the usual definition, i.e., for the irreducible decomposition $D = \sum_{i=1}^r a_iP_i$, the following hold: 
\[
\deg D = \deg_k D = \sum_{i=1}^r a_i \dim_k H^0(P_i, \MO_{P_i}). 
\]
This can be proven by the same proof as the case over an algebraically closed field \cite[Ch. IV, Theorem 1.3]{Har77}. 
If $X$ is a (possibly non-regular) projective curve over $k$ and $D$ is a Cartier divisor on $X$, then 
we have $\deg_k D = \deg_k (\nu^*D)$ (\ref{e3-cond}), 
where $\nu : X^N \to X$ denotes the normalisation of $X$. 
In particular, $\deg_k P = \dim_k H^0(P, \MO_P)$ for a closed point $P$ on $X$ around which $X$ is regular.

We shall often use the following result on 
the Castelnuovo-Mumford regularity.






\begin{prop}\label{p-regularity}
Let $k$ be a field and let $X$ be a 
one-dimensional projective scheme over $k$ with $H^0(X, \MO_X)=k$. 
Let $A$ be an globally generated ample Cartier divisor on $X$. 
Fix $m \in \Z_{\geq 0}$ and assume that 
\[
H^1(X, \MO_X(mA))=0. 
\]
Then the $k$-algebra 
\[
\bigoplus_{n=0}^{\infty}H^0(X, \MO_X(nA))
\]
is generated by $\bigoplus_{n=1}^{m+1}H^0(X, \MO_X(nA))$.  
\end{prop}

\begin{proof}
Note that $\MO_X(A)$ is ${\cred (m+1)}$-regular with respect to $\MO_X(A)$, so that 
\[
H^0(X, \MO_X(A)) \otimes_k H^0(X, \MO_X((m+r+1)A)) \to H^0(X, \MO_X((m+r+2)A))
\]
is surjective for every $r \geq 0$ 
\cite[Section 5.2]{FGI05}, 
\cite[Definition 11.1, Lemma 11.2]{Tan21}. 
\end{proof}

\section{Curves of genus zero}

\begin{dfn}
Let $k$ be a field. 
We say that $X$ is a {\em curve of genus zero} or a {\em genus-zero curve} (over $k$) 
if $X$ is a projective curve over $k$ such that 
$H^1(X, \MO_X)=0$. 
\end{dfn}



\subsection{Divisors}

\begin{prop}\label{p-g0-bpf}
Let $k$ be a field. 
Let $X$ be a one-dimensional projective  scheme over $k$ such that $H^1(X, \MO_X)=0$. 
Let $L$ be an  invertible sheaf on $X$. 
Then the following holds. 
\begin{enumerate}
    \item If $L$ is nef, i.e., $\deg (L|_C) \geq 0$ for every curve $C$ contained in $X$, then $H^1(X, L)=0$ and $L$ is globally generated. 
    \item If $L$ is ample, then $L$ is very ample. 
    \item If $L$ is ample and $H^0(X, \MO_X)=k$, then 
    $\bigoplus_{m =0}^{\infty} H^0(X, L^{\otimes m})$ is generated by $H^0(X, L)$ as a $k$-algebra. 
    \item Assume that $L$ is nef. Consider the morphism induced by $|L|$: 
\[
\varphi_{|L|} : X \xrightarrow{\psi} Y  \hookrightarrow  \P^N_k,  
\]
where 
$N := h^0(X, L) -1$, 
$Y := \varphi_{|L|}(X)$ denotes the scheme-theoretic image of $X$ by $\varphi_{|L|}$, 
and $\psi$ is the induced morphism. 
Then  $\psi_*\MO_X = \MO_Y$ and $H^1(Y, \MO_Y)=0$. 
\end{enumerate}
\end{prop}

\begin{proof}
The assertions (1) and (2) immediately follow from Proposition \ref{p-Lipman}. 
Let us show (3). 
By (1), $L$ is an ample globally generated invertible sheaf. 
Then the assertion follows from 
$H^1(X, \MO_X)=0$ and Proposition \ref{p-regularity}. 
Thus (3) holds.

Let us show (4).  
Let $\psi : X \to Z \to Y$ be the Stein factorisation of $\psi$. 
We obtain $H^1(Z, \MO_Z) \hookrightarrow H^1(X, \MO_X)=0$, which implies $H^1(Z, \MO_Z) =0$. 
Then  the pullback $L_Z$ of $\MO_{\P^N}(1)$ to $Z$ is very ample by (2). 
By $H^0(X, L) \simeq H^0(Z, L_Z)$, we obtain $Z \simeq Y$. 
Thus (4) holds. 
\end{proof}

\begin{prop}\label{p-g0-H1}
Let $k$ be a field and 
let $X$ be a one-dimensional projective scheme over $k$ with $H^1(X, \MO_X)=0$. 
Let $L$ be a nef invertible sheaf. 
Then the following hold. 
\begin{enumerate}
    \item $H^1(X, L)=0$. 
    \item 
$h^0(X, L) = h^0(X, \MO_X) + \deg_k L$. 
\end{enumerate}
\end{prop}

\begin{proof}
The assertion (1) follows from Proposition \ref{p-g0-bpf}. 
Then (1) and the Riemann--Roch theorem imply (2). 
\end{proof}

\begin{lem}
Let $k$ be a field and let $X$ be a genus-zero curve. 
Let $\nu : Y \to X$ be the normalisation. 
Then also $Y$ is a genus-zero curve. 
\end{lem}

\begin{proof}
The assertion immediately follows from  the conductor exact sequence (\ref{e3-cond}). 
\end{proof}





\begin{prop}\label{p-g0-Pic}
Let $k$ be a field and let $X$ be a genus-zero curve. 
Then the following hold. 
\begin{enumerate}
    \item Let $P$ be a closed point on $X$ such that 
    \begin{enumerate}
    \item $X$ is regular around $P$, and 
   \item $\deg_k P = \min_Q \{ \deg_k Q\}$, where $Q$ runs over all the closed points on $X$ around which $X$ is regular. 
    \end{enumerate}
    Then $\Pic\,X \simeq \Z \MO_X(P)$. 
    \item If $k$ is an infinite field, then there exist infinitely many closed points $P$ on $X$ such that {\rm (a)} and {\rm (b)} hold. 
\end{enumerate}
\end{prop}

\begin{proof}
Let us show (1). 
Fix a closed point $P$ satisfying (a) and (b). 
Take a Cartier divisor $D$. 
It suffices to show $D \sim mP$ for some $m \in \Z$. 
Since we have $D = H - H'$ for ample Cartier divisors $H$ and $H'$, 
we may assume that $D$ is ample. 
By Proposition \ref{p-g0-bpf}, $|D|$ is base point free. 
Hence we may assume that $D$ is an effective Cartier divisor such that 
$X$ is regular around $\Supp\,D$. 
Then the problem is further reduced to the case when $Q := D$ is a closed point around which $X$ is  regular. 
If $\deg_k Q$ is divisible by $\deg_k P$, 
then we obtain $Q \sim m P$ by Proposition \ref{p-g0-bpf}. 

Suppose that $\deg_k Q$ is not divisible by $\deg_k P$. 
It suffices to derive a contradiction. 
By taking the greatest common divisor $d$ of $\deg_k P$ and $\deg_k Q$, 
we can find an ample Cartier divisor $E$ on $X$ such that $\deg_k E =d$. 
Automatically, we get $d < \deg_k P$. 
Since $|E|$ is base point free (Proposition \ref{p-g0-bpf}), 
we may assume that $E$ is effective and 
$X$ is regular around $\Supp\,E$. 
Pick a closed point $R \in \Supp\,E$. 
Then 
\[
\deg_k R \leq \deg_k E = d < \deg_k P,
\]
which contradicts the assumption (b) on $P$. Thus (1) holds. 

Let us show (2). 
We can find 
a closed point $P$ satisfying (a) and (b). 
Since $|P|$ is base point free (Proposition \ref{p-g0-bpf}), 
the assertion (2) follows from $\dim_k H^0(X, \MO_X(P)) > \dim_k H^0(X, \MO_X)$ (Proposition \ref{p-g0-H1}).  
\end{proof}



\subsection{A characterisation via $\Delta$-genus}

\begin{dfn}\label{d-Delta}
Let $k$ be a field. 
Let $X$ be a projective variety with $H^0(X, \MO_X)=k$ 
and 
let $D$ be an ample Cartier divisor on $X$. 
We set 
\[
\Delta(X, D) := \dim X + D^{\dim X} - h^0(X, D), 
\]
which is called the {\em $\Delta$-genus} of $(X, D)$. 
Recall that $h^0(X, D) := \dim_k H^0(X, D)$. 
\end{dfn}

\begin{lem}\label{l-Delta-nonneg}
Let $k$ be a field and let $X$ be a projective curve with $H^0(X, \MO_X)=k$ 
and let $D$ be an ample Cartier divisor on $X$. 
Then $\Delta(X, D) \geq 0$. 
\end{lem}

\begin{proof}
Suppose $\Delta(X, D) < 0$. By $\Delta(X, D) = 1 + \deg D  - h^0(X, D)$, 
we obtain $h^0(X, D) > 1+\deg D \geq 2$. 
In particular, we may assume that $D$ is effective. 
By the exact sequence 
\[
0 \to \MO_X \to \MO_X(D) \to \MO_X(D)|_D \to 0, 
\]
we get a surjection $H^1(X, \MO_X) \to H^1(X, \MO_X(D))$, 
{\cred where $H^1(D, \MO_X(D)|_D)=0$ holds by $\dim D = 0$}. 
In particular, $h^1(X, \MO_X) \geq h^1(X, \MO_X(D))$. 

In order to derive a contradiction, 
it suffices to show $h^1(X, \MO_X) < h^1(X, \MO_X(D))$. 
The Riemann--Roch theorem implies 
\[
h^0(X, D) - h^1(X, D) = \chi(X, D) = \chi(X, \MO_X) + \deg D = 1 -h^1(X, \MO_X)+ \deg D. 
\]
Therefore, we get 
\[
h^1(X, D) - h^1(X, \MO_X) = h^0(X, D) -(1+\deg D) = -\Delta(X, D) > 0, 
\]
as required. 
\end{proof}

\begin{prop}\label{p-Delta0-curve}
Let $k$ be a field. 
Let $X$ be a projective curve over $k$ with $H^0(X, \MO_X)=k$. 
Then the following are equivalent. 
\begin{enumerate}
    \item $H^1(X, \MO_X)=0$. 
    \item $\Delta(X, D) = 0$ for every ample Cartier divisor $D$ on $X$. 
    \item $\Delta(X, D) = 0$ for some ample Cartier divisor $D$ on $X$. 
\end{enumerate}
\end{prop}

\begin{proof}
By the Riemann--Roch theorem and Proposition \ref{p-g0-H1}, we get (1) $\Rightarrow$ (2). 
It is clear that (2) $\Rightarrow$ (3).

Let us show that (3) $\Rightarrow$ (1). 
Assume (3), i.e., 
there exists an ample Cartier divisor $D$ on $X$ such that $1 + \deg D- h^0(X, D) = \Delta(X, D) =0$. 
By $h^0(X, D) = 1+ \deg D >0$, we may assume that $D$ is effective. 
The Riemann--Roch theorem  implies 
\[
h^0(X, D) - h^1(X, D) = \chi(X, D) = \chi(X, \MO_X) + \deg D = 1 -h^1(X, \MO_X)+ \deg D. 
\]
This, together with $h^0(X, D) = 1+ \deg D$, implies $h^1(X, D) = h^1(X, \MO_X)$. 
Consider the following exact sequence: 
\[
0 \to H^0(X, \MO_X) \to H^0(X, \MO_X(D)) \xrightarrow{\alpha} H^0(D, \MO_X(D)|_D) 
\]
\[
\to 
H^1(X, \MO_X) \xrightarrow{\simeq} H^1(X, \MO_X(D)) \to 0. 
\]
Then $\alpha$ is surjective and hence $|D|$ is base point free. 
In particular, we can find an effective Cartier divisor $D'$ such that $D \sim D'$ and 
$\Supp\,D \cap \Supp\,D' = \emptyset$. 
For every $n \in \Z_{\geq 0}$, we get the following commutative diagram 
consisting of the induced maps: 
\[
\begin{tikzcd}
H^0(X,\mathcal{O}_X(D)) \arrow[r, "\alpha"] \arrow[d, hook] & H^0(D,\mathcal{O}_X(D)|_D) \arrow[d, "\simeq"] \\
H^0(X,\mathcal{O}_X(D+nD')) \arrow[r, "\beta_n"] & H^0(D,\mathcal{O}_X(D+nD')|_D). 
\end{tikzcd}
\]
Since $\alpha$ is surjective, so is $\beta_n$. 
Recall that $\beta_n$ is obtained from the following exact sequence: 
\[
0 \to H^0(X, \MO_X(nD')) \to H^0(X, \MO_X(D+nD')) \xrightarrow{\beta_n} H^0(D, \MO_X(D+nD')|_D) 
\]
\[
\to 
H^1(X, \MO_X(nD')) \xrightarrow{\theta} H^1(X, \MO_X(D+nD')) \to 0. 
\]
As $\beta_n$ is surjective, we obtain 
\[
H^1(X, \MO_X(nD)) \simeq H^1(X, \MO_X(nD')) \overset{\theta}{\simeq} 
H^1(X, \MO_X(D+nD')) \simeq H^1(X, \MO_X( (n+1)D). 
\]
Since $n$ is chosen to be an arbitrary non-negative integer, 
we get 
\[
H^1(X, \MO_X) \simeq H^1(X, \MO_X(D)) \simeq H^1(X, \MO_X(2D)) 
\simeq \cdots \simeq H^1(X, \MO_X(mD)) =0
\]
for $m \gg 0$, where the last equality follows from the Serre vanishing theorem. 
Thus (1) holds. 
\end{proof}

\subsection{Examples}

Given a genus-zero curve, 
being  Gorestein is equivalent to being a plane conic. 

\begin{ex}\label{e-g0-conic}
Let $k$ be a field and let $X$ be a Gorenstein genus-zero curve over $k$ with 
$H^0(X, \MO_X)=k$. 
By \cite[Lemma 10.6]{Kol13}, one of the following holds. 
\begin{enumerate}
\item $X \simeq \P^1_k$. 
\item $X$ is a conic on $\P^2_k$, $\Pic\,X \simeq \Z$, and $\Pic\,X$ is generated by $\omega_X$. 
\end{enumerate}
\end{ex}


The following construction gives non-Gorenstein genus-zero curves. 

\begin{ex}
Let $\ell$ be a field and set $Y := \P^1_{\ell}$. 
Fix an $\ell$-rational point $P \in Y(\ell)$. 
Let $k \subset \ell$ be a field extension with $[\ell : k] <\infty$. 
Then the pinching $X$ of $Y$ along $k \subset \ell$ is a projective curve over $k$. 
Since $X$ has a $k$-rational point, we have $H^0(X, \MO_X)=k$. 

Let us show $H^1(X, \MO_X)=0$. 
We obtain the conductor exact sequence (\ref{e3-cond}):
\[
0 \to \MO_X \to \nu_*\MO_Y \oplus \MO_{\Spec\,k} \to \nu_*\MO_{\Spec\,\ell} \to 0. 
\]
The induced sequence 
\[
0 \to H^0(\MO_X) \to H^0(\MO_Y) \oplus H^0(\MO_{\Spec\,k}) \to 
H^0(\MO_{\Spec\,\ell}) \to 0
\]
is still exact by 
$\dim_k H^0(\MO_X) = \dim_k H^0(\MO_{\Spec\,k}) =1$ 
and 
\[
\dim_k H^0(\MO_Y) = \dim_k H^0(\MO_{\Spec\,\ell}) = [\ell :k].
\]
We then get an injection 
\[
H^1(X, \MO_X) \hookrightarrow H^1(Y, \MO_Y) \oplus H^1(\Spec\,k, \MO_{\Spec\,k}) =0, 
\]
i.e., $H^1(X, \MO_X)=0$. 
Note that $X$ is Gorenstein if and only if $[k : \ell] =2$ \cite[Proposition 1.1]{Scht}. 
\end{ex}

\section{Gorenstein curves of genus one}\label{s-genus-one-prime}

\subsection{Base point freeness}

\begin{dfn}\label{d-genus1}
Let $k$ be a field. 
We say that $X$ is a {\em curve of genus one} or a {\em genus-one curve} (over $k$)
if $X$ is a projective curve over $k$ such that $
\dim_k H^1(X, \MO_X)=1$. 
\end{dfn}

\begin{rem}
Let $k$ be a field and let $X$ be a genus-one curve over $k$. 
Then it holds that $H^0(X, \MO_X)=k$. 
Indeed, $H^1(X, \MO_X)$ is an $\ell$-vector space for $\ell := H^0(X, \MO_X)$, 
which implies that $[\ell : k] \dim_{\ell} H^1(X, \MO_X) =\dim_k H^1(X, \MO_X)=1$, and hence $\ell = k$. 
\end{rem}

\begin{lem}\label{l-g1-bpf}
Let $k$ be a field. 
Let $X$ be a Gorenstein genus-one curve over $k$. 
Let $D$ be a Cartier divisor with $\deg D \geq 2$. 
Then the following hold. 
\begin{enumerate}
    \item $\Bs\,|D|$ is contained in the non-regular locus of $X$. 
    \item If there exists a regular closed point $P$ of $X$ such that $\deg D > \deg P$, then $|D|$ is base point free. 
\end{enumerate}
\end{lem}

\begin{proof}
Let us show (1). 
    By the Riemann--Roch theorem, we may assume that $D$ is effective. 
Fix a closed point $P \in \Supp\,D$ around which $X$ is regular.  
It is enough to show $P \not\in \Bs\,|D|$. 
We have 
\[
H^0(X, D) \xrightarrow{{\cred \rho}} H^0(P, D|_P) \to H^1(X, D-P) \to H^1(X, D) =0, 
\]
where $H^1(X, D)=0$ follows from Serre duality. 
{\cred If $D \neq P$, then we get $D >P$ as $D$ is effective, 
which implies $H^1(X, D-P)=0$ again by Serre duality.} 
Hence the problem is reduced to the case when $D=P$. 
Then $h^0(P, D|_P) = \deg P =\deg D \geq 2$, whilst $h^1(X, D-P) = h^1(X, \MO_X)=1$. 
Therefore, the restriction map 
\[
\rho : H^0(X, D) \to H^0(P, D|_P)
\]
is nonzero. 
Take an effective divisor $D'$ such that the corresponding element $s_{D'} \in H^0(X, D)$ satisfies $\rho(s_{D'}) \neq 0$. 
We then obtain $P \not\in \Supp\,D'$. 
Therefore, $P \not\in \Bs\,|D|$.  
Thus (1) holds. 

Let us show (2). 
Consider the following exact sequence: 
\[
0 \to H^0(X, \MO_X(D-P)) \to H^0(X, \MO_X(D)) \to H^0(P, \MO_X(D)|_P). 
\]
By 
\[
h^0(X, \MO_X(D)) =\chi(X, \MO_X(D)) =\deg D > \deg P = h^0(P, \MO_X(D)|_P), 
\]
we get $H^0(X, \MO_X(D-P)) \neq 0$. 
In particular, the problem is reduced to the case when $D = P+D'$ for some effective Cartier divisor $D'$. 
By (1), we get 
\[
\Bs\,|D| \subset X_{{\rm non-reg}} \cap \Supp(P+D') \subset \Supp\,D', 
\]
where the latter inclusion holds by the fact that $X$ is regular around $P$. 
Thus it suffices to show that $\Bs\,|D| \cap \Supp\,D' = \emptyset$. 
By the exact sequence 
\[
0 \to \MO_X(D -D') \to \MO_X(D) \to \MO_X(D)|_{D'} \to 0, 
\]
we obtain another exact sequence 
\[
H^0(X, \MO_X(D)) \to H^0(D', \MO_X(D)|_{D'}) 
\to H^1(X, \MO_X(D -D')) = H^1(X, \MO_X(P))=0. 
\]
Therefore $\Bs\,|D|$ is disjoint from $D'$, as required. 
Thus (2) holds 
\end{proof}

\begin{thm}\label{t-g1-bpf}
Let $k$ be a field. 
Let $X$ be a Gorenstein genus-one curve over $k$. 
Let $D$ be a Cartier divisor with $\deg D \geq 2$. 
Then $|D|$ is base point free. 
\end{thm}

\begin{proof}
Taking the base change to the separable closure, 
the problem is reduced to the case when $k$ is separably closed. 
If $X$ is geometrically integral, then the assertion follows from 
\cite[Lemma 11.10(3)]{Tan21}. 
Hence we may assume that $X$ is not geometrically integral. 
In particular, $k$ is an infinite field.

Let $\nu : Y \to X$ be the normalisation. 
Let $C_X \subset X$ and $C_Y \subset Y$ be the closed subschemes defined by the conductor of $\nu$. 
We get the following cartesian diagram: 
\[
\begin{CD}
Y @<<< C_Y\\
@VVV @VVV\\
X @<<< C_X. 
\end{CD}
\]
Set $k_Y := H^0(Y, \MO_Y)$. 
Recall that we have  
\[
K_Y + C_Y \sim \nu^*K_X \sim 0. 
\]
{\cred If $X$ is regular, then the assertion follows from Lemma \ref{l-g1-bpf}(1). Then we may assume that $C_Y \neq 0$ and hence} $Y$ is a regular curve of genus zero.

\begin{step}\label{s01-g1-bpf}
In order to prove   Theorem \ref{t-g1-bpf}, we may assume that 
the following hold. 
\begin{enumerate}
\renewcommand{\labelenumi}{(\alph{enumi})}
\item $D$ is effective. 
    \item $\Supp\,D \subset \Supp\,C_X$. 
    \item $\deg D = \min_P(\deg P)$, where $P$ runs over all the closed points of $Y$. 
    \item $\# D =1$ and $\#(\nu^*D) =1$. 
\end{enumerate}
\end{step}

\begin{proof}[Proof of Step \ref{s01-g1-bpf}]
By the Riemann--Roch theorem, we may assume (a). 
If $\Supp\,D \not\subset \Supp\,C_X$, then $|D|$ is base point free by 
Lemma \ref{l-g1-bpf}{\cred (2)}. 
Hence we may assume (b).

As for (c), we have $\deg D = \deg (\nu^*D) \geq \min_P(\deg P)$ by (a). 
Assume that $\deg D > \min_P(\deg P)$. 
Note that there are infinitely many closed points $Q$ on  $Y$ such that 
$\min_P(\deg P) = \deg Q$ (Proposition \ref{p-g0-Pic}). 
Then we can find a closed point $Q$ on $Y$ 
such that $\min_P(\deg P) = \deg Q$ and $\nu$ is isomorphic around $Q$. 
Set $Q_X := \nu(Q)$. 
Then $X$ is regular around $Q_X$ and $\deg D >  \min_P(\deg P) = \deg Q = \deg Q_X$, 
which implies that $|D|$ is base point free (Lemma \ref{l-g1-bpf}{\cred (2)}). 
Thus we may assume (c). 

Since (c) holds, also (d) holds. 
Indeed, $\#(\nu^*D) =1$ follows from $\deg (\nu^*D) = \deg D \overset{{\rm (c)}}{=} \min_P(\deg P)$. 
We then get $\# D =1$ by $1 \leq \# D \leq \#(\nu^*D) =1$. 
This completes the proof of Step \ref{s01-g1-bpf}. 
\end{proof}

\begin{step}\label{s02-g1-bpf}
The assertion of  Theorem \ref{t-g1-bpf} holds when $\# C_X=1$. 
\end{step}

\begin{proof}[Proof of Step \ref{s02-g1-bpf}]
Assume $\# C_X=1$. 
In this case, $P:=(C_X)_{\red}$ is a unique non-regular closed point on $X$. 
Fix a regular closed point $Q \in X$ such that $\deg Q = \deg D$,
 whose existence is guaranteed by Proposition \ref{p-g0-Pic}. 
We have  $Q \neq (C_X)_{\red} =P$. 
For 
\[
\Spec\,A := X \setminus Q, \qquad \Spec\,\overline{A} := X \setminus P = X \setminus (C_X)_{\red}, 
\]
we get an affine open cover $X = \Spec\,A \cup \Spec\,\overline{A}$. 
Take their inverse images to $Y$: 
\[
\Spec\,B := \nu^{-1}(\Spec\,A), \qquad \Spec\,\overline{B} := \nu^{-1}(\Spec\,\overline{B}). 
\]
In particular, we have $ Y = \Spec\,B \cup \Spec\,\overline{B}$ and $\overline{A} = \overline B$. 
By Step \ref{s01-g1-bpf}, we can find $f \in A$ such that 
\[
\MO_X(-D)|_{\Spec\,A} = fA, \qquad \MO_X(-D)|_{\Spec\,\overline{A}} = \overline{A}
\]
\[
\MO_Y(-\nu^*D)|_{\Spec\,B} = fB, \qquad \MO_Y(-\nu^*D)|_{\Spec\,\overline{B}}=\overline{B}. 
\]
Set $P_Y := (\nu^{-1}(P))_{\red}$ and $Q_Y := \nu^{-1}(Q)$. 
Both $P_Y$ and $Q_Y$ are closed points (Step \ref{s01-g1-bpf}(d)). 
By Step \ref{s01-g1-bpf}, we have an equality 
$\nu^*D|_{\Spec\,B} ={\rm div}(f)|_{\Spec\,B} =  P_Y$. 
We then have ${\rm div}(f) = P_Y + m Q_Y$ for some $m \in \Z$. 
We obtain $m=-1$, i.e., ${\rm div}(f) = P_Y-Q_Y$, 
because  $\deg Q_Y = \deg Q = \deg D = \deg P_Y$ and $0 = \deg ({\rm div}(f)) = \deg P_Y + m \deg Q_Y$. 
In particular, $f^{-1} \in \overline B = \overline A$. 
Set $g := 1+f \in A$. 
Then 
\[
\frac{g}{f} = \frac{1+f}{f} = 1 + f^{-1} \in \overline{B} = \overline{A}. 
\]
Therefore, we obtain $\frac{g}{f} \in H^0(X, \MO_X(D))$. 
Indeed, 
\begin{itemize}
    \item $\frac{g}{f} \in f^{-1}A =\MO_X(D)|_{\Spec\,A}$, and 
    \item $\frac{g}{f} \in \overline{A} =\MO_X(D)|_{\Spec\,\overline{A}}$. 
\end{itemize}
Since $\frac{g}{f}$ does not vanish at the point $P  =D_{\red}$, 
$|D|$ is base point free. 
This completes the proof of Step \ref{s02-g1-bpf}. 
\end{proof}

\begin{step}\label{s03-g1-bpf}
The assertion of  Theorem \ref{t-g1-bpf} holds when $\# C_X \geq 2$. 
\end{step}

\begin{proof}[Proof of Step \ref{s03-g1-bpf}]
Assume that $\# C_X \geq 2$. 
By ${\cred \# C_Y \geq \# C_X \geq 2}$, $K_Y +C_Y \sim 0$, and Example \ref{e-g0-conic}, 
we obtain $Y \simeq \P^1_{k_Y}$ and $C_Y = C_{Y, 1} + C_{Y, 2}$, 
where $k_Y :=H^0(Y, \MO_Y)$  and each of $C_{Y, 1}$ and $C_{Y, 2}$ is a $k_Y$-rational point. 
In particular, 
$\# C_X = \# C_Y =2$. 
For each $i \in \{1, 2\}$,  
\[
k_i := \MO_{C_{X, i}} \hookrightarrow \MO_{C_{Y, i}} =k_Y
\]
is a field extension, 
where each $C_{X, i}$ denotes the connected component of $C_X$ containing $\nu(C_{Y, i})_{\red}$. 
We may assume that $\Supp\,D = \Supp\,C_{X, 1}$. 


Set 
\[
\Spec\,B := Y \setminus C_{Y, 2}, \quad 
\Spec\,\overline{B} := Y \setminus C_{Y, 1}, \quad 
\Spec\,A := X \setminus C_{X, 2}, \quad 
\Spec\,\overline{A} := X \setminus C_{X, 1}. 
\]
We introduce affine coordinates 
\[
B =k_Y[x], \qquad \overline{B} = k_Y[y], 
\]
where the gluing is given by the usual relation $xy=1$. We get 
\[
A = k_1 \oplus xk_Y[x] \subset k_Y[x] =B, \qquad 
\overline{A} = k_2 \oplus yk_Y[y] \subset k_Y[y] =\overline{B}. 
\]
We can write  
\begin{eqnarray*}
    {\cred \MO_X(-D)}|_{\Spec\,A} &=& (b_1x + b_2 x^2 + \cdots +b_Nx^N)A 
\quad \text{and}\quad\\
{\cred \MO_Y(-\nu^*D)}|_{\Spec\,B} &=& (b_1x + b_2 x^2 + \cdots +b_Nx^N)B
\end{eqnarray*}
for some $b_1, b_2, ..., b_N \in k_Y$. 
By Step \ref{s01-g1-bpf}, 
we obtain an equality of principal ideals: $(b_1x + b_2 x^2 + \cdots +b_Nx^N)B = xB$. 
Then $(b_1x + b_2 x^2 + \cdots +b_Nx^N)g(x) = x$ holds for some $g(x) \in k_Y[x]$, so that $b_1x + b_2 x^2 + \cdots +b_Nx^N = b_1x$, i.e., $N=1$. 
Set $b:= b_1 \in k_Y^{\times}$. We get 
\[
\MO_X(-D)|_{\Spec\,A} = bx A 
\quad \text{and}\quad
\MO_Y(-\nu^*D)|_{\Spec\,B} = bxB (=xB).  
\]
Note that we possibly have $bx A \neq xA$. 
We obtain 
\[
\zeta := \frac{b(x+b^{-1})}{bx} \in H^0(X, \MO_X(D)). 
\]
Indeed, the following hold. 
\begin{itemize}
    \item $\zeta = \frac{b(x+b^{-1})}{bx} = \frac{1 +bx}{bx} 
    \in (bx)^{-1}(k_1 \oplus xk_Y[x]) = (bx)^{-1}A = \MO_X(D)|_{\Spec\,A}$. 
    \item $\zeta|_{\Spec\,\overline{A}} = \frac{b(y^{-1}+b^{-1})}{by^{-1}} = 1 + b^{-1}y \in k_2 + yk_Y[y] = \overline{A} = \MO_X(D)|_{\Spec\,\overline{A}}$. 
\end{itemize}
Since $\zeta$ does not vanish at the point $(C_{X,1})_{\red} =D_{\red}$, 
$|D|$ is base point free. 
This completes the proof of Step \ref{s03-g1-bpf}. 
\end{proof}
Step \ref{s02-g1-bpf} and 
Step \ref{s03-g1-bpf} complete the proof of Theorem \ref{t-g1-bpf}. %
\end{proof}


\subsection{Very ampleness and projective normality}

\begin{lem}\label{l-g1-deg2-gene12}
Let $k$ be a field and let $X$ be a Gorenstein genus-one curve over $k$. 
Let $D$ be a Cartier divisor on $X$ with $\deg D \geq 2$. 
Then 
\[
\bigoplus_{m=0}^{\infty} H^0(X, \MO_X(mD))
\] 
is generated by $H^0(X, D) \oplus H^0(X, 2D)$ as a $k$-algebra. 
\end{lem}

\begin{proof}
By Theorem \ref{t-g1-bpf}, $D$ is a globally generated ample Cartier divisor. 
The assertion follows from 
Proposition \ref{p-regularity} and $H^1(X, \MO_X(D))=0$. 
\end{proof}

\begin{nota}\label{n-g1-deg3-birat}
Let $k$ be a field and let $X$ be a Gorenstein genus-one curve over $k$. 
Let $D$ be a Cartier divisor on $X$ with $\deg D \geq 3$. 
Note that $|D|$ is base point free (Theorem \ref{t-g1-bpf}). 
Let 
\[
\varphi_{|D|} : X \xrightarrow{\psi} Y  \overset{j}{\hookrightarrow} \P^{N}_k
\]
be the induced morphisms, where 
$N := h^0(X, D) -1  = \deg D -1$, $Y := \varphi_{|D|}(X)$, and 
$j$ denotes the induced closed immersion. 
\end{nota}

\begin{lem}\label{l-g1-deg3-birat}
We use Notation \ref{n-g1-deg3-birat}. 
Then $\psi : X \to Y$ is birational. 
\end{lem}

\begin{proof}
Fix an ample Cartier divisor $D_Y$ such that $\MO_Y(D_Y) \simeq \MO_{\P^{N}_k}(1)|_Y$. 
We have $D \sim \psi^*D_Y$. 
Since $\psi : X \to Y$ is a surjective $k$-morphism,
we obtain 
\[
k \hookrightarrow H^0(Y, \MO_Y) \hookrightarrow H^0(X, \MO_X)=k, 
\]
which implies $H^0(Y, \MO_Y)=k$. 
In particular, we get 
$0 \leq \Delta(Y, D_Y) = \dim Y + \deg D_Y - h^0(Y, D_Y)$ (Lemma \ref{l-Delta-nonneg}). 
Then we obtain 
\begin{eqnarray*}
0 & \leq & \Delta(Y, D_Y)\\
&=& \dim Y + \deg D_Y - h^0(Y, D_Y)\\
&=& 1 + \frac{\deg D}{\deg \psi } - h^0(X, D)\\
&=& 1 + \frac{\deg D}{\deg \psi } - \deg D\\
&=& 1 + \deg D \left( \frac{1}{\deg \psi } - 1\right)\\
&\leq & 1 + 3 \left( \frac{1}{\deg \psi } - 1\right).
\end{eqnarray*}
This implies that $\deg \psi \leq \frac{3}{2}$, i.e., $\psi$ is birational. 
\end{proof}


\begin{lem}\label{l-im-codim1-prime}
We use Notation \ref{n-g1-deg3-birat}. 
Assume that $k$ is an infinite field. 
Let $H$ be a general hyperplane on $\P^N_k$. 
Then, for the induced $k$-linear map 
\[
\alpha : H^0(\P^N_k, \MO_{\P^N_k}(1)) \to H^0(Y \cap H, \MO_{\P^N_k}(1)|_{Y \cap H}), 
\]
it holds that 
\[
\dim_k (\Im\,\alpha) \geq h^0(Y \cap H, \MO_{\P^N_k}(1)|_{Y \cap H})-1. 
\]
\end{lem}

\begin{proof}
We have the following commutative diagram consisting of the induced maps: 
\[
\begin{tikzcd}
H^0(\P^N_k, \MO_{\P^N_k}(1)) \arrow[r, "\simeq"] \arrow[rr, bend left=15, "\alpha"]& 
H^0(Y, \MO_{\P^N_k}(1)|_Y) \arrow[r]   \arrow[d, "\simeq"] 
& H^0(Y \cap H, \MO_{\P^N_k}(1)|_{Y \cap H})) \arrow[d, "\simeq"]\\
& 
H^0(X, \MO_X(D)) \arrow[r, "\beta"] 
& H^0(\psi^{-1}(Y \cap H), \MO_X(D)|_{\psi^{-1}(Y \cap H)})). 
\end{tikzcd}
\]
Note that we have $\psi^{-1}(Y \cap H) \xrightarrow{\simeq} Y \cap H$, because 
$\psi : X \to Y$ is birational (Lemma \ref{l-g1-deg3-birat}) and 
$H$ is chosen to be a general hyperplane, so that $Y \cap H$ is contained in the isomorphic locus of $\psi : X \to Y$. 
{\cred Here we used the fact that 
if $k$ is infinite, then a general hyperplane $H \subset \P^N_k$ can avoid given finitely many closed points on $\P^N_k$.} 
Since $D' := \psi^{-1}(Y \cap H)$ is a member of $|D|$, we obtain an exact sequence: 
\[
H^0(X, \MO_X(D)) \xrightarrow{\beta} 
H^0(D', \MO_X(D)|_{D'}) \to H^1(X, \MO_X(D -D')). 
\]
Hence the assertion follows from 
$\dim_k H^1(X, \MO_X(D-D')) = \dim_k H^1(X, \MO_X) =1$. 
\end{proof}

\begin{thm}\label{t-g1-va-prime}
We use Notation \ref{n-g1-deg3-birat}. 
Then $\bigoplus_{m=0}^{\infty} H^0(X, \MO_X(mD))$ is generated by 
$H^0(X, \MO_X(D))$ as a $k$-algebra. 
In particular, $|D|$ is very ample. 
\end{thm}

\begin{proof}
We may assume that $k$ is separably closed. 
In particular, $k$ is an infinite field. 
If $X$ is geometrically integral, then there is nothing to show \cite[Proposition 11.11(3)]{Tan21}. 
Hence we may assume that $X$ is not geometrically integral. 
Then any regular closed point $P$ of $X$ is not $k$-rational \cite[Corollary 2.14]{Tan21}, 
i.e., $\deg P \geq 2$. 
In particular, $|P|$ is base point free (Theorem \ref{t-g1-bpf}). 
Fix a general hyperplane $H$ on $\P^N_k$, 
so that $Y$ is regular around $Y \cap H$ and $\psi : X \to Y$ is isomorphic around $Y \cap H$ 
(Lemma \ref{l-g1-deg3-birat}). 



Set $\MO_Y(\ell) := \MO_{\P^N_k}(\ell)|_Y$ for every $\ell \in \Z$. 
It is enough to prove (A)-(E) below. 
\begin{enumerate}
\item[(A)] $H^0(\P^N_k, \MO_{\P^N_k}(n)) 
\to H^0(Y \cap H, \MO_{\P^N_k}(n)|_{Y \cap H})$ is surjective for every $n \geq 2$. 
\item[(B)] $H^1(Y, \MO_Y(m))=0$ for every $m \in \Z_{>0}$. 
\item[(C)] $\psi: X \to Y$ is an isomorphism. 
\item[(D)] $\alpha_n : H^0(\P^N, \MO_{\P^N}(n))  \to H^0(Y, \MO_{Y}(n))$ is surjective for every $n \geq 0$. 
\item[(E)] $H^0(Y, \MO_{Y}(m-1)) \otimes_k H^0(Y, \MO_{Y}(1)) \to H^0(Y, \MO_{Y}(m))$ is surjective for every $m \geq 1$. 
\end{enumerate}
Indeed, 
{\cred (C) and (E) imply} 
that $\bigoplus_{m=0}^{\infty} H^0(X, \MO_X(mD))$ is generated by $H^0(X, \MO_X(D))$. 

\medskip

Let us show (A). 
By Lemma \ref{l-im-codim1-prime} and Lemma \ref{l-key}, 
\[
H^0(\P^N_k, \MO_{\P^N_k}(2)) 
\to H^0(Y \cap H, \MO_{\P^N_k}(2)|_{Y \cap H})
\]
is surjective. 
In particular, 
\[
H^0(\P^N_k, \MO_{\P^N_k}(n)) \to H^0(Y \cap H, \MO_Y(n)|_{Y \cap H}) 
\]
is surjective for every $n \geq 2$. 
Thus (A) holds.

Let us show (B). 
For an integer $n \geq 2$, we have  the exact sequence: 
\[
H^0(Y, \MO_Y(n)) \to H^0(Y \cap H, \MO_Y(n)|_{Y \cap H}) \to 
H^1(Y, \MO_Y(n-1)) \to H^1(Y, \MO_Y(n)) \to 0. 
\]
It follows from (A) and the Serre vanishing theorem that $H^1(Y, \MO_Y(m))=0$ for every $m \in \Z_{>0}$. 
Thus (B) holds.

Let us show (C). 
By (B), we have  
\begin{equation}\label{e01-g1-va1}
\chi(Y, \MO_Y(1)) = h^0(Y, \MO_Y(1)) = h^0(X, \MO_X(D)) = \chi(X, \MO_X(D)). 
\end{equation}
For the closed subschemes $C_X \subset X$ and $C_Y \subset Y$ defined by the conductor of $\psi: X \to Y$, 
we have the following conductor exact sequence (\ref{e3-cond}): 
\begin{equation}\label{e02-g1-va1}
0 \to \MO_Y \to \psi_* \MO_X \oplus \MO_{C_Y} \to \psi_*\MO_{C_X} \to 0. 
\end{equation}
We then obtain 
\[
\chi(Y, \MO_Y) - \chi(X, \MO_X) = \chi(C_Y, \MO_{C_Y}) - \chi(C_X, \MO_{C_X}) 
\overset{{\rm (i)}}{=} \chi(Y, \MO_Y(1)) -\chi(X, \MO_X(D))\overset{{\rm (ii)}}{=} 0, 
\]
where (ii) follows from (\ref{e01-g1-va1}) and (i) is obtained by applying the  tensor 
product $(-) \otimes \MO_Y(1)$ to the exact sequence (\ref{e02-g1-va1}). 
By $H^0(X, \MO_X)=H^0(Y, \MO_Y) =k$, we get $h^1(X, \MO_X)=h^1(Y, \MO_Y)$. 
Then (\ref{e02-g1-va1}) induces the following exact sequence: 
\[
0 \to  H^0(Y, \MO_Y) \to H^0(X, \MO_X) \oplus H^0(C_Y, \MO_{C_Y}) \to H^0(C_X, \MO_{C_X}) \to 0.
\]
Again by $H^0(X, \MO_X)=H^0(Y, \MO_Y) =k$, 
we get the induced $k$-linear isomorphism: 
\[
 H^0(C_Y, \MO_{C_Y}) \xrightarrow{\simeq} H^0(C_X, \MO_{C_X}), 
\]
which implies $\MO_{C_Y} \xrightarrow{\simeq} \psi_*\MO_{C_X}$. 
Again by the conductor exact sequence (\ref{e02-g1-va1}), 
we obtain $\MO_Y \xrightarrow{\simeq} \psi_*\MO_X$, i.e., $\psi$ is an isomorphism. 
Thus (C) holds.

\medskip

Let us show (D). 
If $n \in \{0, 1\}$, then $\alpha_n$ is surjective by construction. 
Fix $n \in \Z_{\geq 2}$. 
We have the following commutative diagram in which each horizontal sequence is exact 
(note that $H^1(\MO_Y(n-1))=0$  by (B)): 
\[
\begin{CD}
0 @>>> H^0(\MO_{\P^{N}}(n-1)) @>>>  H^0(\MO_{\P^{N}}(n)) @>>>  H^0(\MO_H(n)) @>>> 0\\
@. @VV\alpha_{n-1}V @VV\alpha_n V @VV\gamma_n V\\
0 @>>>  H^0(\MO_{Y}(n-1)) @>>>  H^0(\MO_{Y}(n)) @>\rho_n >>  H^0(\MO_{Y \cap H}(n)) @>>> 0. 
\end{CD}
\]
By $n \geq 2$, the restriction map 
\[
H^0(\MO_{\P^{N}}(n)) \to H^0(\MO_{Y \cap H}(n))
\]
is surjective by (A). In particular, $\gamma_n$ is surjective. 
By the snake lemma, the surjectivity of $\alpha_{n-1}$ implies 
the surjectivity of $\alpha_n$. 
By induction on $n$, $\alpha_n$ is surjective. 
Thus (D) holds.

Let us show (E).  Fix $m \geq 1$. 
We have the following commutative diagram: 
\[
\begin{CD}
H^0(\MO_{\P^{N}}(m-1)) \otimes_k H^0(\MO_{\P^{N}}(1)) @>\mu>> H^0(\MO_{\P^{N}}(m))\\
@VV\alpha_{m-1} \otimes \alpha_1 V @VV\alpha_m V\\
H^0(\MO_{Y}(m-1)) \otimes_k H^0(\MO_{Y}(1)) @>\mu'>> H^0(\MO_{Y}(m)).\\
\end{CD}
\]
Since $\mu$ and $\alpha_m$ are surjective by (D), also $\mu'$ is surjective. 
Thus (E) holds. 
\qedhere

\end{proof}

\begin{lem}\label{l-key}
Let $k$ be a  field. 
Let $Z$ be a zero-dimensional closed subscheme of $\P^N_k$. 
Assume that the image of the $k$-linear map 
\[
\rho_1 : H^0(\P^N_k, \MO_{\P^N_k}(1)) \to H^0(Z, \MO_{\P^N_k}(1)|_{Z}), 
\]
is at least condimension one, i.e.,  $\dim (\Im\,\rho_1) \geq h^0(Z, \MO_{\P^N_k}(1)|_{Z}) -1$.  
Then the induced $k$-linear map 
\[
\rho_2 : H^0(\P^N_k, \MO_{\P^N_k}(2)) 
\to H^0(Z, \MO_{\P^N_k}(2)|_{Z})
\]
is surjective. 
\end{lem}

\begin{proof}
Taking the base change to the algebraic closure of $k$, 
the problem is reduced to the case when $k$ is algebraically closed. 
We have 
\[
Z = Z_1 \amalg \cdots \amalg Z_m, 
\]
where each $(Z_i)_{\red}$ is one point. 

We introduce a homogeneous coordinate $\P^N_k = \Proj\,k[X_0, ..., X_N]$ 
and fix an affine space: 
\[
\A^N_k := \Spec\,k[x_1, ..., x_N] \subset \P^N_k, \qquad \text{where}\qquad x_i := X_i/X_0. 
\]
We may assume that $Z \subset \A^N_k$. 
Set $R := \MO_Z(Z)$, which is an artinian  ring, equipped with the induced surjective 
$k$-algebra homomorphism  
\[
\pi : k[x_1, ..., x_N] \to R. 
\]
Set $\overline{x}_i := \pi(x_i)$. 
We obtain  the following commutative diagram: 
\[
\begin{tikzcd}
H^0(\P^n_k, \MO_{\P^N_k}) \arrow[r, "\rho_0"] \arrow[d, "\times X_0"] & 
H^0(Z, \MO_Z) \arrow[d, "{\simeq}"', "\times X_0"] \\
H^0(\P^n_k, \MO_{\P^N_k}(1)) \arrow[r, "\rho_1"] \arrow[d, "\times X_0"] & H^0(Z, \MO_{\P^N_k}(1)|_Z) \arrow[d, "{\simeq}"', "\times X_0"] \\
H^0(\P^n_k, \MO_{\P^N_k}(2)) \arrow[r, "\rho_{2}"] & H^0(Z, \MO_{\P^N_k}(2)|_Z).
\end{tikzcd}
\]

We now compute under the affine local setting. 
Via the above right vertical isomorphisms, we identify the following three $k$-vector spaces: 
\[
(R=)H^0(Z, \MO_Z),\qquad H^0(Z, \MO_{\P^N_k}(1)|_Z),\qquad  H^0(Z, \MO_{\P^N_k}(2)|_Z). 
\]
Then $\rho_1$ can be rewritten as  
\[
\rho_1^{\aff} : k \oplus k x_1 \oplus \cdots \oplus k x_N \to R, \qquad 
1 \mapsto 1, \quad x_i \mapsto \overline{x}_i. 
\]
Similarly, $\rho_2$ is given by 
\[
\rho_2^{\aff} : k \oplus k x_1 \oplus \cdots \oplus k x_N \oplus 
\bigoplus_{1 \leq i<j \leq N} k x_ix_j \to R, \qquad 
1 \mapsto 1, \quad x_i \mapsto \overline{x}_i, \quad x_ix_j \mapsto \overline{x}_i\overline{x}_j. 
\]
Since the induced ring homomorphism $\pi: k[x_1, ..., x_N] \to R$ is a surjective $k$-linear map, 
$R$ is generated by the monomials $\{\overline{x}_1^{d_1} \cdots \overline{x}_N^{d_N}\}_{d_1, ..., d_N}$ 
as a $k$-linear space. 
Let $I_1^{\aff}$ and $I_2^{\aff}$ be the images of $\rho_1^{\aff}$ and $\rho_2^{\aff}$, respectively. 
In particular, $I_1^{\aff} \subset I_2^{\aff}$. 
If $I_1^{\aff} = R$, then there is {\cred nothing} to show. 
Assume $I_1^{\aff} \subsetneq R$. 
It suffices to show that $I_1^{\aff} \subsetneq I_2^{\aff}$. 
Since $I_1^{\aff}$ contains $1, \overline{x}_1, ..., \overline{x}_N$, 
we can find a unique integer $d \geq 2$ such that 
\begin{enumerate}
    \item[(i)] $I_1^{\aff}$ contains all the monomials $\overline{x}_1^{e_1} \cdots \overline{x}_N^{e_N}$ satisfying $\sum_i e_i <d$, and 
    \item[(ii)] there exists a monomial $\overline{x}_1^{d_1} \cdots \overline{x}_N^{d_N} \not\in I_1^{\aff}$ with $\sum_i d_i =d$. 
\end{enumerate}
By (ii), it is enough to show that $I_2^{\aff}$ contains $\overline{x}_1^{d_1} \cdots \overline{x}_N^{d_N}$. 
Since we have $d>0$, we may assume, after permuting the indices, that $d_1 >0$. 
By (i), we can write 
\[
\rho_1^{\aff}( a + b_1x_1 + \cdots + b_Nx_N) = 
\overline{x}_1^{d_1-1} \overline{x}_2^{d_2}\cdots \overline{x}_N^{d_N} 
\]
for some $a, b_1, ..., b_N \in k$. 
Then we obtain 
\[
\rho_2^{\aff}( x_1(a + b_1x_1 + \cdots + b_Nx_N)) = 
\overline{x}_1^{d_1} \overline{x}_2^{d_2}\cdots \overline{x}_N^{d_N}, 
\]
as required. 
This completes the proof. 
\end{proof}

\section{Two-dimensional elliptic singularities over imperfect fields}\label{s-ell}

The purpose of this section is to study elliptic singularities defined over an arbitrary field. 
If the base field is an algebraically closed field of characteristic zero, 
then all the results in this section are well known \cite{Wag70}, \cite{Lau77}, \cite[Section 4.4]{KM98}.

\subsection{Minimal resolutions and canonical models}

\begin{dfn}\label{d-ell-sing}
Let $k$ be a field. 
We say that $P \in X$ is an {\em elliptic singularity} (over $k$) if 
\begin{enumerate}
\item $X$ is a normal surface over $k$, 
\item $P$ is a closed point of $X$, 
\item $\MO_{X, P}$ is not regular, 
\item $K_X$ is Cartier, and 
\item $\dim_k R^1f_* \MO_Y=1$ for the minimal resolution $f: Y \to X$ of $P \in X$. 
\end{enumerate}
Recall that the minimal resolution $f: Y \to X$ of $P \in X$ 
is the  birational morphism from a normal surface $Y$ such that 
\begin{itemize}
\item $f$ is an isomorphism over $X \setminus P$ and 
\item there is an open neighbourhood $X'$ of $P \in X$ such that 
$f|_{Y'} : Y' \to X'$ is the minimal resolution of $X'$ for $Y' := f^{-1}(X')$. 
\end{itemize}
\end{dfn}

\begin{lem}\label{l-ell-fund}
Let $k$ be a field and let $P \in X$ be an elliptic singularity over $k$. 
Assume that $P$ is a unique non-regular point of $X$. 
Let $f : Y \to X$ be the minimal resolution of $X$ 
and let $Z$ be the $f$-exceptional effective $\Z$-divisor on $Y$ satisfying $K_Y +Z \sim f^*K_X$. 
Then the following hold. 
\begin{enumerate}
\item $P$ is a $k$-rational point and $\dim_k H^1(Z, \MO_Z)=1$. 
\item If 
$Z'$ is a nonzero effective Cartier divisor on $Y$ satisfying $Z' <Z$, then $H^1(Z', \MO_{Z'})=0$.
\item 
$Z$ is the fundamental cycle of $f$. 
\item 
There is the following commutative commutative diagram in which each horizontal sequence is exact and all the vertical arrows are isomorphisms. 
\[
\begin{tikzcd}
0 \arrow[r] & f_*\MO_Y(-Z) \arrow[r] &f_*\MO_Y \arrow[r] &f_*\MO_Z \arrow[r]  & 0\\
0 \arrow[r] & \m_P \arrow[r]\arrow[u, "\simeq"]  & \MO_X \arrow[r] \arrow[u, equal] & \MO_P \arrow[r] \arrow[u, "\simeq"]  & 0
\end{tikzcd}
\]
In particular, $H^0(Z, \MO_Z) = k$. 
\end{enumerate}
\end{lem}

\begin{proof}
Replacing $X$ by a suitable affine open neighbourhood of $P \in X$, 
we may assume that $X$ is affine and $K_X \sim 0$. 

Let us show (1). 
We have the exact sequence 
\[
0 \to \MO_Y(-Z) \to \MO_Y \to \MO_Z \to 0. 
\]
By the Kawamata--Viehweg vanishing theorem for the birational morphism $f : Y \to X$ \cite[Theorem 10.4]{Kol13}, 
we obtain $R^1f_*\MO_Y(-Z) \simeq R^1f_*\MO_Y(K_Y) =0$. 
Hence 
\[
k \simeq R^1f_*\MO_Y \xrightarrow{\simeq} H^1(Z, \MO_Z). 
\]
Since $H^1(Z, \MO_Z)$ is a $\kappa(P)$-vector space satisfying 
\[
1 = \dim_k H^1(Z, \MO_Z) =  [\kappa(P):k] \dim_{\kappa(P)} H^1(Z, \MO_Z), 
\]
we obtain $[\kappa(P):k]=1$, i.e., $P$ is a $k$-rational point. 
Thus (1) holds.

Let us show (2). 
Take $Z'$ as in the statement of (2). 
For the nonzero effective divisor $Z'' := Z - Z'$, we obtain 
an exact sequence 
\[
0 \to \MO_Y(K_Y+Z') \xrightarrow{\alpha} \MO_Y(K_Y+Z) \xrightarrow{\beta} \MO_Y(K_Y+Z)|_{Z''} \to 0. 
\]
By $\MO_Y(K_Y+Z) \simeq \MO_Y$, the map 
\[
H^0(\beta) : H^0(Y, \MO_Y(K_Y+Z)) \to H^0(Z'', \MO_Y(K_Y+Z)|_{Z''})
\]
is nonzero. Hence $H^0(\alpha) : H^0(\MO_Y(K_Y+Z')) \hookrightarrow H^0(\MO_Y(K_Y+Z))$ is injective but not bijective. We have the following commutative diagram 
in which each horizontal sequence is exact: 
\[
\begin{CD}
0 @>>> \MO_Y(K_Y) @>>> \MO_Y(K_Y+Z') @>>> \MO_Y(K_Y+Z')|_{Z'} (\simeq \omega_{Z'}) @>>>0\\
@. @| @VV\alpha V @VV\gamma V\\
0 @>>> \MO_Y(K_Y) @>>> \MO_Y(K_Y+Z) @>>> \MO_Y(K_Y+Z)|_{Z} (\simeq \omega_Z) @>>>0.
\end{CD}
\]
By $R^1f_*\MO_Y(K_Y)=0$ and the snake lemma, 
$H^0(\gamma): H^0(Z', \omega_{Z'}) \to H^0(Z, \omega_Z)$ is 
injective but not bijective. 
It follows from Serre duality that $H^1(Z, \MO_Z) \to H^1(Z', \MO_{Z'})$ is surjective but not bijective. 
By (1), we get $H^1(Z', \MO_{Z'})=0$. 
Thus (2) holds.

Let us show (3). 
By $K_Y + Z \sim f^*K_X$, $-Z$ is $f$-nef. 
Hence we get $Z_f \leq Z$. 
Suppose $Z_f <Z$. Then (2) implies $H^1(Z_f, \MO_{Z_f})=0$. 
However, this leads to a contradiction: $R^1f_*\MO_Y = 0$ (Proposition \ref{p-Artin-rat}). 
Thus (3) holds. 

Let us show (4). 
By $R^1f_*\MO_Y(-Z)=0$, 
there is the following commutative commutative diagram in which each horizontal sequence is exact: 
\[
\begin{tikzcd}
0 \arrow[r] & f_*\MO_Y(-Z) \arrow[r] &f_*\MO_Y \arrow[r] &f_*\MO_Z \arrow[r]  & 0\\
0 \arrow[r] & \m_P \arrow[r]\arrow[u, "\zeta"]  & \MO_X \arrow[r] \arrow[u, equal] & \MO_P \arrow[r] \arrow[u, "\xi"]  & 0,
\end{tikzcd}
\]
where $\zeta$ is induced by 
$H^0(X, \m_P) \subset H^0(Y, \m_P \MO_Y) \subset H^0(Y, \MO_Y(-Z))$ 
(Proposition \ref{p-fund-mult}(1)) and $\zeta$ induces $\xi$. 
By diagram chase, $\xi$ is surjective. 
Since $\xi : \MO_P \to f_*\MO_Z$ is a ring homomorphism from a field $\MO_P$, 
$\xi$ is injective. Thus $\xi$ is bijective. 
By the snake lemma, also $\zeta$ is bijective. 
Thus (4) holds. 
\end{proof}


\begin{nota}\label{n-cano-model}
Let $k$ be a field and let $P \in X$ be an elliptic singularity over $k$. 
Assume that $X$ is affine, $K_X \sim 0$, and $X \setminus P$ is regular. 
Let $f : Y \to X$ be the minimal resolution of $X$. 
Take the canonical model $\overline Y$ of $Y$ over $X$: 
\[
f : Y \xrightarrow{g} \overline Y \xrightarrow{\overline f} X. 
\]
More explicitly, $g: Y \to \overline Y$ is the projective birational morphism to a normal surface such that $\Ex(g)$ consists of all the $f$-exceptional prime divisors $C$ with $K_Y \cdot C=0$. 
Note that the existence of $g$ is guaranteed by \cite[Theorem 4.2] {Tan18m}. 
In this case, both $\overline Y$ and $\overline f : \overline Y \to X$ are called {\em the canonical model} (over $X$). 
Set $Z :=Z_f$ and $\overline Z := g_*Z$. 
\end{nota}

\begin{prop}\label{p-cano-model}
We use Notation \ref{n-cano-model}. 
Then the following hold. 
\begin{enumerate}
    \item $\overline Z$ is a nonzero effective Cartier divisor on $\overline Y$ such that $K_{\overline Y} +\overline Z \sim 0$ and $Z = g^*\overline Z$. 
    \item There exists a morphism $g|_{Z} : Z \to \overline Z$ such that the following diagram is cartesian: 
    \begin{equation}\label{e1-cano-model}
    \begin{tikzcd}
    Z \arrow[r, hook] \arrow[d, "g|_Z"] & Y \arrow[d, "g"]\\
    \overline Z \arrow[r, hook] & \overline Y,
\end{tikzcd}
    \end{equation}
    where each horisontal arrow is the induced closed immersion. 
\item $(g|_Z)_*\MO_Z = \MO_{\overline Z}$. 
\end{enumerate}
\end{prop}

\begin{proof}
Let us show  (1). 
If $\overline {Z}=0$, then $Y$ would be canonical, which contradicts $R^1f_*\MO_Y \neq 0$ \cite[Theorem 10.4]{Kol13}. 
Hence $\overline Z \neq 0$. 
By $K_Y + Z \sim 0$ and $\overline Z =g_*Z$, we obtain $K_{\overline Y} + \overline Z =g_*(K_Y + Z) \sim 0$. 
Since $K_{\overline Y}$ is Cartier, $\overline Z$ is an effective Cartier divisor. 
Moreover, we obtain $\overline Z = g^*Z$ by the negativity lemma \cite[Lemma 2.11(1)]{Tan18m} and 
\[
K_{\overline Y} + \overline Z  \sim 0 \sim g^*(K_Y + Z)  \sim K_{\overline Y} + g^*Z. 
\]
Thus (1) holds. 

Let us show (2).
We have $g^*\overline Z = Z$ by (1). 
This means that $Z$ is the scheme-theoretic inverse image of $Z'$, 
so that the diagram (\ref{e1-cano-model}) is cartesian. 
Thus (2) holds. 

Let us show (3). 
By $g^*\overline Z = Z$, we get $g_*\MO_Y(-Z) = \MO_{\overline Y}(-\overline Z)$. 
Therefore, we obtain the following commutative diagram in which each horisontal sequence is exact:  
\[
\begin{CD}
0 @>>> g_*\MO_Y(-Z) @>>> g_*\MO_Y @>>> (g|_Z)_*\MO_Z @>>> R^1g_*\MO_Y(-Z) =0\\
@. @| @| @AAA\\
0 @>>> \MO_Y(-\overline Z) @>>> \MO_{\overline Y} @>>> \MO_{\overline Z} @>>> 0,
\end{CD}
\]
where $R^1g_*\MO_Y(-Z) \simeq R^1g_*\MO_Y(K_Y)=0$. 
Therefore, the snake lemma implies $(g|_Z)_*\MO_Z = \MO_{\overline Z}$. 
Thus (3) holds. 
\end{proof}

\subsection{Base point freeness}





\begin{prop}\label{p-bpf-smaller}
We use Notation \ref{n-cano-model}. 
Let $L$ be a nef invertible sheaf on $Z$ such that $L \cdot Z >0$. 
Then the following hold. 
\begin{enumerate}
\item Let $C$ be a curve on $Z$ (i.e., $C$ is a prime divisor on $Y$ with $C \subset \Supp\,Z$) such that $L \cdot C>0$. 
For $Z' := Z -C$, the restriction map 
\[
H^0(Z, L) \to H^0(Z', L|_{Z'})
\]
is surjective. 
    \item $H^1(Z, L^{\otimes n})=0$ for every $n \in \Z_{>0}$. 
\end{enumerate}
\end{prop}

\begin{proof}
Let us show (1). 
We have 
\[
0 \to L(-Z')|_C \to L \to L|_{Z'} \to 0, 
\]
where $L(-Z') := \MO_Y(-Z')|_Z \otimes L$. 
It suffices to show $H^1(C, L(-Z')|_C)=0$. 
By 
\[
Z' \cdot C = ( (K_Y+Z) -(K_Y+C)) \cdot C = -\deg \omega_C, 
\]
we obtain 
\[
\deg (L(-Z')|_C) = L \cdot C -Z' \cdot C > \deg \omega_C. 
\]
Hence $H^1(C, L(-Z')|_C)=0$ holds by Serre duality. Thus (1) holds. 

Let us show (2). 
To this end, it suffices to show 
$H^1(C, L^{\otimes n}(-Z')|_C)=0$ and 
$H^1(Z', L^{\otimes n}|_{Z'})=0$. 
The former one follows from 
\[
\deg (L^{\otimes n}(-Z')|_C) \geq \deg (L(-Z')|_C) > \deg \omega_C. 
\]
The latter one holds by  Proposition \ref{p-g0-bpf} and Lemma \ref{l-ell-fund}. Thus (2) holds. 
\end{proof}

\begin{lem}\label{l-Z'-H0}
We use Notation \ref{n-cano-model}. 
Let $C$ be a prime divisor on $Y$ such that $C \subset \Supp\,Z$ and $Z' := Z - C \neq 0$. 
Then 
$h^0(Z', \MO_{Z'}) = h^0(C, \MO_C)$. 
\end{lem}

\begin{proof}
By the exact sequence 
\[
0 \to \MO_Y(-Z')|_C \to \MO_Z \to \MO_{Z'} \to 0, 
\]
we obtain the following exact sequence: 
\[
H^0(C, \MO_Y(-Z')|_C) \to H^0(Z, \MO_Z) \xrightarrow{\rho} H^0(Z', \MO_{Z'})
\]
\[
\to H^1(C, \MO_Y(-Z')|_C) \to H^1(Z, \MO_Z) \to H^1(Z', \MO_{Z'}). 
\]
We have 
\[
\MO_Y(-Z')|_C \simeq \MO_Y(-Z +C)|_C \simeq \MO_Y(K_Y +C)|_C \simeq \omega_C.  
\]
By Serre duality, we get $h^1(C, \MO_Y(-Z')|_C) = h^0(C, \MO_C)$. 
Recall that $h^0(Z, \MO_Z)= h^1(Z, \MO_Z)=1$  (Lemma \ref{l-ell-fund}(1)(4)). 
Since $\rho$ is nonzero, $h^0(Z, \MO_Z)=1$ implies that 
$\rho$ is injective and $H^0(C, \MO_Y(-Z')|_C)=0$. 
We are done by  $H^1(Z', \MO_{Z'})=0$ (Lemma \ref{l-ell-fund}). 
\end{proof}

\begin{thm}\label{t-bpf}
We use Notation \ref{n-cano-model}. Then the following hold. 
\begin{enumerate}
\item If $M$ is an ample invertible sheaf on $\overline Z$ with $\deg M \geq 2$, then $M$ is globally generated. 
\item If $-Z^2 \geq 2$, then $\MO_Y(-Z)|_Z$ and $\MO_{\overline Y}(- \overline Z)|_{\overline Z}$ 
are globally generated.
\end{enumerate}
\end{thm}

\begin{proof}
Let us show (1). 
If $\overline Z$ is a prime divisor, 
then the assertion follows from Theorem \ref{t-g1-bpf}. 
Therefore, the problem is reduced to the case when 
$Z = C_1 + \cdots +C_r$, 
where $r \geq 2$, each $C_i$ is a prime divisor, $g^*M \cdot C_1 >0$, and $g^*{\cred M} \cdot C_2 >0$ (we possibly have $C_i = C_j$ even if $i \neq j$). 
Set $Z_i := Z - C_i$ for each $i \in \{1, 2\}$. 
By Proposition \ref{p-g0-bpf}(1) and Proposition \ref{p-bpf-smaller}(1), 
we get $\Bs\,|g^*M| \cap Z_i  = \emptyset$ for each $i \in \{1, 2\}$. 
Hence $g^*M$ is globally generated. 
By $(g|_Z)_*\MO_Z = \MO_{\overline Z}$ (Proposition \ref{p-cano-model}), 
also $M$ is globally generated. 
Thus (1) holds.

Let us show (2). 
Assume $-Z^2 \geq 2$. 
By (1), $\MO_{\overline Y}(- \overline Z)|_{\overline Z}$ is globally generated. 
Then its pullback $\MO_Y(-Z)|_Z$ is globally generated. 
Thus (2) holds. 
\end{proof}

For later usage, we now treat the case when $-Z^2 =1$. 

\begin{prop}\label{p-bpf-1}
We use Notation \ref{n-cano-model}. 
Assume that $-Z^2 = 1$. 
Then the base scheme $\Bs\,|\MO_Y(-Z)|_Z|$ is scheme-theoretically equal to a $k$-rational point $Q$ on $Z$. 
Furthermore, $Z$ is regular around $Q$. 
\end{prop}


\begin{proof}
Fix a prime divisor $C$ on $Y$ with $-Z \cdot C >0$, whose existence is guaranteed by $Z^2 <0$. 
We can uniquely write $ Z = aC +Z'$ for some $a \in \Z_{>0}$ and the effective $\Z$-divisor $Z'$ on $Y$ satisfying $C \not\subset \Supp\,Z'$. 
Since $-Z$ is $f$-nef, we get 
\[
1 = -Z^2 = -Z \cdot (aC +Z') \geq a (-Z) \cdot C \geq a. 
\]
Therefore, we get $a=1, -Z \cdot C = 1$,  and $-Z \cdot Z'=0$. 
Note that 
\[
H^0(Z, \MO_Y(-Z)|_Z) \to H^0(Z', \MO_Y(-Z)|_{Z'})
\]
is surjective (Proposition \ref{p-bpf-smaller}). 
Since $\MO_Y(-Z)|_{Z'}$ is globally generated (Proposition \ref{p-g0-bpf}), 
the base locus $\Bs\,|\MO_Y(-Z)|_Z|$ is contained in $Z \setminus Z'$. 
As $Z$ is connected, the restriction map  
\[
H^0(Z, \MO_Y(-Z)|_Z) \to H^0(C, \MO_Y(-Z)|_C)
\]
is nonzero. 
Pick a nonzero element $s \in H^0(Z, \MO_Y(-Z)|_Z)$. 
We have $h^0(Z, \MO_Y(-Z)|_Z)= 
-Z^2 +\chi(Z, \MO_Z)=1$, 
where the first equality follows from the Riemann--Roch theorem 
and the second one holds by $h^0(Z, \MO_Z) = h^1(Z, \MO_Z)=1$ 
(Lemma \ref{l-ell-fund}(1)(4)).  
By $h^0(Z, \MO_Y(-Z)|_Z)=1$ and $\Bs\,|\MO_Y(-Z)|_Z| \subset Z \setminus Z'$, 
the closed subscheme $V(s)$ of $Z$ defined by $s$ is contained in $Z \setminus Z'$. 
In particular, $V(s)$ is an effective Cartier divisor on $Z$. 
Furthermore, $V(s) = Q$ for some $k$-rational point $Q \in Z \setminus Z'$, 
because the restriction $V(s)|_C$ to $C$, which is the closed subscheme defined by $s|_C$, 
is an effective Cartier divisor on $C$ with $\deg (V(s)|_C) = \deg (\MO_Y(-Z)|_C)=-Z \cdot C =1$. 
Since $V(s) =Q$ is a regular effective Cartier divisor on $Z$, 
$Z$ is regular around $Q$. 
\qedhere





\end{proof}

\subsection{Birationality}

\begin{nota}\label{n-cano-model2}
Let $k$ be a field and let $P \in X$ be an elliptic singularity over $k$. 
Assume that $X$ is affine, $K_X\sim 0$, 
and $X \setminus P$ is regular. 
Let $f : Y \to X$ be the minimal resolution of $X$ and let $Z \subset Y$ be the fundamental cycle of $P \in X$. 
Take the canonical model $\overline Y$ of $Y$ over $X$ (cf. Notation \ref{n-cano-model}): 
\[
f : Y \xrightarrow{g} \overline Y \xrightarrow{\overline f} X. 
\]
Set $\overline Z := g_*Z$. 
Assume that $-Z^2 \geq 3$, so that $L_Z := \MO_Y(-Z)|_Z$ and 
$L_{\overline Z} := \MO_{\overline{Y}}(-\overline{Z})|_{\overline Z}$ are globally generated (Theorem \ref{t-bpf}). 
Set $N :=h^0(Z, L_Z)-1 =h^0(\overline Z, L_{\overline Z})-1$ and 
let $W$ be the scheme-theoretic image of $\varphi_{|L_Z|} : Z \to \P_k^N$, which coincides with the scheme-theoretic image of 
$\varphi_{|L_{\overline Z}|} : \overline Z \to \P^N_k$ by $(g|_Z)_*\MO_Z = \MO_{\overline Z}$. 
Let 
\[
\psi : Z \xrightarrow{g|_Z} \overline Z \xrightarrow{\overline \psi} W
\]
be the induced morphisms. For $L_W := \MO_{\P^N}(1)|_W$, the following isomorphisms hold: 
$L_{\overline Z} \simeq \overline \psi^*L_W$ and 
$ L_Z \simeq (g|_Z)^*L_{\overline Z} \simeq \psi^*L_W$. 
\end{nota}

\begin{prop}\label{p-birat}
We use Notation \ref{n-cano-model2}. 
Then $\overline{\psi} : \overline Z \to W$ is birational, i.e., there exists an open dense subset $W^{\circ}$ of $W$ 
such that 
$\overline{\psi}^{-1}(W^{\circ})$ is an open dense subset of $\overline Z$ and the induced morphism 
$\overline \psi|_{\overline \psi^{-1}(W^{\circ})} : \overline \psi^{-1}(W^{\circ}) \to W^{\circ}$ an isomorphism. 
\end{prop}

\setcounter{step}{0}
\begin{proof}
 Note that $\psi$ is surjective.

\begin{step}\label{s01-birat}
The assertion of Proposition \ref{p-birat} holds if 
$\overline Z$ is a prime divisor on $\overline Y$. 
\end{step}

\begin{proof}[Proof of Step \ref{s01-birat}]
In this case, 
the assertion follows from Theorem \ref{t-g1-va-prime}. 
\end{proof}

In what follows, we assume that $\overline Z$ is not a prime divisor. 
Fix a prime divisor $C$ on $Y$ such that $L_Z \cdot C = -Z \cdot C >0$. 
In particular, $C \subset \Supp\,Z$. 
Set $Z' := Z - C$ and let $W' := \psi(Z')$ be the scheme-theoretic image of $W'$ by $\psi$. 
We have the following commutative diagram in which all the horizontal arrows are closed immersions: 
\[
\begin{tikzcd}[column sep=huge, row sep=large]
Y & Z \arrow[l, hook'] 
\arrow[d, "\psi"] & Z' \arrow[l, hook', "j_Z"'] \arrow[d, "\psi'"] \arrow[ll, bend right, hook']\\
& W& W'. \arrow[l, hook', "j_W"']
\end{tikzcd}
\]
Set $L_{W'} := j_W^*L_W$ and $L_{Z'} := \psi'^*L_{W'} = j_Z^*L_Z$. 
We have the decomposition $Z' = Z'_{>0} + Z'_{=0}$ such that 
\begin{itemize}
    \item $Z'_{>0}$ and $Z'_{=0}$ are effective $\Z$-divisors on $Y$, 
    \item $L_Z \cdot \Gamma >0$ for every prime divisor $\Gamma \subset \Supp\,Z'_{>0}$, and 
    \item $L_Z \cdot \Gamma =0$ for every prime divisor $\Gamma \subset \Supp\,Z'_{=0}$. 
\end{itemize}
We have $Z'_{>0} \neq 0$, because $\overline Z$ is not a prime divisor.

\begin{step}\label{s02-birat}
The following hold. 
\begin{enumerate}
    \item $\MO_W \hookrightarrow \psi_*\MO_Z$, i.e., the induced homomorphism 
    $\MO_W \to \psi_*\MO_Z$ is injective. 
    \item $H^0(Z, \MO_Z) = H^0(W, \MO_W)=k$. 
    \item $\psi'_*\MO_{Z'} = \MO_{W'}$. 
    \item $H^1(Z', \MO_{Z'})=H^1(W', \MO_{W'})=0$. 
    \item 
    The composite morphism 
    \[
    Z'_{>0} \hookrightarrow Z' \to W'
    \]
    is birational. 
\end{enumerate}
\end{step}

\begin{proof}[Proof of Step \ref{s02-birat}]
The assertion (1) follows from the fact that    $W$ is a scheme-theoretic image of $Z$. 
Since $\psi : Z\to W$ is a $k$-morphism, 
we obtain ring homomorphisms $k \to H^0(W, \MO_W) \hookrightarrow H^0(Z, \MO_Z) =k$ 
(Lemma \ref{l-ell-fund}(4)), which imply $H^0(Z, \MO_Z) = H^0(W, \MO_W) =k$, i.e., (2) holds. 
The assertion (3) follows from Proposition \ref{p-g0-bpf}(4) 
and the surjectivity of $H^0(Z, L_Z) \to H^0(Z', L_{Z'})$ (Proposition \ref{p-bpf-smaller}). 

Let us show (4). By Proposition \ref{p-g0-bpf}, we get $H^1(Z', \MO_{Z'})=0$. 
Then $H^1(W', \MO_{W'})=0$ follows from the injection $H^1(W', \psi'_*\MO_{Z'}) \hookrightarrow H^1(Z', \MO_{Z'})$, which is guaranteed by the corresponding Leray spectral sequence. Thus (4) holds. 
The assertion (5) follows from (3). 
\end{proof}

\begin{step}\label{s03-birat}
The assertion of Proposition \ref{p-birat} holds if either 
\begin{enumerate}
    \item $\overline Z$ has at least three irreducible components, or 
    \item $\overline Z$ has exactly two irreducible components and $\overline Z$ is not reduced. 
\end{enumerate}
\end{step}

\begin{proof}[Proof of Step \ref{s03-birat}]
Since the proofs are very similar, we only treat the case when (2) holds. 
Let $\overline Z = c_1 \overline C_1 + c_2 \overline C_2$ be the irreducible decomposition, 
where $c_i \in \Z_{>0}$ and $\overline{C}_i$ is a prime divisor. 
Since $\overline Z$ is not reduced, we may assume that $c_1 \geq 2$. 
We can write 
\[
Z = c_1 C_1 + c_2 C_2 + Z_{=0}, 
\]
where each $C_i$ denotes the proper transform of $\overline C_i$ and 
$Z_{=0} := Z -(c_1 C_1 + c_2 C_2)$. 
Note that $g_*Z_{=0}=0$. 

By applying Step \ref{s02-birat}(2) for $C := C_1$ and $Z' := Z-C_1$, 
we see that $W$ has exactly two irreducible components. 
For the irreducible decomposition $W_{\red}= W_1 \cup W_2$, we have $\psi(C_1) = W_1$ and $\psi(C_2) = W_2$ by  permuting $W_1$ and $W_2$ if necessary. 
Furthermore, $\overline{\psi}:\overline Z \to W$ is birational around the generic point $\xi_{\overline{C}_2}$ of $\overline{C}_2$, i.e., there exists an open neighbourhood $U$ of  $\xi_{\overline{C}_2} \in W$ such that $\overline{\psi}|_{\overline{\psi}^{-1}(U)}:\overline{\psi}^{-1}(U) \to U$ is an isomorphism. 
It suffices to show that $\overline{\psi}:\overline Z \to W$ is birational around the generic point of $\overline{C}_1$. 
This holds by applying Step \ref{s02-birat} for $C := C_2$ and $Z' : Z -C_2$. 
This completes the proof of Step \ref{s03-birat}. 
\end{proof}

By Step \ref{s03-birat}, we may assume that one of (I)--(III) holds. 
\begin{enumerate}
    \item[(I)] $\overline Z = \overline{C}_1 + \overline{C}_2$, where $\overline{C}_1$ and 
    $\overline{C}_2$ are distinct prime divisors on $\overline Y$. 
\item[(II)] $\overline Z = 2 \overline C$, 
where $\overline C$ is a prime divisor on $\overline Y$. 
\item[(III)] $\overline Z = n \overline C$, 
where $n \in \Z_{\geq 3}$ and $\overline C$ is a prime divisor on $\overline Y$. 
\end{enumerate}

\begin{step}\label{s04-birat}
If (I) or (II) holds, then the following $(*)$ does not hold. 
\begin{itemize}
\item[$(*)$] $W$ is irreducible and $W$ is generically reduced.  
\end{itemize}
\end{step}

\begin{proof}[Proof of Step \ref{s04-birat}]
Assume that (I) or (II) holds.  
Recall that the prime divisor $C$ on $Y$ has been chosen so that $L_Z \cdot C =-Z \cdot C >0$. 
If (II) holds, then such a prime divisor is unique. 
If (I) holds, then there are exactly two such prime divisors, which are the proper transforms $C_1$ and $C_2$ of $\overline{C}_1$ and $\overline{C}_2$, respectively. 
Possibly after permuting $C_1$ and $C_2$, we may assume that $C = C_1$. 
In order to treat two cases (I) and (II) simultaneously, 
we set $C_1 := C_2 :=C$ for the case when (II) holds. 
We get $Z'_{>0} 
=C_2$, which  is a prime divisor. 

Suppose that $(*)$ holds. Let us derive a contradiction. 
By Step \ref{s02-birat}(5), $W'$ is irreducible and generically reduced. 
Then the closed immersion $j_W : W' \hookrightarrow W$ is birational by $(*)$.

We now we prove that 
\begin{equation}\label{e1-birat}
\deg L_Z = 2\deg L_W = 2\deg L_{W'}. 
\end{equation}
Since  $j_W : W' \hookrightarrow W$ and $C_2 = Z'_{>0} \to W'$ are birational, 
we obtain $\deg L_W = \deg L_{W'}$ and $L_Z \cdot C_2 = \deg L_{W'}$ by 
\cite[Definition 1.15 and Lemma 1.18]{Bad01}. 
By symmetry, we get $L_Z \cdot C_1 = \deg L_{W'}$. 
Then 
\[
\deg L_Z = 
L_Z \cdot C_1 + L_Z \cdot C_2 =2\deg L_{W'}. 
\]
This completes the proof of (\ref{e1-birat}).







We obtain  
\[
h^0(Z, L_Z) \overset{\rm (a)}{=} \chi(Z, L_Z) \overset{\rm (b)}{=} \chi(Z, \MO_Z) + \deg L_Z 
\overset{\rm (c)}{=} \deg L_Z \overset{\rm (d)}{=} 2 \deg L_W. 
\]
Here (a) follows from Proposition \ref{p-bpf-smaller}(2), 
(b) holds by the Riemann--Roch theorem, 
(c) is obtained from Lemma \ref{l-ell-fund}(1)(4), and we get (d) by (\ref{e1-birat}). 
For the time being, we finish the proof of Step \ref{s04-birat} by assuming Claim \ref{cl-birat}. 

\begin{claim}\label{cl-birat}
$H^1(W, \MO_W)=0$ and $H^1(W, L_W)=0$. 
\end{claim}

We obtain  
\[
2 \deg L_W = h^0(Z, L_Z )  \overset{{\rm (i)}}{=}h^0(W, L_W) \overset{{\rm (ii)}}{=}
\chi(W, L_W)
\]
\[
\overset{{\rm (iii)}}{=} \chi(W, \MO_W) + \deg L_W \overset{{\rm (iv)}}{=} 1 + \deg L_W, 
\]
where (i) follows from the fact that the composite map 
$H^0(\P^N, \MO_{\P^n}(1)) \to H^0(W, L_W)  \hookrightarrow H^0(Z, L_Z)$ is bijective (and hence both the maps are bijective), 
(ii) and (iv) hold by Claim \ref{cl-birat}, and (iii) is an obtained by the Riemann--Roch theorem. 
We then get $\deg L_W = 1$, which implies 
$\deg L_Z = 2 \deg L_W = 2$ (\ref{e1-birat}). 
This contradicts the assumption $\deg L_Z =-Z^2 \geq 3$ (Notation \ref{n-cano-model2}). 

In order to complete the proof of Step \ref{s04-birat}, it is enough to prove Claim \ref{cl-birat}. 
Since $H^1(W, \MO_W)=0$ implies $H^1(W, L_W)=0$ (Proposition \ref{p-g0-bpf}), 
let us prove $H^1(W, \MO_W)=0$. 
We have a birational closed immersion $j_W : W' \to W$ of irreducible one-dimensional projective schemes. 
Hence we obtain the following exact sequence 
\[
0 \to I \to \MO_{W} \to \MO_{W'} \to 0
\]
for the coherent ideal sheaf $I$ defining $W'$. 
Then $H^1(W, \MO_{W}) \simeq H^1(W', \MO_{W'}) =0$, 
where the isomorphism follows from $\dim(\Supp\,I)=0$ and the equality holds by Step \ref{s02-birat}(4). 
This completes the proof of Step \ref{s04-birat}. 
\qedhere

\end{proof}

\begin{step}\label{s05-birat}
The assertion of Proposition \ref{p-birat} holds if (I) holds.     
\end{step}

\begin{proof}[Proof of Step \ref{s05-birat}]
In this case, $W$ is reduced, because so is $\overline Z$. 
By $\overline Z = \overline{C}_1 + \overline{C}_2$, $W$ has at most two irreducible components. 
Since $(*)$ does not hold, 
$W$ has exactly two irreducible components. 
Let $W = W_1 \cup W_2$ be the irreducible decomposition, where each $W_i$ is equipped with the reduced scheme structure. 
Possibly after permuting $W_1$ and $W_2$, we may assume that $\overline{\psi}(\overline{C}_i) = W_i$ for each $i \in \{1, 2\}$. 
By Step \ref{s02-birat}(5), each $\overline C_i \to W_i$ is birational. 
Therefore, also $\overline Z \to W$ is birational. 
This completes the proof of Step \ref{s05-birat}. 
\end{proof}

In what follows, we assume that $\overline Z$ is irreducible and non-reduced. 
We can write 
\begin{itemize}
    \item $Z = nC + \widetilde{Z}$ for some $n \in \Z_{\geq 2}$, prime divisor $C$, and effective $\Z$-divisor $\widetilde Z$, 
    \item $\overline Z = n\overline{C}$ for $\overline C := g_*C$. 
\end{itemize}
We set 
\begin{itemize}
    \item $Z' := (n-1)C + \widetilde{Z}$, and  
    \item $W' := \psi(Z')$, which denotes the scheme-theoretic image of $Z'$. 
\end{itemize}
Let $\eta_Z, \eta_{\overline{Z}}, \eta_W$ be the generic points of 
$C, \overline{Z}, W$, respectively. 
Consider the following local rings at the generic points 
\[
A := \MO_{W, \eta_W},\quad 
B := \MO_{Z, \eta_Z} = \MO_{\overline{Z}, \eta_{\overline{Z}}},\quad 
R := \MO_{Y, \eta_Z} = \MO_{\overline{Y}, \eta_{\overline{Z}}}, 
\]
\[
A' := \MO_{W', \eta_{W}},\qquad 
B' := \MO_{Z', \eta_{Z}}. 
\]
Note that $R$ is a discrete valuation ring. 
Pick a generator $t \in R$ of the maximal ideal $\m_R$ of $R$, i.e., 
$\m_R = tR$. 
In particular, $B = R/t^nR$ and $B' = R/ t^{n-1}R$. 
We have the following commutative diagram consisting of the induced ring homomorphisms: 
\[
\begin{tikzcd}[column sep=huge]
R \arrow[r, twoheadrightarrow, "\pi"] \arrow[rr, twoheadrightarrow, bend left, "\pi'"] & B= R/t^nR 
\arrow[r, twoheadrightarrow, "\pi''"] & B' = R/ t^{n-1}R\\
& A \arrow[r, twoheadrightarrow] \arrow[u, hook] & A'.\arrow[u, "\simeq"'] \\
\end{tikzcd}
\]
In order to show that $\psi : \overline Z \to W$ is birational, it suffices to show $A=B$, because both $\overline Z$ and $W$ are irreducible.

\begin{step}\label{s06-birat}
The assertion of Proposition \ref{p-birat} holds if (II) holds, i.e., $n =2$. 
\end{step}

\begin{proof}[Proof of Step \ref{s06-birat}]
Assume  $n=2$. 
Set $\kappa := R/\m_R = R/tR$. 
Then $\kappa = B' \simeq A'$, i.e., $W'$ is generically reduced. 
Then $(A, \m_A)$ is an artinian local ring with $A/\m_A =A' \simeq  \kappa$. 
There exists a coefficient field $\kappa_0 \subset A$ \cite[Theorem 28.3]{Mat86}, 
i.e., $\kappa_0$ is a subring of $A$ such that 
the composite ring homomorphism $\kappa_0 \hookrightarrow  A \to A/\m_A$ is an isomorphism. 
Via the composite inclusion $\kappa_0 \subset A \subset B$, $\kappa_0$ is also a coefficient field of $B$. 
By $B = R/t^2 R$, we get $\dim_{\kappa_0} B = 2$. 
In particular, $\dim_{\kappa_0} A =1$ or $\dim_{\kappa_0} A =2$. 
If $\dim_{\kappa_0} A =1$, then 
the induced closed immersion $j_W : W' \hookrightarrow W$ is birational. 
However, this contradicts the fact that $(*)$ does not hold (Step \ref{s05-birat}). 
Hence we get $\dim_{\kappa_0} A =2$, i.e., 
$A=B$, which implies that $\overline Z \to W$ is birational. 
This completes the proof of Step \ref{s06-birat}. 
\qedhere
\end{proof}

\begin{step}\label{s07-birat}
The assertion of Proposition \ref{p-birat} holds if (III) holds, i.e., $n\geq 3$. 
\end{step}

\begin{proof}[Proof of Step \ref{s07-birat}]
Assume $n \geq 3$. In what follows, we identify $A'$ and $B'$ via the above isomorphism. 
For $r \in R$, we set $\overline{r} := \pi(r) \in B$.  
By the above commutative diagram, the following holds: 
\begin{enumerate}
\item[(\#)] For an element $r \in R$, 
there exists $s \in R$ such that $\overline{r} - \overline{s}\overline{t}^{n-1} \in A$. 
\end{enumerate}
By {\cred (\#) and} $B =R +R\overline{t} + R\overline{t}^2 + \cdots + R \overline{t}^{n-1}$, it suffices to show $(\star)_m$ for every $m \in \Z_{>0}$. 
\begin{enumerate}
\item[$(\star)_m$]  $\overline{r}\overline{t}^m \in A$ for every $r \in R$. 
\end{enumerate}
Take $r \in R$. By (\#), we obtain $\overline{r} =\overline{s}\overline{t}^{n-1} +a$ for some $a \in A$ and $s \in R$. 
Again by (\#), 
we have $\overline{t} = \overline{s'} \cdot \overline{t}^{n-1} +a'$
for some $a' \in A$ and $s' \in R$. 
It follows from $\overline{t}^n=0$ and $n-1 \geq 1$ 
that  $\overline{t}^{\ell} = (\overline{t} - \overline{s'} \cdot \overline{t}^{n-1})^{\ell} =a'^{\ell} \in A$ for every integer $\ell \geq 2$. 
By $n -1 \geq 2$, we get $\overline{t}^{n-1} \in A$, 
which implies 
$\overline{r} \overline{t}^{n-1} = (\overline{s}\overline{t}^{n-1} +a)\overline{t}^{n-1} =a \overline{t}^{n-1} \in A$. 
Thus  $(\star)_{n-1}$ holds. 
Fix $m \in \Z_{>0}$. Let us show $(\star)_m$. 
By $\overline{t}^n=0$ and $n-1 \geq 1$, we obtain 
\[
\overline{r} \overline{t}^m = 
(\overline{s}\overline{t}^{n-1} +a)(\overline{s'}\overline{t}^{n-1} +a')^m
= (\overline{s}\overline{t}^{n-1} +a)
(m\overline{s'}\overline{t}^{n-1}a'^{m-1} +a'^m)
\]
\[
= \overline{s}\overline{t}^{n-1} \cdot a'^m + 
a \cdot m\overline{s'}\overline{t}^{n-1}a'^{m-1} + aa'^m. 
\]
By $(\star)_{n-1}$, we obtain 
$\overline{s}\overline{t}^{n-1} \cdot a'^m \in A$ and 
$a \cdot m\overline{s'}\overline{t}^{n-1}a'^{m-1} \in A$. 
Therefore, we get $\overline{r} \overline{t}^m \in A$, i.e., 
$(\star)_m$ holds. 
This completes the proof of Step \ref{s07-birat}. 
\qedhere




\end{proof}
Step \ref{s03-birat}, Step \ref{s05-birat}, Step \ref{s06-birat}, and Step \ref{s07-birat} complete  the proof of Proposition \ref{p-birat}. 
\end{proof}

\subsection{Very ampleness and projective normality}

\begin{lem}\label{l-deg12-generation}
Let $k$ be a field and let $P \in X$ be an elliptic singularity over $k$. 
Assume that $X$ is affine and $X \setminus P$ is regular. 
Let $f : Y \to X$ be the minimal resolution of $X$ and let $Z \subset Y$ be the fundamental cycle of $P \in X$. 
Set $L_Z := \MO_Y(-Z)|_Z$. 
Then the following hold. 
\begin{enumerate}
\item 
$\bigoplus_{m=0}^{\infty} H^0(Z, L_Z^{\otimes m})$  is generated by $H^0(Z, L_Z) \oplus H^0(Z, L_Z^{\otimes 2}) \oplus H^0(Z, L_Z^{\otimes 3})$ as a $k$-algebra. 
\item 
If $\deg {\cred L_Z} \geq 2$, then 
$\bigoplus_{m=0}^{\infty} H^0(Z, L_Z^{\otimes m})$  is generated by $H^0(Z, L_Z) \oplus H^0(Z, L_Z^{\otimes 2})$ as a $k$-algebra. 
\end{enumerate}
\end{lem}


\begin{proof}
We may assume that $K_X \sim 0$. 
We use Notation \ref{n-cano-model}. 
Set $L_{\overline Z} := \MO_{\overline Y}(-\overline Z)$. 
By $\bigoplus_{m=0}^{\infty} H^0(Z, L_Z^{\otimes m}) = \bigoplus_{m=0}^{\infty} H^0(\overline Z, L_{\overline Z}^{\otimes m})$ 
(Proposition \ref{p-cano-model}), 
it is enough to show 
{\cred the corresponding statements for} 
$\bigoplus_{m=0}^{\infty} H^0(\overline Z, L_{\overline Z}^{\otimes m})$. 

Recall that we have 
\[
H^1(\overline Z, L_{\overline Z}) \hookrightarrow H^1(Z, L_Z)=0, 
\]
where $H^1(Z, L_Z)=0$ holds by Proposition \ref{p-bpf-smaller} 
and the inclusion is obtained by 
$(g|_Z)_*\MO_Z = \MO_{\overline Z}$ (Proposition \ref{p-cano-model}(3)) and $g^*L_Z \simeq L_{\overline Z}$. 
Note that $H^0(Z, \MO_Z) = H^0(\overline{Z}, \MO_{\overline Z})=k$ 
(Lemma \ref{l-ell-fund}(4)). 

Let us  show (2). 
By $\deg L_{\overline Z} \geq 2$, 
$L_{\overline Z}$ is an ample globally generated invertible sheaf 
(Theorem \ref{t-bpf}). 
Then the assertion (2) follows from 
Proposition \ref{p-regularity}. 

Let us show (1). 
In this case, $L_{\overline Z}^{\otimes 2}$ is an ample globally generated invertible sheaf (Theorem \ref{t-bpf}). 
By the same argument as in (2), 
we see that $\bigoplus_{n=0}^{\infty} H^0(\overline Z, L_{\overline Z}^{\otimes 2n})$ 
is generated by $H^0(\overline Z, L_{\overline Z}^{\otimes 2})$ as a $k$-algebra. 
Thus it suffices to show that 
\[
H^0(\overline Z, L_{\overline Z}^{\otimes 2}) \otimes 
H^0(\overline Z, L_{\overline Z}^{\otimes (2n+3)}) \to 
H^0(\overline Z, L_{\overline Z}^{\otimes (2n+5)}) 
\]
is surjective for every $n \geq 0$. 
By $H^1(\overline Z, L_{\overline Z}^{\otimes 3} \otimes L_{\overline Z}^{\otimes -2}) =0$ 
(Proposition \ref{p-bpf-smaller}), 
$L^{\otimes 3}_{\overline Z}$ is $0$-regular with respect to $L^{\otimes 2}_{\overline Z}$, 
which implies the required surjectivity 
\cite[Definition 11.1 and Proposition 11.2]{Tan21} (cf. \cite[Section 5.2]{FGI05}). 
\qedhere



\end{proof}

Although the proofs of Lemma \ref{l-im-codim1} and Theorem \ref{t-g1-va} below are 
almost identical to 
the ones of Lemma \ref{l-im-codim1-prime} and Theorem \ref{t-g1-va-prime}, we include proofs for the sake of completeness.

\begin{lem}\label{l-im-codim1}
We use Notation \ref{n-cano-model2}. 
Assume that $k$ is an infinite field. 
Let $H$ be a general hyperplane on $\P^N_k$. 
Then, for the induced $k$-linear map 
\[
\alpha : H^0(\P^N_k, \MO_{\P^N_k}(1)) \to H^0(W \cap H, \MO_{\P^N_k}(1)|_{W \cap H}), 
\]
it holds that 
 \[ 
\dim_k (\Im\,\alpha) \geq h^0(W \cap H, \MO_{\P^N_k}(1)|_{W \cap H})-1. 
\]
\end{lem}

\begin{proof}
We have the following commutative diagram consisting of the induced maps: 
\[
\begin{tikzcd}
H^0(\P^N_k, \MO_{\P^N_k}(1)) \arrow[r, "\simeq"] \arrow[rr, bend left=15, "\alpha"]& 
H^0(W, \MO_{\P^N_k}(1)|_W) \arrow[r]   \arrow[d, "\simeq"] 
& H^0(W \cap H, \MO_{\P^N_k}(1)|_{W \cap H})) \arrow[d, "\simeq"]\\
& 
H^0(\overline Z, L_{\overline Z}) \arrow[r, "\beta"] 
& H^0(\psi^{-1}(W \cap H), L_{\overline Z}|_{\overline{\psi}^{-1}(W \cap H)})).
\end{tikzcd}
\]
Recall that $\overline \psi :\overline{Z} \to W$ is birational (Proposition \ref{p-birat}). 
Hence we have $\overline{\psi}^{-1}(W \cap H) \xrightarrow{\simeq} W \cap H$, because 
$H$ is chosen to be a general hyperplane, so that $W \cap H$ is contained in the isomorphic locus of $\overline{\psi} : \overline{Z} \to W$. 
Since $D := \overline{\psi}^{-1}(W \cap H)$ is a member of $|L_{\overline Z}|$, we obtain an exact sequence: 
\[
H^0(\overline{Z}, L_{\overline Z}) \xrightarrow{\beta} 
H^0(D, L_{\overline Z}|_{D}) \to H^1(\overline Z, L_{\overline Z} \otimes \MO_{\overline Z}(-D)). 
\]
By $\MO_{\overline Z}(-D) \simeq L_{\overline Z}^{-1}$, the assertion follows from 
\[
\dim_k H^1(\overline Z, L_{\overline Z} \otimes \MO_{\overline Z}(-D)) 
= \dim_k H^1(\overline Z, \MO_{\overline Z}) 
\overset{{\rm (i)}}{\leq} \dim_k H^1(Z, \MO_{Z})\overset{{\rm (ii)}}{\leq}1, 
\]
where (i) and (ii) hold by Proposition \ref{p-cano-model}(3) 
and 
Lemma \ref{l-ell-fund}(1), respectively. 
\end{proof}

\begin{thm}\label{t-g1-va}
We use Notation \ref{n-cano-model2}. 
Then $\bigoplus_{m=0}^{\infty} H^0(Z, L_Z^{\otimes m})$ 
(resp. $\bigoplus_{m=0}^{\infty} H^0(\overline Z, L_{\overline Z}^{\otimes m})$) is generated by 
$H^0(Z, L_{Z})$ (resp. $H^0(\overline Z, L_{\overline Z})$) as a $k$-algebra. 
In particular, $|L_{\overline Z}|$ is very ample. 
\end{thm}

\begin{proof}
By $\bigoplus_{m=0}^{\infty} H^0(Z, L_Z^{\otimes m}) = \bigoplus_{m=0}^{\infty} H^0(\overline Z, L_{\overline Z}^{\otimes m})$  (Proposition \ref{p-cano-model}(3)), 
it is enough to show that 
$\bigoplus_{m=0}^{\infty} H^0(\overline Z, L_{\overline Z}^{\otimes m})$ 
is generated by $H^0(\overline Z, L_{\overline Z})$. 
We may assume that $k$ is separably closed. 
In particular, $k$ is an infinite field. 

Set $\MO_W(\ell) := \MO_{\P^N_k}(\ell)|_W =L_W^{\otimes \ell}$ for every $\ell \in \Z$. 
It is enough to prove the  following assertions (A)--(E). 
\begin{enumerate}
\item[(A)] $H^0(\P^N_k, \MO_{\P^N_k}(n)) 
\to H^0(W \cap H, \MO_{\P^N_k}(n)|_{W \cap H})$ is surjective for every $n \geq 2$. 
\item[(B)] $H^1(W, \MO_W(m))=0$ for every $m \in \Z_{>0}$. 
\item[(C)] $\overline{\psi}: \overline Z \to W$ is an isomorphism. 
\item[(D)] $\alpha_n : H^0(\P^N, \MO_{\P^N}(n))  \to H^0(W, \MO_{W}(n))$ is surjective for every $n \geq 0$. 
\item[(E)] $H^0(W, \MO_{W}(m-1)) \otimes_k H^0(W, \MO_{W}(1)) \to H^0(W, \MO_{W}(m))$ is surjective for every $m \geq 1$. 
\end{enumerate}
Indeed, 
{\cred (C) and (E) imply} 
that $\bigoplus_{m=0}^{\infty} H^0(\overline Z, \MO_{\overline Z}(mD))$ is generated by $H^0(\overline Z, \MO_X(D))$. 

\medskip

Let us show (A). 
By Lemma \ref{l-key} and Lemma \ref{l-im-codim1}, 
\[
H^0(\P^N_k, \MO_{\P^N_k}(2)) 
\to H^0(W \cap H, \MO_{\P^N_k}(2)|_{W \cap H})
\]
is surjective. 
Since $W \cap H$ is zero-dimensional, also  
\[
H^0(W, \MO_W(n)) \to H^0(W \cap H, \MO_W(n)|_{W \cap H}) 
\]
is surjective for every $n \geq 2$. 
Thus (A) holds. 

Let us show (B). 
For an integer $n \geq 2$, we have  the exact sequence: 
\[
H^0(W, \MO_W(n)) \to H^0(W \cap H, \MO_W(n)|_{W \cap H}) \to 
H^1(W, \MO_W(n-1)) \to H^1(W, \MO_W(n)) \to 0. 
\]
By (A), it follows from the Serre vanishing theorem that $H^1(W, \MO_W(m))=0$ for every $m \in \Z_{>0}$. 
Thus (B) holds.

Let us show (C). 
By (B) and $H^1(\overline Z, L_{\overline Z})=0$ 
(Proposition \ref{p-cano-model}(3), Proposition \ref{p-bpf-smaller}(2)), we have  
\begin{equation}\label{e01-g1-va}
\chi(W, \MO_W(1)) = h^0(W, \MO_W(1)) = h^0(\overline Z, 
L_{\overline Z}) = \chi(\overline Z, L_{\overline Z}). 
\end{equation}
For the closed subschemes $C_{\overline Z} \subset \overline Z$ and $C_W \subset W$ defined by the conductor of $\overline{\psi}: \overline Z \to W$, 
we have the following conductor exact sequence (\ref{e3-cond}): 
\begin{equation}\label{e02-g1-va}
0 \to \MO_W \to \overline{\psi}_* \MO_{\overline Z} \oplus \MO_{C_W} \to \overline{\psi}_*\MO_{C_{\overline Z}} \to 0. 
\end{equation}
We then obtain 
\[
\chi(W, \MO_W) - \chi(\overline Z, \MO_{\overline Z}) 
= \chi(C_W, \MO_{C_W})  -\chi(C_{\overline Z}, \MO_{C_{\overline Z}})  
\overset{{\rm (i)}}{=} 
\chi(W, \MO_W(1)) - \chi(\overline Z, L_{\overline Z}) \overset{{\rm (ii)}}{=} 0, 
\]
where (ii) follows from (\ref{e01-g1-va}) and (i) is obtained by 
applying the  tensor $(-) \otimes \MO_W(1)$ to the exact sequence (\ref{e02-g1-va}). 
By $H^0(\overline Z, \MO_{\overline Z})=H^0(W, \MO_W) =k$ 
(Lemma \ref{l-ell-fund}(1), Proposition \ref{p-cano-model}(3)), 
we get $h^1(\overline Z, \MO_{\overline Z})=h^1(W, \MO_W)$. 
Then (\ref{e02-g1-va}) induces the following exact sequence: 
\[
0 \to  H^0(W, \MO_W) \to H^0(\overline Z, \MO_{\overline Z}) \oplus H^0(C_W, \MO_{C_W}) \to H^0(C_{\overline Z}, \MO_{C_{\overline Z}}) \to 0.
\]
Again by $H^0(\overline Z, \MO_{\overline Z})=H^0(W, \MO_W) =k$, 
we get the induced $k$-linear isomorphism: 
\[
 H^0(C_W, \MO_{C_W}) \xrightarrow{\simeq} H^0(C_{\overline Z}, \MO_{C_{\overline Z}}), 
\]
which implies $\MO_{C_W} \xrightarrow{\simeq} \overline{\psi}_*\MO_{C_{\overline Z}}$. 
Again by the conductor exact sequence (\ref{e02-g1-va}), 
we obtain $\MO_W \xrightarrow{\simeq} \overline{\psi}_*\MO_{\overline Z}$, i.e., 
$\overline{\psi} : \overline Z \to W$ is an isomorphism. 
Thus (C) holds.

\medskip

Let us show (D). 
If $n \in \{0, 1\}$, then $\alpha_n$ is surjective by construction. 
Fix $n \in \Z_{\geq 2}$. 
We have the following commutative diagram in which each horizontal sequence is exact 
(note that $H^1(W, \MO_W(n-1))=0$  by (B)): 
\[
\begin{CD}
0 @>>> H^0(\MO_{\P^{N}}(n-1)) @>>>  H^0(\MO_{\P^{N}}(n)) @>>>  H^0(\MO_H(n)) @>>> 0\\
@. @VV\alpha_{n-1}V @VV\alpha_n V @VV\gamma_n V\\
0 @>>>  H^0(\MO_{W}(n-1)) @>>>  H^0(\MO_{W}(n)) @>\rho_n >>  H^0(\MO_{W \cap H}(n)) @>>> 0. 
\end{CD}
\]
By $n \geq 2$, the composite map 
\[
H^0(\MO_{\P^{N}}(n)) \to H^0(\MO_{W \cap H}(n))
\]
is surjective by (A). In particular, $\gamma_n$ is surjective. 
By the snake lemma, the surjectivity of $\alpha_{n-1}$ implies 
the surjectivity of $\alpha_n$. 
By induction on $n$, $\alpha_n$ is surjective. 
Thus (D) holds.

\medskip

Let us show (E).  Fix $m \geq 1$. 
We have the following commutative diagram: 
\[
\begin{CD}
H^0(\MO_{\P^{N}}(m-1)) \otimes_k H^0(\MO_{\P^{N}}(1)) @>\mu>> H^0(\MO_{\P^{N}}(m))\\
@VV\alpha_{m-1} \otimes \alpha_1 V @VV\alpha_m V\\
H^0(\MO_{W}(m-1)) \otimes_k H^0(\MO_{W}(1)) @>\mu'>> H^0(\MO_{W}(m)).\\
\end{CD}
\]
Since $\mu$ and $\alpha_m$ are surjective by (D), also $\mu'$ is surjective. 
Thus (E) holds. 
\qedhere
\end{proof}

\subsection{Case study}


\begin{nota}\label{n-ell-cases}
Let $k$ be a field and let $P \in X$ be an elliptic singularity over $k$. 
Assume that $X \setminus P$ is regular, $X$ is affine, and $K_X \sim 0$. 
Let $f : Y \to X$ be the minimal resolution of $X$ and let $Z \subset Y$ be the fundamental cycle of $P \in X$. 
Take the canonical model $\overline Y$ of $Y$ over $X$ (cf. Notation \ref{n-cano-model}): 
\[
f : Y \xrightarrow{g} \overline Y \xrightarrow{\overline f} X. 
\]
Set $\overline Z := g_*Z$, $L :=\MO_Y(-Z), L_{\overline Y} := \MO_{\overline Y}(-\overline Z), L_Z := \MO_Y(-Z)|_Z$, and 
$L_{\overline Z} := \MO_{\overline{Y}}(-\overline{Z})|_{\overline Z}$. 
By $K_Y +Z \sim f^*K_X \sim 0$, we obtain 
\[
\overline{Y} = \Proj\,\bigoplus_{n=0}^{\infty} H^0(Y, nK_Y) \simeq \Proj\,\bigoplus_{n=0}^{\infty} H^0(Y, \MO_Y(-nZ)). 
\]
Let $\m_P$ be the maximal ideal of $H^0(X, \MO_X)$ corresponding to $P$. 
For every $n \in \Z_{\geq 0}$, we set 
\[
I_n :=H^0(Y, \MO_Y(-nZ)) \subset H^0(Y, \MO_Y) = H^0(X, \MO_X), 
\]
which is an ideal of $H^0(X, \MO_X)$. 
By 
\[
0 \to H^0(Y, {\cred \MO_Y}(-(n+1)Z)) \to H^0(Y, {\cred \MO_Y}(-nZ)) \to H^0(Z, \MO_Y(-nZ)|_Z) 
\]
\[
\to 
H^1(Y, {\cred \MO_Y}(-(n+1)Z))=0, 
\]
we obtain $I_n/I_{n+1} \simeq H^0(Z, \MO_Y(-nZ)|_Z) \simeq H^0(Z, L_Z^{\otimes n})$, and hence 
\[
\bigoplus_{n=0}^{\infty} I_n/I_{n+1} \simeq \bigoplus_{n=0}^{\infty} H^0(Z, L_Z^{\otimes n}). 
\]
\end{nota}

\begin{lem}\label{l-I1=m}
We use {\cred Notation} \ref{n-ell-cases}. Then $I_1 =\m_P$. 
\end{lem}

\begin{proof}
We have $\m_P = f_*\MO_Y(-Z)$ (Lemma \ref{l-ell-fund}), which implies $I_1 = \m_P$. 
\end{proof}

\begin{rem}
We use {\cred Notation} \ref{n-ell-cases}. 
We have $I_1 \supset I_2 \supset \cdots$ and $I_1^n \subset I_n$ for every $n \in \Z_{>0}$. 
Hence the inclusions $I_1^n \subset I_n \subset I_1$ imply $\sqrt{I_n} = I_1 = \m$. 
\end{rem}

\begin{lem}\label{l-Artin-Rees}
We use Notation \ref{n-ell-cases}. 
Then there exists $a \in \Z_{>0}$ such that the induced $k$-linear map
\[
H^0(Y, L^{\otimes n}) \otimes_k H^0(Y, L^{\otimes a}) \to H^0(Y, L^{\otimes (n+a)})
\]
is surjective for every $n \in \Z_{\geq 0}$. 
In particular, $I_{n} \cdot I_a = I_{n+a}$ for every $n \in \Z_{\geq 0}$.  
\end{lem}

\begin{proof}
This map is the same as 
\[
H^0(\overline Y, L^{\otimes n}_{\overline Y}) \otimes_k  H^0(\overline Y, L^{\otimes a}_{\overline Y}) \to H^0(\overline Y, L^{\otimes (n+a)}_{\overline Y}). 
\]
Then the assertion follows from the fact that $L_{\overline Y}$ is (relatively) ample. 
\end{proof}

\subsubsection{Case $-Z^2 \geq 3$}

\begin{thm}\label{t-ell-3}
We use Notation \ref{n-ell-cases}. 
Assume $-Z^2 \geq 3$. 
Then $\overline Y \simeq {\rm Bl}_PX$, i.e., 
the canonical model $\overline{f} : \overline{Y} \to X$ coincides with the blowup 
${\rm Bl}_P\,X \to X$ at the maximal ideal $\m_P$ corresponding to the singularity $P$. 
\end{thm}

\begin{proof}
Recall that we have 
\[
\overline Y = 
\Proj\,\bigoplus_{n=0}^{\infty} H^0(Y, \MO_Y(nK_Y)) 
\simeq 
\Proj\,\bigoplus_{n=0}^{\infty} H^0(Y, \MO_Y(-nZ)) =  
\Proj\,\bigoplus_{n=0}^{\infty} I_n. 
\]
By ${\rm Bl}_P\,X =\Proj\,\bigoplus_{n=0}^{\infty} \m^n$, it suffices to show that $I_n = \m^n$ for every $n \geq 0$. 
We have 
\[
\bigoplus_{n=0}^{\infty} I_n/I_{n+1} \simeq \bigoplus_{n=0}^{\infty} H^0(Z, nL_Z).
\]
Since the right hand side is generated by $H^0(Z, L_Z)$ as a $k$-algebra (Theorem \ref{t-g1-va}), 
the left hand side is generated by $I_1/I_2$. 
Hence we get $I_n = I_1^n + I_{n+1}$ for every $n \in \Z_{>0}$. 
By $I_n = I_1^n + I_{n+1}$ and $I_{n+1}= I_1^{n+1} + I_{n+2}$, 
we get 
\[
I_n = I_1^n + I_{n+1}  = I_1^n + (I_1^{n+1} + I_{n+2}) = I_1^n + I_{n+2}. 
\]
Repeating this, it holds that 
\[
I_n = I_1^n + I_{N}
\]
whenever $N >n \geq 1$. 
Recall that there exists $a \in \Z_{>0}$ such that $I_{m+a} = I_m \cdot I_a$ 
for every $m \in \Z_{>0}$ (Lemma \ref{l-Artin-Rees}). 
For $N \gg 0$, we get  
\[
I_N  = I_{N-a}I_a = I_{N-2a} I_a^2 = \cdots = I_{N -na} I_a^n \subset I_a^n \subset I_1^n. 
\]
Therefore, we obtain $I_1^n \subset I_n = I_1^n + I_{N} \subset I_1^n$, as required. 
\end{proof}

\begin{cor}\label{c-ell-3}
We use Notation \ref{n-ell-cases}. 
Assume $-Z^2 \geq 3$. 
Then $\m_P \MO_{\overline Y} = \MO_{\overline Y}(-\overline Z)$, 
{\cred where} 
$\m_P$ denotes the maximal ideal corresponding to $P$. 
\end{cor}

\begin{proof}
By Theorem \ref{t-ell-3},  $\m_P \MO_{\overline Y}$ is an invertible sheaf. 
Thus we can write $\m_P \MO_{\overline Y} = \MO_{\overline Y}(-F)$ for some effective Cartier divisor $F$ on $\overline Y$. In particular, $F$ is $\overline f$-exceptional for $\overline f : \overline Y \to X$. 
We have $\m_P \MO_Y = \MO_Y(-Z)$ by Proposition \ref{p-fund-mult}(5) and  Theorem \ref{t-bpf}. 
Hence $g^*F \sim Z \sim g^*\overline{Z}$. 
Hence $F =g_*g^*F \sim g_*g^*\overline{Z} = \overline Z$. 
Since both $F$ and $\overline Z$ are $\overline f$-exceptional, 
we obtain $F = \overline Z$ by the negativity lemma \cite[Lemma 2.11]{Tan18m}. 
\end{proof}

\subsubsection{Case $-Z^2 \leq 2$}




\begin{lem}\label{l-sym2}
We use Notation \ref{n-ell-cases}. 
Assume $-Z^2 =1$ or $-Z^2 =2$. 
Then the induced $k$-linear map 
\[
\mu_Z : S^2H^0(Z, L_Z) \to H^0(Z, L_Z^{\otimes 2})
\]
is injective, where $S^2(-)$ denotes the symmetric product. 
\end{lem}

\begin{proof}
By $(g|_Z)_*\MO_Z = \MO_{\overline Z}$ (Proposition \ref{p-cano-model}(3)) and $L_Z \simeq g^*L_{\overline Z}$, 
we have the following commutative diagram consisting of the induced $k$-linear maps: 
\[
\begin{tikzcd}
S^2H^0(Z, L_Z) \arrow[r, "\mu_Z"] \arrow[d, "\simeq"] & H^0(Z, L_Z^{\otimes 2}) \arrow[d, "\simeq"]\\
S^2H^0(\overline{Z}, L_{\overline{Z}}) \arrow[r, "\mu_{\overline Z}"] & H^0(\overline{Z}, L^{\otimes 2}_{\overline{Z}}).
\end{tikzcd}
\]
Therefore, if $\overline Z$ is a prime divisor, then $\mu_{\overline Z}$ is injective, and hence 
also $\mu_Z$ is injective. 

In what follows, we treat the case when $\overline Z$ is not a prime divisor. 
In particular, 
we get $L \cdot Z = -Z^2 =2$. 
Furthermore, the following holds. 
\begin{enumerate}
\item[$(*)$] $Z = C+C'+D$, where $L \cdot C = L \cdot C'=1$, $L \cdot D =0$, both $C$ and $C'$ are prime divisors, 
and $D$ is an effective $\Z$-divisor on $Y$ (note that we possibly have $C =C'$). 
\end{enumerate}
Set $Z' := Z - C = C'+D$. 
By $L \cdot C = L \cdot C'=1$, we obtain 
\[
h^0(Z', \MO_{Z'}) =h^0(C, \MO_C)=1, \qquad  h^0(C', \MO_{C'})=1, 
\]
where $h^0(Z', \MO_{Z'}) =h^0(C, \MO_C)$ follows from Lemma \ref{l-Z'-H0}. 
We obtain $H^1(Z', \MO_{Z'})=0$ (Lemma \ref{l-ell-fund}(2)). 
By $h^0(Z', L_{Z'}) =  L \cdot Z' +h^0(Z', \MO_{Z'}) = 2$ 
(Proposition \ref{p-g0-H1}), 
we have the morphism 
\[
\varphi : Z' \xrightarrow{\psi} W \hookrightarrow \P^1_k
\]
induced by $|L|_{Z'}|$, where $W$ denotes its image and $\psi_*\MO_{Z'} = \MO_W$ (Proposition \ref{p-g0-bpf}). 
We then get $W = \P^1_k$, since $\P^1_k$ is a unique one-dimensional closed subscheme of $\P^1_k$. 
In particular, $\varphi_*\MO_{Z'} = \MO_W$. 
For the time being, let us finish the proof by assuming Claim \ref{cl-sym2}. 

\begin{claim}\label{cl-sym2}
The restriction map $\rho: H^0(Z, L_Z) \to H^0(Z', L_Z|_{Z'})$ is an isomorphism. 
\end{claim}

For $L_W := \MO_{\P^1}(1)$, we obtain  the following commutative diagram: 
\[
\begin{tikzcd}[column sep=large, row sep=large]
S^2H^0(Z, L_Z) \arrow[r, "{S^2\rho,\,\,\simeq}"] \arrow[d, "\mu_Z"] & S^2H^0(Z', L_Z|_{Z'}) \arrow[d, "\mu_{Z'}"] \arrow[r, "\simeq"] & S^2H^0(W, L_W) \arrow[d, hook, "\mu_W"]\\
H^0(Z, L^{\otimes 2}_Z) \arrow[r, "\rho'"] & H^0(Z', L_Z^{\otimes 2}|_{Z'}) \arrow[r, "\simeq"] & H^0(W, L_W^{\otimes 2}). 
\end{tikzcd}
\]
Note that $S^2\rho$ is an isomorphism by Claim \ref{cl-sym2}. 
It follows from $\varphi_*\MO_{Z'} = \MO_W$ that both the right horizontal arrows are isomorphisms. 
Finally, $\mu_W$ is injective by  $W=\P^1_k$. 
By diagram chase, $\mu_Z$ is injective. 

It suffices to show  Claim \ref{cl-sym2}. 
We have the following exact sequence: 
\[
0 \to \MO_{Y}(-Z')|_C \to \MO_Z \to \MO_{Z'} \to 0, 
\]
which induces another exact sequence: 
\[
H^0(C, (L \otimes \MO_{Y}(-Z'))|_C) \to H^0(Z, L_Z) \to H^0(Z', L_{Z}|_{Z'}) 
\to H^1(C, (L \otimes \MO_{Y}(-Z'))|_C). 
\]
In order to prove Claim \ref{cl-sym2}, it is enough to show 
$H^0(C, (L \otimes \MO_{Y}(-Z'))|_C) = H^1(C, (L \otimes \MO_{Y}(-Z'))|_C) = 0$. 
It holds that 
\[
-Z' \cdot C = (-Z +C) \cdot C = (-(K_Y+Z) +(K_Y+C)) \cdot C =(K_Y+C)\cdot C = \deg \omega_{C}.  
\]
Thus $H^1(C, (L \otimes \MO_{Y}(-Z'))|_C)=0$ by $L \cdot C=1>0$ and Serre duality. 
The Riemannn--Roch theorem and Serre duality imply  
\[
\deg \omega_C = \chi(C, \omega_C) -\chi(C, \MO_C) = 
2h^1(C, \MO_C) -2h^0(C, \MO_C) = -2h^0(C, \MO_C) \leq -2. 
\]
By $L \cdot C = 1$, we obtain 
\[
\deg ((L \otimes \MO_{Y}(-Z'))|_C) = L \cdot C  + (-Z') \cdot C = 1 + \deg \omega_C \leq -1, 
\]
which implies $H^0(C, (L \otimes \MO_{Y}(-Z'))|_C)=0$, as required.
\qedhere


\end{proof}

\begin{thm}\label{t-ell-12}
We use Notation \ref{n-ell-cases}. 
Assume $-Z^2 =1$ or $-Z^2 =2$. 
Then the following hold. 
\begin{enumerate}
    \item $\dim (\m_P/\m_P^2)=3$. 
\item $\mult\,\MO_{X, P} =2$. 
\end{enumerate}
\end{thm}

\begin{proof}
Let us show (1). 
We only treat the case when $-Z^2 =2$, as both  proofs are very similar. 
We have 
\[
\bigoplus_{n=0}^{\infty} I_n/I_{n+1} \simeq 
\bigoplus_{n=0}^{\infty} H^0(Z, L_Z^{\otimes n}). 
\]
Fix a $k$-linear basis $\overline x, \overline y$ of  $I_1/I_2$: 
\[
H^0(Z, L_Z) \simeq I_1/I_2 = k \overline{x} \oplus k \overline{y}.
\]
By Lemma \ref{l-sym2}, we obtain   
\[
H^0(Z, L_Z^{\otimes 2}) \simeq I_2/I_3 = k \overline{x}^2 \oplus k \overline{x}\overline{y} \oplus k\overline{y}^2 \oplus k \overline z
\]
for some $\overline z \in I_2/I_3$. 
Since $\bigoplus_{n=0}^{\infty} H^0(Z, L_Z^{\otimes n})$ is generated by $H^0(Z, L_Z) \oplus H^0(Z, {\cred L_Z^{\otimes 2}})$ as a graded $k$-algebra (Lemma \ref{l-deg12-generation}{\cred )}, 
the following hold. 
\begin{itemize}
\item $H^0(Z, L_Z) \simeq I_1/I_2 = k\overline{x}\oplus k\overline y$. 
\item $H^0(Z, L_Z^{\otimes 2}) \simeq I_2/I_3 = k\overline{x}^2\oplus k\overline{xy} \oplus k\overline y \oplus k\overline{z}$. 
\item $H^0(Z, L_Z^{\otimes 3}) \simeq I_3/I_4 = 
k\overline{x}^3 + k\overline{x}^2\overline{y} + k\overline{x}\overline{y}^2 
+ k\overline y^3 + k\overline{x}\overline{z} + k\overline{y}\overline{z}$. 
\item $H^0(Z, L_Z^{\otimes 4}) \simeq I_4/I_5 = 
k\overline{x}^4 + k\overline{x}^3\overline{y} + k\overline{x}^2\overline{y}^2 
+ k\overline{x}\overline{y}^3 
+ k\overline y^4 + k\overline{x}^2\overline{z} + k\overline{xy}\overline{z} + k\overline{y}^2\overline{z} +  \overline{z}^2$. 
\end{itemize}
Take lifts  $x, y \in I_1, z \in I_2$ of $\overline x, \overline y, \overline z$. 
Set 
\[
A := \Gamma(X, \MO_X),\qquad \m := \m_P \subset A, 
\qquad \m' := Ax +Ay+Az \subset \m \subset A, 
\]
\[
J_1 := Ax + Ay, \qquad J_2 := Az.
\]
We then obtain 
\begin{itemize}
\item $I_1 = I_2+ J_1$. 
\item $I_2 =  I_3+J_1^2 + J_2$. 
\item $I_3 = I_4+J_1^3 + J_1J_2$. 
\item $I_4 = I_5 + J_1^4 + J_1^2J_2 + J_2^2$. 
\end{itemize}
Similarly, the following hold for every $n \in \Z_{>0}$: 
\[
I_{2n} = I_{2n+1} + \sum_{\ell=0}^n J_1^{2n-2\ell}J_2^{\ell}, 
\qquad 
I_{2n-1} = I_{2n} +  \sum_{\ell=0}^{n-1} J_1^{2n-1-2\ell}J_2^{\ell}. 
\]
Hence 
\[
I_{2n} \subset I_{2n+1} + J_1 + J_2 \subset I_{2n+1} + \m',
\qquad 
I_{2n-1} \subset I_{2n} + J_1 + J_2 \subset I_{2n} + \m'. 
\]
By applying these inclusions repeatedly, we obtain 
\[
I_1 \subset I_2 +\m' \subset I_3 + \m' \subset \cdots \subset I_N + \m'
\]
for every $N \in \Z_{>0}$. 
It follows from Lemma \ref{l-Artin-Rees} that $I_N = I_a^2 I_{N-2a} \subset I_a^2 \subset I_1^2$ for some $a \in \Z_{>0}$ and $N \gg 0$. 
By $\m =I_1$ (Lemma \ref{l-I1=m}), we get 
\[
\m = I_1 \subset I_N + \m' \subset I_1^2 + \m' \subset \m^2 + \m'. 
\]
Therefore, the maximal ideal $\mathfrak n := \m \cdot A_{\m}/\m'A_{\m}$ of the local ring $A_{\m}/\m'A_{\m}$ satisfies $\mathfrak n^2 = \mathfrak n$. 
Then Nakayama's lemma implies $\mathfrak n =0$, i.e., $\m A_{\m} = \m'A_{\m}$. 
Hence $\m A_{\m}$ is generated by at most three elements, i.e., 
$\dim (\m/\m^2) \leq 3$. 
On the other hand, $A_{\m}$ is not regular, and hence $\dim (\m/\m^2)>2$. 
Thus (1) holds.

Let us show (2). 
If $-Z^2 = 2$, then $\MO_Y(-Z)$ is globally generated (Theorem \ref{t-bpf}), and hence 
we get $\mult\,\MO_{X, P} = -Z^2 = 2$ by Proposition \ref{p-fund-mult}(3)(5). 
Assume $-Z^2 =1$. 
Note that the base scheme $\Bs\,|\MO_Y(-Z)|$ is equal to a $k$-rational point $Q$ (Proposition \ref{p-bpf-1}). 
For the time being, we finish the proof by assuming the following Claim. 

\begin{claim}\label{cl-ell-12}
$\m_Q\MO_Y(-Z_f) \subset \m_P \MO_Y$. 
\end{claim}

Claim \ref{cl-ell-12} and Proposition \ref{p-fund-mult}(1) imply 
\[
\m_Q\MO_Y(-Z) \subset \m_P \MO_{Y} \subset \MO_Y(-Z). 
\]
Then either $\m_P \MO_{Y} = \MO_Y(-Z)$ or 
$\m_P \MO_{Y} =\m_Q\MO_Y(-Z)$ holds. 
If the former one  $\m_P \MO_{Y} = \MO_Y(-Z)$ is true, 
then  Proposition \ref{p-fund-mult}(3) would imply $\mult\,\MO_{X, P} = -Z^2 =1$, 
which is a contradiction 
(since $X$ is Cohen-Macaulay, 
$\mult\,\MO_{X, P} =1$ implies that $\MO_{X, P}$ is a regular local ring \cite[Theorem 40.6]{Nag62}). 
Therefore, we get $\m_P \MO_{Y} =\m_Q\MO_Y(-Z)$. 
Let $\sigma: Y' \to Y$ be the blowup at $Q$. 
For $E := \Ex(\sigma)$, 
we obtain 
\[
\m_P\MO_{Y'} = \m_Q\MO_Y(-Z) \cdot \MO_{Y'} = \MO_{Y'}(-\sigma^*Z-E). 
\]
By Proposition \ref{p-fund-mult}(3), we get 
\[
\mult\,\MO_{X, P} = -(\sigma^*Z+E)^2 = -(Z^2 +E^2) =2. 
\]

It is enough to show Claim \ref{cl-ell-12}. 
By the scheme-theoretic equality $\Bs\,|\MO_Y(-Z)| = Q$, 
there exists a surjective $\MO_Y$-module homomorphism 
\[
\MO_Y^{\oplus N} \to \m_Q \MO_Y(-Z), \qquad (a_1, ..., a_N) \mapsto a_1 \zeta_1+ \cdots + a_N\zeta_N,  
\]
\[
\text{where}\qquad \zeta_1, ..., \zeta_N \in H^0(Y, \m_Q \MO_Y(-Z)). 
\]
Via the identification $f_*\MO_Y = \MO_X$, 
we get $\zeta_1, ..., \zeta_N \in H^0(Y, \m_Q \MO_Y(-Z)) \subset \m_P$. 
Therefore, we get $\m_Q \MO_Y(-Z) \subset \m_P \MO_Y$. 
\qedhere 



\end{proof}


\subsection{Miscellanies}

\begin{lem}\label{l-ell-omega-char}
We work over a field $k$. 
Let $X$ be a Gorenstein normal affine surface with a unique non-regular point $P$. 
Let $f: Y \to X$ be the minimal resolution of $X$. 
Let $D$ be the effective Cartier divisor defined by $K_Y + D \sim f^*K_X$. 
Assume that $D \neq 0$. 
Then the following hold. 
\begin{enumerate}
\item 
$H^1(Y, \MO_Y) \xrightarrow{\simeq} H^1(D, \MO_D)$. 
\item $\omega_X/f_*\omega_Y \simeq H^0(D, \omega_D)$. 
\item $\dim_k R^1f_*\MO_Y = \dim_k (\omega_X/f_*\omega_Y)$. 
\end{enumerate}
\end{lem}

\begin{proof}
Let us show (1). 
Consider the exact sequence 
\[
0 \to \MO_Y(-D) \to \MO_Y \to \MO_D \to 0. 
\]
Since $-D-K_Y$ is $f$-nef, we obtain $R^if_*\MO_Y(-D)=0$ for every $i>0$ \cite[Theorem 10.4]{Kol13}. 
Thus (1) holds. 

Let us show (2). 
We have an exact sequence 
\[
0 \to \omega_Y \to \omega_Y(D) \to \omega_D \to 0, \qquad \text{where}\qquad \omega_Y(D) := \MO_Y(K_Y+D). 
\]
By $R^1f_*\omega_Y=0$ \cite[Theorem 10.4]{Kol13}, we obtain 
\[
0 \to f_*\omega_Y \to f_*\omega_Y(D) \to f_*\omega_D (\simeq H^0(D, \omega_D)) \to 0. 
\]
We have $f_*\omega_Y(D) = f_*\MO_Y(K_Y+D) \simeq f_*f^*\MO_X(K_X) \simeq \omega_X$. 
Thus (2) holds. 
The assertion (3) follows from 
\[
\dim_k (\omega_X/f_*\omega_Y) \overset{{\rm (2)}}{=} h^0(D, \omega_D) =h^1(D, \MO_D) 
\overset{{\rm (1)}}{=} h^1(Y, \MO_Y) = \dim_k R^1f_*\MO_Y. 
\]
\end{proof}

\begin{prop}\label{p-ell-omega-char}
We work over a field $k$. 
Let $X$ be a Gorenstein normal surface with a unique non-regular point $P$. 
Let $f : Y \to X$ be a projective birational morphism from a regular surface. 
Then the following hold. 
\begin{enumerate}
\item 
The following are equivalent. 
\begin{enumerate}
\item $P$ is a canonical singularity of $X$. 
\item $f_*\omega_Y = \omega_X$
\end{enumerate}
\item 
The following are equivalent. 
\begin{enumerate}
\item $P$ is an elliptic singularity of $X$. 
\item $f_*\omega_Y = \m_P \cdot \omega_X$
\end{enumerate}
\end{enumerate}
\end{prop}

\begin{proof}
Taking a suitable compactification after taking an affine open neighbourhood of $P \in X$, 
we may assume that $X$ is projective over $k$. 
The assertion (1) follows from the same argument as in \cite[Theorem 5.22 and Corollary 5.24]{KM98}. 
The assertion (2) holds by Lemma \ref{l-ell-omega-char}(3). 
\end{proof}

\section{Results for threefolds}

{\cred Based} on 
results in Section \ref{s-ell}, 
we prove that an arbitrary flop of a smooth threefold is again smooth (Theorem \ref{t-flop} in Subsection \ref{ss-flop}). 
To this end, we shall study three-dimensional Gorenstein terminal singulairities 
(Theorem \ref{t-terminal3} in Subsection \ref{ss-term3}). 
In Subsection \ref{ss-generic} and Subsection \ref{ss-vanishing}, 
we shall summarise some results which will be needed in Subsection \ref{ss-term3} and Subsection \ref{ss-flop}, respectively. 
The strategy in this section is almost identical to that of \cite[Section 6.2]{KM98}.

\subsection{Generic members}\label{ss-generic}

The purpose of this section is to recall some terminologies and results on generic members \cite{Tana}. 

Let $k$ be a field and let $X$ be a variety over a field $k$. 
Take a Cartier divisor $D$ on $X$ and a finite-dimensional nonzero $k$-vector subspace $V \subset H^0(X, \MO_X(D))$. 
For its function field $\kappa := K(\P(V))$, 
we obtain the following commutative diagram in which all the squares are cartesian: 
\[
\begin{tikzcd}
	X^{\gen}_{L, V}  & X^{\univ}_{L, V} \\
	X \times_k \kappa & X \times_k \P(V) & X\\
	\Spec\,\kappa & \P(V) & \Spec\,k.
	\arrow[from=3-2, to=3-3]
	\arrow[from=3-1, to=3-2]
	\arrow[from=2-3, to=3-3]
	\arrow[from=2-2, to=3-2, "{\rm pr}_2"']
	\arrow[from=2-1, to=3-1]
	\arrow[from=2-1, to=2-2]
	\arrow[hook, from=1-2, to=2-2, "j"']
	\arrow[hook, from=1-1, to=2-1]
	\arrow[from=1-1, to=1-2, "\gamma"]
	\arrow[from=1-2, to=2-3, "\alpha"]
	\arrow[from=2-2, to=2-3, "{\rm pr}_1"']
 \arrow[from=1-1, to=2-3, "\beta", bend left=15mm]
\end{tikzcd}
\]
where $X^{\univ}_{L, V}$ denotes the universal family. 
In this case, 
\begin{itemize}
    \item $\kappa$ is called the {\em function field} of $\Lambda$ if $\Lambda$ denotes the linear system corresponding to $V$. 
    \item $X^{\gen}_{L, V}$ is called the {\em generic member  of} $V$ (or $\Lambda$). 
    \item If $\Lambda$ is the complete linear system of $D$, then $X^{\gen}_{L, V}$ is called the {\em generic member of} $D$. 
\end{itemize}

\begin{thm}\label{t-generic}
We use the same notation as above. 
If $X$ is regular (resp. normal), then so is $X^{\gen}_{L, V} \setminus \beta^{-1}(\Bs\,\Lambda)$.  
\end{thm}

\begin{proof}
By \cite[Proposition 5.10]{Tana}, we may replace $X$ 
by 
$X \setminus \Bs\,\Lambda$. 
Hence the problem is reduced to the case when $\Lambda$ is base point free. 
Then the assertion follows from \cite[Theorem 4.9 and Remark 5.8]{Tana}. 
\end{proof}



\subsection{Gorenstein terminal singularities}\label{ss-term3}


\begin{lem}\label{l-push-bup-ideal}
We work over a field $k$. 
Let $f: V \to W$ be a projective birational morphism 
between normal varieties. 
Let $P \in W$ be a $k$-rational point. 
Assume that $\m_P\MO_V$ is an invertible sheaf, where 
$\m_P$ denotes the coherent ideal sheaf on $W$ corresponding to $P$. 
Then $f_*(\m_P\MO_V) = \m_P$. 
\end{lem}

\begin{proof}
The problem is reduced to the case when $W$ is affine. 
Fix a closed embedding $W \subset \A^N_k$ such that $P$ is the origin of $\A^N_k$. 
By the universal property of blowups, 
we may assume that 
$f: V \to W$ can be written as 
\[
f : V \xrightarrow{\nu} V_0 \xrightarrow{f_0} W, 
\]
where $f_0$ is the blowup at $P$ and $\nu$ is the normalisation of $V_0$. 
In particular, $f$ is an isomorphism over $W \setminus P$. 
By $f_*(\m_P\MO_V) \subsetneq f_*\MO_V = \MO_W$, we obtain 
$f_*(\m_P\MO_V) \subset \m_P$. 

It suffices to show $f_*(\m_P\MO_V) \supset \m_P$. 
We can write $\m_P\MO_V=\MO_V(-F)$ for some effective Cartier divisor $F$ on $V$. 
Fix an  affine coordinate: $\A^N_k = \Spec\,k[x_1, ..., x_N]$. 
Then $\m_P = (\overline{x}_1, ..., \overline{x}_N)$, where 
each $\overline{x}_i$ denotes the image of $x_i$ to $\Gamma(X, \MO_X)$. 
By construction, we have the following scheme-theoretic equality 
\[
F = f^{-1}(Z(\overline{x}_1)) \cap \cdots \cap f^{-1}(Z(\overline{x}_N)),
\]
where each $Z(\overline{x}_i)$ is the closed subscheme on $W$ defined by $\overline{x}_i$. 
This implies the inequality $f^{-1}(Z(\overline{x}_i)) \geq F$ of effective Cartier divisors. In other words, we get $\overline{x}_i \in H^0(V, \MO_V(-F))$ via the identification: $\MO_W = f_*\MO_V = H^0(V, \MO_V)$. 
Hence we obtain $\m_P \subset f_*\MO_V(-F) = f_*(\m_P\MO_V)$. 
\end{proof}

\begin{prop}\label{p-cano-or-ell}
Let $k$ be an algebraically closed field and 
let $X$ be a Gorenstein canonical affine threefold over $k$. 
Fix a closed point $P \in X$ and a closed embedding $X \subset \A^N_k$. 
Let $\Lambda$ be the linear system on $X$ 
consisting of all the hyperplane sections passing through $P$.  
Then the generic member $H^{\gen}$ of $\Lambda$ is a normal surface 
over the function field $\kappa$ of $\Lambda$ such that 
{\cred either} 
\begin{enumerate}
    \item $H^{\gen}$ has at worst canonical singularities, or 
    \item $H^{\gen}$ is not canonical and has at worst elliptic singularities. 
\end{enumerate}
\end{prop}


\begin{proof}
By abuse of notation, we set $P := P \times_k \kappa$. 
Note that $H^{\gen} \setminus P$ is a normal affine surface over $\kappa$ which is an 
effective Cartier divisor on $X \times_k \kappa$ (Theorem \ref{t-generic}). 
In particular, $H^{\gen}$ is Gorenstein and regular in codimension one. 
Then $H^{\gen}$ is normal by Serre's criterion. 
Thus $H^{\gen}$ is a Gorenstein normal affine surface over $\kappa$.

For a $k$-scheme $Z$, set $Z_{\kappa} := Z \times_k \kappa$. 
Fix a hyperplane section $H$ on $X$ {\cred passing through $P$}, which is a Cartier divisor satisfying $H_{\kappa} \sim H^{\gen}$. 
Let $f_0 : Y_0 \to X$ be the blowup at $P$ and let $Y \to Y_0$ be a resolution of singularities: 
\[
f : Y \to  Y_0 \xrightarrow{f_0} X, \qquad 
f_{\kappa} : Y_{\kappa}\to 
(Y_0)_{\kappa} \xrightarrow{(f_0)_{\kappa}} X_{\kappa}. 
\]
Since the base scheme $\Bs\,\Lambda$ of $\Lambda$ is scheme-theoretically equal to $P$, $f_0$ coincides with the resolution of the indeterminacies of 
the induced rational map $\varphi_{\Lambda} : X \dashrightarrow \P_k^{N-1}$. 
Hence we obtain 
\[
f^*H \sim \widetilde{H} + F, 
\]
where $F$ is the fixed part of $f^*\Lambda$ and $|\widetilde H|$ is base point free. 
Note that we have $\m_P \MO_Y = \MO_Y(-F)$. 
For the generic members $H^{\gen}$ and $\widetilde H^{\gen}$ of $H$ and $\widetilde H$, 
we obtain 
\[
f_{\kappa}^*H^{\gen} \sim \widetilde{H}^{\gen} + F_{\kappa}. 
\]
Since $Y$ is regular, also $\widetilde{H}^{\gen}$ is regular (Theorem \ref{t-generic}). 
As $X$ is Gorenstein and canonical, 
we have $K_Y= f^*K_X + E$ for some $f$-exceptional effective Cartier divisor $E$ on $Y$. 
For the induced morphism $g : \widetilde{H}^{\gen} \to H^{\gen}$, the following holds:   
\begin{eqnarray*}
\omega_{\widetilde{H}^{\gen}} 
&=& \MO_{Y_{\kappa}}(K_{Y_{\kappa}} +\widetilde{H}^{\gen})|_{\widetilde{H}^{\gen}}\\
&=& \MO_{Y_{\kappa}}(f^*_{\kappa}K_{X_{\kappa}} +E_{\kappa} + f_{\kappa}^*H^{\gen} - F_{\kappa})|_{\widetilde{H}^{\gen}}\\
&=& \MO_{\widetilde{H}^{\gen}}(g^*K_{H^{\gen}} +E_{\kappa}|_{\widetilde{H}^{\gen}}- F_{\kappa}|_{\widetilde{H}^{\gen}}). 
\end{eqnarray*}
By applying $g_*$, we obtain 
\[
g_*\omega_{\widetilde{H}^{\gen}}  = 
g_*(\MO_{\widetilde{H}^{\gen}}(g^*K_{H^{\gen}} +E_{\kappa}|_{\widetilde{H}^{\gen}}- F_{\kappa}|_{\widetilde{H}^{\gen}}) 
\supset 
\]
\[
g_*(\MO_{\widetilde{H}^{\gen}}(g^*K_{H^{\gen}} - F_{\kappa})|_{\widetilde{H}^{\gen}}) 
=\omega_{H^{\gen}} \otimes g_*\MO_{\widetilde{H}^{\gen}}(-F_{\kappa}|_{\widetilde{H}^{\gen}}) =\m_P \cdot \omega_{H^{\gen}}, 
\]
where the last equality holds by Lemma \ref{l-push-bup-ideal}. 
Then we are done by Proposition \ref{p-ell-omega-char}. 
\end{proof}

\begin{thm}\label{t-terminal3}
Let $X$ be a Gorenstein terminal threefold over an algebraically closed field $k$. 
Let $P \in X$ be a singular point and let $\m_P$ be the coherent ideal sheaf on $X$ 
corresponding to ${\cred P}$. 
Then $\dim_k \m_P/\m_P^2 = 4$ and $\mult\,\MO_{X, P} =2$. 
\end{thm}

\begin{proof}
We may assume that $X$ is affine and $K_X \sim 0$. 
Fix a closed embedding $X \subset \A^N_k$. 
Let $\Lambda$ be  the linear system consisting of all the hyperplane sections $H \subset X$ passing through $P$. 
Let $H^{\gen}$ be the generic member of $\Lambda$. 
By Proposition \ref{p-cano-or-ell}, the following hold. 
\begin{itemize}
\item $H^{\gen}$ is a normal prime divisor on $X \times_k \kappa$ for the function field $\kappa$ of $\Lambda$
\item $P \in H^{\gen}$  is either a canonical singularity or an elliptic singularity, 
where $P$ denotes $P \times_k \kappa$ by abuse of notation. 
\end{itemize}
There is nothing to show for the former case, i.e.,  $P \in H^{\gen}$  is  a canonical singularity, because 
both  assertions follow from $\dim \m_{H^{\gen}, P}/\m_{H^{\gen}, P}^2 =3$ and $\mult\,\MO_{H^{\gen}, P} =2$. 
Hence we may assume that $P \in H^{\gen}$ is an elliptic singularity. 
Let $Z$ be the fundamental cycle of $P \in H^{\gen}$ (i.e., 
of its minimal resolution of $P \in H^{\gen}$). 
Recall that $-Z^2 \in \Z_{>0}$. 
If $-Z^2=1$ or $-Z^2=2$, then we get 
 $\dim \m_{H^{{\cred \gen}}, P}/\m_{H^{{\cred \gen}}, P}^2 =3$ and $\mult\,\MO_{H^{{\cred \gen}}, P} =2$  (Theorem \ref{t-ell-12}), which imply what we want. 

 Hence the problem is reduced to the case when 
 $-Z^2 \geq 3$. 
We shall prove that this case does not occur 
by deriving a contradiction. 
In what follows, we set $V_{\kappa} := V \times_k \kappa$ when $V$ is a $k$-scheme. 
Let $f_0 : Y_0 \to X$ be the blowup at $P$. 
Let $\nu : Y \to Y_0$ be the normalisation of $Y_0$: 
\[
f : Y \xrightarrow{\nu} Y_0 \xrightarrow{f_0} X, \qquad 
f_{\kappa} : Y_{\kappa} \xrightarrow{\nu_{\kappa}} Y_{0, \kappa} \xrightarrow{f_{0, \kappa}} X_{\kappa}.
\]
Then 
$f_0$ coincides with the resolution of the indeterminacies of the 
rational map $\varphi_{\Lambda}$ induced by $\Lambda$. 
For a hyperplane section $H$ on $X$ and the effective Cartier divisor $F_0$ satisfying $\m_P \MO_{Y_0} = \MO_{Y_0}(-F_0)$, we obtain 
\[
f_0^*H \sim H_{Y_0} + F_0, 
\]
where $H_{Y_0}$ is a Cartier divisor such that $|H_{Y_0}|$ is base point free. 
For the generic members  $H^{\gen}$ and $H^{\gen}_{Y_0}$  of $H$ and $\widetilde H$, 
we obtain 
\[
f_{0, \kappa}^*H^{\gen} \sim H_{Y_0}^{\gen} + F_{0, \kappa}. 
\]
Note that both $H_{Y_0}^{\gen}$ and $F_{0, \kappa}$ are effective Cartier divisors on $Y_{0, \kappa}$. 
Since the induced morphism $g: H^{\gen}_{Y_0} \to H^{\gen}$ is the blowup at $P$, 
$g$ (and $H^{\gen}_{Y_0}$) coincides with the canonical model over $H^{\gen}$ (Theorem \ref{t-ell-3}). 
In particular, $H^{\gen}_{Y_0}$ is normal. 
Since $H^{\gen}_{Y_0}$ is a normal effective Cartier divisor on $Y_{0, \kappa}$, 
$Y_{0, \kappa}$ is normal around $H^{\gen}_{Y_0}$, i.e., 
$\nu_{\kappa}: Y_{\kappa} \to Y_{0, \kappa}$ is an isomorphism around $H^{\gen}_{Y_0}$. 
Set $H^{\gen}_Y := \nu_{\kappa}^*H^{\gen}_{Y_0}$ and $F := \nu^*F_0$, so that we obtain 
\[
f_{\kappa}^* H^{\gen} = H^{\gen}_{Y} + F_{\kappa}. 
\]
Since $X$ is Gorenstein and terminal, we have $K_Y = f^*K_X +E$ for some effective $\Z$-divisor $E$. 
We then get 
\[
-F_{\kappa}|_{H^{\gen}_Y} \overset{{\rm (i)}}{\sim} 
K_{H^{\gen}_Y} \overset{{\rm (ii)}}{\sim}  (K_{Y_{\kappa}}+H^{\gen}_Y)|_{H^{\gen}_Y} 
\]
\[
= (f_{\kappa}^*K_{X_{\kappa}}+E_{\kappa} +f_{\kappa}^*H^{\gen} -F_{\kappa})|_{H^{\gen}_Y} 
\sim (E_{\kappa}-F_{\kappa})|_{H^{\gen}_Y}, 
\]
where (i) and (ii) will be proven below. 
Note that $E_{\kappa}$ is Cartier around $H^{\gen}_Y$ (cf. the proof of (ii) below).

(i) 
We identify $H^{\gen}_{Y_0}$ and $H^{\gen}_{Y}$. 
Since  $g : H^{\gen}_{Y_0} (\simeq H^{\gen}_{Y}) \to H^{\gen}$ is the canonical model over $H^{\gen}$ for an elliptic singularity $P \in H^{\gen}$, 
we have $K_{H^{\gen}_Y} + Z_g \sim g^*K_{H^{\gen}}$ (Lemma \ref{l-ell-fund}(2)) 
and $\MO_{H^{\gen}_Y}(-Z_g) = \m_P \MO_{H^{\gen}_Y} = \MO_Y(-F)|_{H^{\gen}_Y}$ (Corollary \ref{c-ell-3}), where $Z_g$ denotes the push-forward of the fundamental cycle $Z$ of $P \in H^{\gen}$. 

(ii) 
Note that $Y_{\kappa}$ is Gorenstein around $H^{\gen}_Y$, 
because $Y_{\kappa}$ is a variety, $H^{\gen}_Y$ is Gorenstein, 
and $H^{\gen}_Y$ is an effective Cartier divisor on $Y_{\kappa}$. 
Therefore, we obtain $K_{H^{\gen}_Y} =  (K_{Y_{\kappa}}+H^{\gen}_Y)|_{H^{\gen}_Y}$. 

\medskip

Hence $E_{\kappa}|_{H^{\gen}_Y} \sim 0$. 
Note that either $E_{\kappa}=0$ or $f_{\kappa}(E_{\kappa}) = P$ as sets. 
In any case, $E_{\kappa}|_{H^{\gen}_Y}$ is $g$-exceptional. 
By the negativity lemma, we obtain the equality $E_{\kappa}|_{H^{\gen}_Y}=0$ of divisors. 
As  $H^{\gen}_{Y_0}$ is $f_{0, \kappa}$-ample, $H^{\gen}_Y=\nu_{\kappa}^*H^{\gen}_{Y_0}$ is $f_{\kappa}$-ample. 
Then $H_Y^{\gen}$ intersects all the irreducible components of $E_{\kappa}$, 
and hence $E_{\kappa} = 0$. 
Therefore, $K_{Y_{\kappa}} = f_{\kappa}^*K_{X_{\kappa}}$, which implies 
$K_{Y} = f^*K_{X}$. 
Then $a(E', X, 0)=0$ for 
the discrepancy $a(E', X, 0)$ of a prime divisor $E'$ contained in $E$. 
Therefore, $X$ is not terminal. This is absurd. 
\qedhere



\end{proof}



\subsection{Vanishing theorems}\label{ss-vanishing}

The purpose of this subsection is to establish a vanishing theorem for flopping contractions (Corollary \ref{c-reg-van}). 
Although the argument is very similar to those of \cite{Tan15} and \cite{Tan18m}, we include the proofs for the sake of completeness. 
We start by recalling $F$-finiteness and strongly $F$-regular pairs for the reader's convenience. 

\begin{nasi}[$F$-finiteness] \label{n-Ffinite}
We say that an $\F_p$-scheme $X$ is {\em $F$-finite} if 
the absolute Frobenius morphism 
$F : X \to X$ is finite. 
We say that an $\F_p$-algebra $R$ is {\em $F$-finite} if $\Spec R$ is $F$-finite. 
In particular, a field of characteristic $p>0$ is $F$-finite field if and only if $[k : k^p] <\infty$. 
It is easy to see that 
if $X$ is a scheme of finite type over an $F$-finite noetherian $\F_p$-scheme $Y$, then also $X$ is $F$-finite. 
\end{nasi}

\begin{nasi}[Strongly $F$-regular pairs]\label{n-SFR}
We work over an  $F$-finite field $k$. 
We say that   $(X, \Delta)$ is  {\em strongly $F$-regular} 
if $X$ is a normal quasi-projective variety (over $k$) and $\Delta$ is an effective $\Q$-divisor on $X$ such that
\begin{itemize}
    \item $K_X+\Delta$ is $\Q$-Cartier, and 
    \item given a point $x \in X$ and an effective $\Q$-divisor $E$ on $X$, 
there exists $e \in \Z_{>0}$ such that 
\[
\MO_{X, x} \xrightarrow{F^e} F_*^e\MO_{X, x} \hookrightarrow F^e(\MO_X( (p^e-1)\Delta +E))_x
\]
splits as an $\MO_{X, x}$-module homomorphism. 
\end{itemize}
{\cred Although it is not so common to impose the condition that $K_X+\Delta$ is $\Q$-Cartier, we assume it for simplicity.} 
It is known that the following properties hold. 
\begin{enumerate}
\item If $X$ is a regular variety, then $(X, 0)$ is strongly $F$-regular. 
\item If $(X, \Delta)$ is strongly $F$-regular 
and $\Delta'$ is a {\cred $\Q$-Cartier} $\Q$-divisor 
such that $0 \leq \Delta' \leq \Delta$ {\cred and $K_X+\Delta'$ is $\Q$-Cartier}, 
then $(X, \Delta')$ is strongly $F$-regular. 
\item Let $(X, \Delta)$ be a strongly $F$-regular pair and let $F$ be an effective {\cred $\Q$-Cartier} $\Q$-divisor on $F$. 
Then there exists a rational number $\epsilon >0$ such that $(X, \Delta + \epsilon F)$ is strongly $F$-regular \cite[Theorem 3.9]{SS10}. 
\item If $(X, \Delta)$ is strongly $F$-regular, then 
there is a surjective $\MO_X$-module homomorphism (cf. \cite[Proposition 2.6]{CTX15}): 
\[
\Tr^e_X : F_*^e\MO_X( \ulcorner -(p^e-1)(K_X+\Delta)\urcorner ) \to \MO_X.  
\]
\end{enumerate}
For more foundational results, we refer to \cite{Sch08}, \cite{SS10}. 
\end{nasi}

\begin{thm}\label{t-Ffin-van}
We work over an $F$-finite  field $k$ of characteristic $p>0$. 
Let $f: X \to Y$ be a projective morphism from a normal variety $X$ to a quasi-projective $k$-scheme $Y$. 
Assume that $f$ is not a finite morphism. 
Let $L$ be a Cartier divisor on $X$. 
Assume that there exists an effective $\Q$-divisor $\Delta$ on $X$ such that  
\begin{enumerate}
\item $(X, \Delta)$ 
is strongly $F$-regular and 
\item $L-(K_X+\Delta)$ is $f$-nef and $f$-big. 
\end{enumerate} 
Then $R^if_*\MO_X(L)=0$ for every ${\displaystyle i \geq \max_{y \in Y} \dim f^{-1}(y)}$. 
\end{thm}

\begin{proof}
We first reduce the problem to the case when 
\begin{enumerate}
\item[(3)] $L-(K_X+\Delta)$ is $f$-ample. 
\end{enumerate}
Since $L-(K_X+\Delta)$ is $f$-big, 
we can write $L - (K_X+\Delta) = A +E$ for some $\Q$-Cartier $\Q$-divisors $A$ and $E$ such that 
$A$ is $f$-ample and $E$ is effective. 
For 
a rational number $0 < \epsilon < 1$, we get 
\[
L- (K_X + \Delta) -\epsilon E = A+E -\epsilon E = \epsilon A + (1-\epsilon)(A+E), 
\]
which is $f$-ample. 
Since $(X, \Delta + \epsilon E)$ is strongly $F$-regular for any $0 < \epsilon \ll 1$ (\ref{n-SFR}), 
the problem is reduced, after replacing $\Delta$ by $\Delta + \epsilon E$, 
to the case when (1) and (3) hold (note that (3) implies (2)). 

Furthermore, we may assume that 
\begin{enumerate}
\item[(4)] $(p^{e_1}-1)(K_X+\Delta)$ is Cartier for some $e_1 \in \Z_{>0}$. 
\end{enumerate}
Indeed, this condition is achieved by applying the same argument as in 
\cite[Lemma 4.1]{CTX15}. 

\medskip

Fix an integer $i$ satisfying ${\displaystyle i \geq \max_{y \in Y} \dim f^{-1}(y)}$. 
By (3) and the Serre vanishing theorem, 
we can find $e \in e_1 \Z_{>0}$ such that 
\[
R^if_* \MO_X(p^eL- (p^e-1)(K_X+\Delta))=0. 
\]
It follows from (4) and $e \in e_1 \Z_{>0}$ that $(p^e-1)(K_X+\Delta)$ is Cartier. 
By (1), 
we have a short exact sequence  (\ref{n-SFR}): 
\[
0 \to B_e \to F_*^e \MO_X(- (p^e-1)(K_X+\Delta)) \to \MO_X \to 0
\]
for every $e \in e_1\Z_{>0}$, where $B_e$ is a coherent sheaf depending on $e$. 
By ${\displaystyle i \geq \max_{y \in Y} \dim f^{-1}(y)}$, we obtain 
\[
R^if_*F_*^e \MO_X(p^eL- (p^e-1)(K_X+\Delta)) \to R^if_*\MO_X(L) \to R^{i+1}f_*(B_e \otimes \MO_X(L))=0. 
\]
It follows from $R^if_*F_*^e \simeq F_*^eR^if_*$ that 
\[
R^if_*F_*^e \MO_X(p^eL- (p^e-1)(K_X+\Delta)) 
\simeq 
F_*^eR^if_* \MO_X(p^eL- (p^e-1)(K_X+\Delta)) =0. 
\]
Hence $R^if_*\MO_X(L) =0$. 
\end{proof}

\begin{cor}\label{c-reg-van}
We work over a field $k$. 
Let $f: X \to Y$ be a projective morphism from a regular variety $X$ to a quasi-projective $k$-scheme $Y$ such that $f$ is not a finite morphism. 
Assume that $L$ is a Cartier divisor on $X$ such that $L-K_X$ is $f$-nef and $f$-big. 
Then $R^if_*\MO_X(L)=0$ for every ${\displaystyle i \geq \max_{y \in Y} \dim f^{-1}(y)}$. 
\end{cor}

\begin{proof}
If $k$ is of characteristic zero, then the assertion is well known. 
We may assume that $k$ is of characteristic $p>0$. 
Taking a suitable intermediate field $\F_p \subset k' \subset k$ such that $k'$ is finitely generated over $\F_p$, the problem is reduced to the case when $k$ is an $F$-finite field of characteristic $p>0$. 
Then the assertion follows from Theorem \ref{t-Ffin-van}. 
\end{proof}

\subsection{Flops}\label{ss-flop}

\begin{nota}\label{n-flop-exist}
Let $k$ be an algebraically closed field of characteristic $p>0$ and 
let $Y$ be a smooth threefold over $k$. 
Let $\psi : Y \to Z$ be a flopping contraction, 
i.e., 
$\psi$ is a projective birational morphism to a normal threefold $Z$ over $k$ such that 
\begin{enumerate}
\item $\dim \Ex(\psi)=1$, 
\item $\rho(Y/Z)=1$, and 
\item $K_Y$ is $\psi$-numerically trivial.  
\end{enumerate}
\end{nota}

\begin{prop}\label{p-flop-van}
We use Notation \ref{n-flop-exist}. 
Then the following hold. 
\begin{enumerate}
\item $R^i\psi_*\MO_Y =0$ for every $i >0$. 
\item $R^i\psi_*\omega_Y =0$ for every $i >0$. 
\item $K_Z$ is Cartier and $K_Y \sim \psi^*K_Z$. 
\item $Z$ is Gorenstein and terminal. 
\end{enumerate}
\end{prop}

\begin{proof}
The assertions (1) and (2) follow from Corollary \ref{c-reg-van}. 
Let us show (3). 
Since $K_Y$ is $\psi$-nef and $\dim \Ex(\psi)=1$, $K_Y$ is $\psi$-free (Proposition \ref{p-Lipman}(1)). 
Therefore, we can write $K_Y \sim f^*D$ for some Cartier divisor $D$ on $Z$. 
We then get $K_Z \sim f_*K_Y \sim f_*f^*D = D$. Thus (3) holds. 

Let us show (4). 
We first prove that $Z$ is terminal. 
Let $\mu : V \to Z$ be a projective birational morphism and let $E$ be a prime divisor on $V$. 
We can write $K_V \sim \mu^*K_Z+ \sum_E a_E E$ for some $a_E \in \Z$. 
It is enough to prove $a_E >0$. 
Since we may replace $V$ by a higher model,  we may assume that $\mu : V \to Z$ factors through $\psi: Y \to Z$: 
\[
\mu : V \xrightarrow{\mu'} Y \xrightarrow{\psi} Z. 
\]
By $\dim \Ex(\psi)=1$ and $K_Y \sim \psi^*K_Z$, we see that $a_E>0$, as $Y$ is terminal. 
Therefore, $Z$ is terminal.

Since $K_Z$ is Cartier, it is enough to show that $Z$ is Cohen--Macaulay. 
By taking suitable compactifications {\cred and desingularisations} of $Y$ and $Z$ (after taking an affine open cover of $Z$), 
we may assume that $Y$ and $Z$ are projective. 
Fix an ample Cartier divisor $H_Z$ on $Z$. 
Set $H_Y := \psi^*H_Z$. 
By \cite[Ch. III, Theorem 7.6(b)]{Har77}, it is enough to show that 
$H^i(Z, \MO_Z(-qH_Z))=0$ for $i<3$ and $q \gg 0$. 
Recall the following Leray spectral sequence: 
\[
E_2^{i, j} = H^i(Z, R^j\psi_*\MO_Y(-qH_Y)) \Rightarrow H^{i+j}(Y, \MO_Y(-qH_Y)) =E^{i+j}. 
\]
We have  
\[
R^j\psi_*\MO_Y(-qH_Y) 
\simeq R^j \psi_* \psi^*\MO_Z(-qH_Z) 
\simeq (R^j\psi_*\MO_Y) \otimes \MO_Z(-qH_Z) \overset{{\rm (1)}}{=}0
\]
for every $j>0$, which implies 
\[
H^i(Z, \MO_Z(-qH_Z)) \simeq H^i(Z, \psi_*\MO_Y(-qH_Y)) = E_2^{i, 0} 
\simeq E^i =  H^{i}(Y, \MO_Y(-qH_Y)). 
\]
By Serre duality, we get 
\[
h^i(Y, \MO_Y(-qH_Y))  = h^{3-i}(Y, \MO_Y( K_Y +qH_Y)). 
\]
Therefore, it is enough to prove $H^i(Y, \MO_Y(K_Y +qH_Y))=0$ for $i>0$ and $q \gg 0$. 
Again by a Leray spectral sequence 
\[
H^i(Z, R^j\psi_*\MO_Y(K_Y +qH_Y)) \Rightarrow H^{i+j}(Y, \MO_Y(K_Y +qH_Y)), 
\]
the same argument as above implies 
\[
H^i(Z, \MO_Z(K_Z + qH_Z)) \simeq 
H^i(Z, \psi_*\MO_Z(K_Y + qH_Y)) \simeq 
H^{i}(Y, \MO_Y(K_Y +qH_Y)). 
\]
Since $H_Z$ is ample, 
we get $H^i(Z, \MO_Z(K_Z + qH_Z))  =0$ for $i>0$ and $q \gg 0$ by the Serre vanishing theorem. 
Thus (4) holds. 
\end{proof}

\begin{lem}[Weierstrass Preparation Theorem]\label{l-WPT}
Let $(A, \m, \kappa)$ be a noetherian complete local ring. 
Take an element $f = \sum_{n=0}^{\infty}a_nt^n \in A[[t]]$ and 
let $\overline{f} 
\in \kappa [[t]]$ be the image of $f$. 
Assume that $\overline{f} \neq 0$, i.e., 
$a_r \not\in \m$ for some $r \in \Z_{\geq 0}$. 
Let $s$ be the smallest non-negative integer such that $a_s \not\in \m$. 
Then there exist $u \in (A[[t]])^{\times}$ and 
\[
F =t^s +b_{s-1}t^{s-1}+ \cdots + b_0 \in A[t]
\]
such that $b_0, ..., b_{s-1} \in \m$ and $f =uF$. 
In particular, $A[[t]]/(f) = A \oplus At \oplus \cdots \oplus A t^{s-1}$. 
\end{lem}

\begin{proof}
See \cite[Chapter VII, §3, no. 8, Proposition 6]{Bou89}. 
\end{proof}

\begin{lem}\label{l-double-cover}
We use Notation \ref{n-flop-exist}. 
Let $P \in Z$ be a singular point and let $\m_P$ be the maximal ideal of $\widehat{\MO}_{Z, P}$, where $\widehat{\MO}_{Z, P}$ denotes the completion of $\MO_{Z, P}$. 
Then there exist $w \in \m_P \setminus \m^2_P$ and an injective $k$-algebra homomorphism 
\[
A := k[[x, y, z]] \hookrightarrow \widehat{\MO}_{Z, P}
\]
such that 
$\widehat{\MO}_{Z, P} = A \oplus A w$. 
In particular, $\widehat{\MO}_{Z, P}$ is a free $A$-module of rank two. 
\end{lem}

\begin{proof}
By Theorem \ref{t-terminal3} and Proposition \ref{p-flop-van}(4), we can write
\[
\widehat{\MO}_{Z, z} = k[[x, y, z, w]]/(f) 
\qquad \text{ for some }\quad 
f = f_2 + f_3 + \cdots \quad \text{ with } f_2 \neq 0,
\]
where each $f_d \in k[x, y, z, w]$ is a homogeneous polynomial of degree $d$ {\cred (note that 
the equality $\dim_k \m_P/\m_P^2=4$ implies that 
$\m_P=(x, y, z, w)$ 
for some $x, y, z, w$, 
which induces a surjective $k$-algebra homomorphism
$k[[x, y, z, w]] \to \widehat{\MO}_{Z, z}$)}. 
Applying a  suitable coordinate change, 
we may assume that 
\[
f_2 = w^2 + w\ell(x, y, z) + q(x, y, z)
\]
for some homogeneous polynomials $\ell(x, y, z)$ and $q(x, y, z)$ of degree $1$ and $2$, respectively. 
We now apply Lemma \ref{l-WPT} for $A := k[[x, y, z]]$ and $\m := (x, y, z)k[[x, y, z]]$. 
We can uniquely write 
\[
f = a_0 + a_1 w + a_2w^2 + a_3w^3 +\cdots \qquad\text{with}\qquad a_i \in A = k[[x, y, z]]. 
\]
We have $a_0 \in \m,  a_1 \in \m$, and $a_2 \not\in \m$. 
Therefore, we obtain $\widehat{\MO}_{Z, z} = A \oplus A w$ by Lemma \ref{l-WPT}. 
\qedhere


\end{proof}

\begin{thm}\label{t-flop}
We use Notation \ref{n-flop-exist}. 
Then the flop $Y^+$ of $\psi : Y \to Z$ exists. Furthermore, $Y^+$ is smooth. 
\end{thm}

\begin{proof}
Fix a singular point $P \in Z$ {\cred and a $\psi$-ample Cartier divisor $A$ on $Y$}. 
Set $\widetilde Z := \Spec\,\widehat{\MO}_{Z, P}$ and 
$\widetilde Y := Y \times_Z \widetilde Z$. 
By Lemma \ref{l-double-cover}, there exists an injective $k$-algebra homomorphism $k[[x, y, z]] \hookrightarrow \widehat{\MO}_{Z, P}$ such that $\widehat{\MO}_{Z, P}$ is a free 
$k[[x, y, z]]$-module of rank $2$. 
For $W := \Spec\,k[[x, y, z]]$, we obtain the induced morphisms of noetherian integral normal separated schemes: 
\[
\alpha: \widetilde{Y} \xrightarrow{\widetilde \psi} \widetilde{Z} \xrightarrow{\beta} W. 
\]
Then $\beta$ is a {\cred finite} surjective morphism 
with $[K(\widetilde Z) : K(W)]=2$. 
{\cred Set $A_{\widetilde Y}$ to be the pullback of $A$ on $\widetilde{Y}$, 
which is a $\widetilde{\psi}$-ample Cartier divisor.} 
Note that the 
{\cred Weil} divisor  
{\cred $A_{\widetilde Z} := \widetilde{\psi}_*A_{\widetilde Y}$}
on $\widetilde Z$ 
is not $\Q$-Cartier, as otherwise we {\cred would} get 
$A_{\widetilde Y} = \widetilde{\psi}^*A_{\widetilde Z}$,  
which contradicts the $\widetilde \psi$-ampleness of {\cred $A_{\widetilde Y}$}. 
In particular, $\widetilde Z$ is not $\Q$-factorial.  
Then $K(W) \subset K(\widetilde Z)$ is not purely inseparable by \cite[Lemma 2.5]{Tan18b}. 
Hence $K(W) \subset K(\widetilde Z)$ is a Galois extension of degree $2$. 
Since $W$ is regular, 
the 
{\cred Weil} divisor 
$\beta_*A_{\widetilde Z}$ 
is Cartier. 
{\cred It is easy to see that 
\[
\beta^*\beta_* A_{\widetilde Z} = A_{\widetilde Z} + A'_{\widetilde Z}. 
\]
for $A'_{\widetilde Z} := \iota^*A_{\widetilde Z}$ and 
the Galois involution  $\iota: \widetilde Z \xrightarrow{\simeq} \widetilde Z$ of the Galois extension $K(\widetilde Z)/K(W)$ of degree $2$.} 
{\cred 
For the proper transform $A'_{\widetilde Y}$ of 
$A'_{\widetilde Z}$ on $\widetilde Y$, 
 $-A'_{\widetilde Y}$ is $\widetilde \psi$-ample by $A_{\widetilde Y} +A'_{\widetilde Y} = \widetilde \psi^*(A_{\widetilde Z} + A'_{\widetilde Z})$.}

{\cred 
Recall that 
$\widetilde Y \simeq \Proj\, \bigoplus_{d \geq 0} H^0(\widetilde Z, \MO_{\widetilde Z}(d{\cred A_{\widetilde Z}}))$ and 
$\bigoplus_{d \geq 0} H^0(\widetilde Z, \MO_{\widetilde Z}(d{\cred A_{\widetilde Z}}))$ is a finitely generated 
${\Gamma(\widetilde{Z}, \MO_{\widetilde{Z}})}$-algebra. 
Moreover, 
we have an isomorphism 
\[
\bigoplus_{d \geq 0} H^0(\widetilde Z, \MO_{\widetilde Z}(d{\cred A_{\widetilde Z}})) \simeq 
\bigoplus_{d \geq 0} H^0(\widetilde Z, \MO_{\widetilde Z}(d{\cred A'_{\widetilde Z}}))
\]
induced by the Galois involution $\iota : \widetilde{Z} \xrightarrow{\simeq} \widetilde{Z}$. 
Then the flop $\psi^+ : Y^+ \to Z$ of $\psi : Y \to Z$ exists 
 \cite[Lemma 6.2 and Proposition 6.6]{KM98}.  
For $\widetilde{Y}^+ := Y^+ \times_Z \widetilde Z$, 
we get 
\[
\widetilde Y \simeq 
\Proj\, \bigoplus_{d \geq 0} H^0(\widetilde Z, \MO_{\widetilde Z}(d{\cred A_{\widetilde Z}}))
\simeq 
\Proj\, \bigoplus_{d \geq 0} H^0(\widetilde Z, \MO_{\widetilde Z}(d{\cred A'_{\widetilde Z}})) \simeq \widetilde Y^+,  
\]
where the middle $k$-isomorphism is induced by the Galois involution $\iota : \widetilde{Z} \xrightarrow{\simeq} \widetilde{Z}$. 
Since $Y$ is smooth, $\widetilde Y$ is regular. 
Then also $\widetilde{Y}^+$ is regular, and hence $Y^+$ is smooth \cite[Theorem 23.7(i)]{Mat86}.}
\end{proof}

\bibliographystyle{skalpha}
\bibliography{reference.bib}

\end{document}